\documentclass[a4paper,11pt]{article}
\usepackage{amsfonts,amssymb,eucal,amsmath,array,url}
\usepackage{amsfonts,amssymb,eucal,amsmath,array,url,comment,float}

\newcommand{\Cay}{\mbox{Cay}}
\newcommand{\adm}{\'Ad\'am}
\newcommand*\samethanks[1][\value{footnote}]{\footnotemark[#1]}

\newcommand{\Z}{\mathbb{Z}}

\newcommand{\Aut}{\mbox{Aut}}
\newcommand{\Aff}{\mbox{Aff}}
\newcommand{\M}{\mathcal{M}}

\usepackage{setspace}
\usepackage{enumitem}
\usepackage{lscape}
\usepackage{amsmath, amsfonts, amssymb, graphicx, textcomp, amsthm}
\usepackage{tikz}
\usepackage{appendix}
\usepackage{doc} 
\usepackage{index}
\usepackage[tc]{titlepic}
\usepackage{fancyhdr}
\usepackage[nokeyprefix]{refstyle}
\usepackage{listings}
\usepackage{changepage}
\usepackage[left=5cm, bottom=5cm]{geometry}
\usepackage[labelsep=period]{caption}
\numberwithin{equation}{section}

\newcommand\dedicationname{}  
\makeatletter
\if@titlepage
     {\par\vfil\null
}
\fi
\makeatother

\newcommand\Acknowledgementsname{Acknowledgements}  
\makeatletter
\if@titlepage
     {\par\vfil\null
}
\fi
\makeatother

\makeatletter
\g@addto@macro\bfseries{\boldmath}
\makeatother

\newtheorem{Theorem}{Theorem}
\newtheorem{Lemma}{Lemma}
\newtheorem{Proposition}{Proposition}

\newtheorem{Corollary}{Corollary}
\newtheorem{Example}{Example}

\begin{document}









\title{Constructive and analytic enumeration of
circulant graphs with $p^3$ vertices; $p=3,5$}
\author{Victoria Gatt\\ University of Malta, Malta \\ victoria.gatt@um.edu.mt \and Mikhail Klin\thanks{Supported by Project Mobility Grant (see Acknowledgements)} \\ Ben Gurion University of the Negev, Beer Sheva, 84105, Israel\\
Matej Bel University, Bansk{\'a} Bystrica, Slovakia \\ klin@cs.bgu.ac.il  \and Josef Lauri \\ University of Malta, Malta \\ josef.lauri@um.edu.mt \and Valery Liskovets\samethanks \\ National Academy of Sciences of Belarus, Minsk, Belarus \\
liskov@im.bas-net.by}
\date{ }

\maketitle

\begin{abstract}
Two methods, structural (constructive) and multiplier (analytical),
of exact enumeration of undirected and directed circulant graphs of
orders 27 and 125 are elaborated and represented in detail here together with
intermediate and final numerical data. The first method is based on the known
useful classification of circulant graphs in terms of $S$-rings and results
in exhaustive listing (with the use of COCO and GAP) of all corresponding $S$-rings of
the indicated orders. The latter method is conducted in the framework of a
general approach developed earlier for counting circulant graphs of prime-power
orders. It is a Redfield--P\'olya type of enumeration based on an isomorphism criterion for circulant graphs of such orders. In particular, five intermediate enumeration subproblems
arise, which are refined further into eleven subproblems of this type (5 and 11 are,
not accidentally, the 3d Catalan and 3d little Schr\"oder numbers, resp.).
All of them are resolved for the four cases under consideration (again with
the use of GAP). We give a brief survey of some background theory of the results which form the basis of our computational approach.

Except for the case of undirected circulant graphs of orders 27,
the numerical results obtained here are new. In particular  the number 
(up to isomorphism) of directed
circulant graphs of orders 27, regardless of valency, is shown to be equal to 3,728,891 while 457 of these are self-complementary.
Some curious and rather unexpected identities are established between intermediate
valency-specified enumerators (both for undirected and directed circulant graphs) and their
validity is conjectured for arbitrary cubed odd prime $p^3$.
 
We believe that this research can serve as the crucial step towards 
explicit uniform enumeration formulae for circulant graphs of orders $p^3$ for
arbitrary prime $p>2$.
\end{abstract}


\newpage

\tableofcontents

\newpage

\section{Introduction}

The present research is carried out in the framework of the general program
outlined in the paper \cite{KLP96} for counting circulant graphs of prime-power orders.
We refer to this paper for details concerning two approaches to the exact
enumeration of circulant graphs, namely, constructive and analytical. Recall
that the former is based on the known useful classification of circulant
graphs in terms of $S$-rings. This not only counts the nonisomorphic circulant graphs but enables us, in principle at least, actually to list them. 

The analytic approach is based on the familiar isomorphism
theorem \cite{KP80} for circulant graphs of prime-power orders. In analytic enumeration we are
guided also by the subsequent adaptation of this theorem to the enumeration of circulant graphs as
developed in \cite{LP2000}. In particular, for circulant graphs of order $p^3$ their
analytical (that is, formula-wise) enumeration has been reduced in this paper
to five well-specified (and rather sophisticated) enumeration subproblems of
Redfield--P\'olya type. In order to obtain the solutions we refined them further
into more elementary eleven subproblems.

We restrict ourselves to orders 27 and 125 only.
There are several arguments for our choice. First of all, there are almost no
numerical results for the number of isomorphism classes of circulant graphs of prime-cubed orders, including $p=3$ and $5$.
There exist huge numbers of circulant graphs even of these orders so that
it hardly makes much sense to enumerate constructively circulant graphs of larger
orders. In principle, no such obstacle arises for analytical enumeration, but even here these
two orders require much effort. Presumably, the main difficulties of analytical
enumeration should become apparent already on these least prime-cubed orders.
Moreover our aim is also to compare both approaches on the same classes of objects and to obtain
confirmation of numerical results obtained in both ways using COCO and GAP and also partially
by brute force.

As in several enumeration problems, in order to count the number of non-isomorphic structures of a certain type one needs to have a criterion for determining when two such structures are ``the same'', in our case, are isomorphic. The importance and difficulty of counting circulant graphs stems from the falsity of a very natural conjecture of \adm\ giving a condition on the connecting sets for isomorphism to hold, which however turned out to be false: \adm's condition is sufficient for isomorphism but not necessary. The falsity of \adm's conjecture led to some beautiful results which characterised completely the conditions on the order of the circulant graph for \adm's Conjecture to hold. This led to two very important threads of research in algebraic graph theory: discovering what necessary and sufficient conditions on the connecting sets give isomorphism when \adm's Conjecture fails, and the question of when it is possible to determine the isomorphism of general Cayley graphs from conditions on the connecting sets. 

In this paper, since we are enumerating circulant graphs, we shall use results from the first line of research, which we shall describe below.  We shall adopt two very different methods which have been used for enumerating non-isomorphic circulant graphs: the structural and the multiplier approach. We first introduce these approaches in the following subsections. Then, in the next sections, we shall present some results obtained for enumeration of circulant graphs of order $p^k$ mainly for $p=3, 5$ and $k=2,3$ using these methods, giving more detail for the multiplier approach. The results we present for $p=5$ are new as are the results for directed circulants of order $3^3$. Both our numerical results and the generating functions which we obtain are important. In fact we point out several relations which arise between the intermediate terms which form these generating functions. We conjecture that the relations which emerge from our generating functions for $k=3$ and $p=3,5$ hold for all odd prime $p$. In an appendix we give some theoretical support for these conjectures by proving some similar relations for $k=2$ for all odd prime $p$. We then conclude with the enumeration of self-complementary circulant graphs. 

Fuller details including all case-by-case analysis, all generating functions produced and all the GAP programmes used, can be found in \cite{arXivVg}. 

For standard graph theoretic terms we refer the reader to the two texts \cite{klin&al91} and \cite{lauri&sca16}. 

\subsection{First definitions}
A circulant graph is a Cayley graph of a cyclic group. That is, let $G$ be a cyclic group (which we shall represent as the group $\mathbb{Z}_n$ of integers with addition modulo $n$, the size of the group) and let $S\subseteq G$ (called the \emph{connecting set} of the Cayley graph) such that $0\not\in S$. Then a circulant is a Cayley graph $\Cay(G,S)$ which has $G$ as vertex-set and two vertices $g, h$ are adjacent if $g=h+s$ for some $s\in S$. If the set $S$ generates $G$ then the circulant graph $\Cay(G,S)$ is connected. In the special case when  $-S=S$ (that is, $s\in S$ if and only if $-s\in S$), the circulant graph is also referred to as an \emph{undirected graph}. For brevity we shall sometimes refer to ``circulants'' instead of ``circulant graphs''. 

An edge $\{a,b\}$ is considered to be the union of the two arcs $(a,b)$ and $(b,a)$. A graph is said to be \emph{undirected} if, for every pair of vertices $a$ and $b$, the graph either contains both arcs $(a,b)$ and $(b,a)$ or none of them; otherwise the graph is said to be \emph{directed}. Therefore our undirected graphs are special cases of directed graphs in which every arc is accompanied by its opposite, and directed graphs could also be ``mixed'', in the sense that they could contain both arcs and edges.

The \emph{valency} of a vertex $v$ in a directed graph is equal to the number of arcs of the form $(v,x)$; for an undirected graph this is equal to the number of edges containing that vertex.  In our generating functions, we usually denote valency by the letter $r$.  

Finally, $I_{(G,X)}$ will denote the cycle index of the permutation group $G$ acting on the set $X$.

\subsection{The structural approach: an introduction}

The group ring $\langle\mathbb{Z}\lbrack\mathbb{Z}_{n}\rbrack;+, \cdot\rangle$ of $\mathbb{Z}_{n}$ over $\mathbb{Z}$, consists of the set of all formal linear combinations of elements of $\mathbb{Z}_{n}$ with integral coefficients, that is, all formal sums $\sum_{h\in\mathbb{Z}_{n}}\alpha_{h}\underline{h}$ with $\alpha_{h}\in\mathbb{Z}, h\in\mathbb{Z}_{n}$, together with addition
\[\sum_{h\in\mathbb{Z}_{n}}\alpha_{h}\underline{h}+\sum_{h\in\mathbb{Z}_{n}}\beta_{h}\underline{h}:=\sum_{h\in\mathbb{Z}_{n}}(\alpha_{h}+\beta_{h})\underline{h}\]
and formal multiplication
\[\left(\sum_{h\in\mathbb{Z}_{n}}\alpha_{h}\underline{h}\right)\cdot\left(\sum_{k\in\mathbb{Z}_{n}}\beta_{k}\underline{k}\right):=\sum_{h,k\in\mathbb{Z}_{n}}\alpha_{h}\beta_{k}\left(\underline{h+k}\right)=\sum_{h\in\mathbb{Z}_{n}}\left(\sum_{k\in\mathbb{Z}_{n}}\alpha_{h-k}\beta_{k}\right)\underline{h}.\]

Note that we are writing $\underline{h}$ for $h\in \mathbb{Z}_n$ in order to distinguish clearly between elements of $\mathbb{Z}_n$ and $\mathbb{Z}$.

The elements of $\mathbb{Z}[\mathbb{Z}_n]$ also satisfy the \emph{Schur-Hadamard product} defined as follows

\[\left(\sum_{h\in\mathbb{Z}_{n}}\alpha_{h}\underline{h}\right)\circ\left(\sum_{h\in\mathbb{Z}_{n}}\beta_{h}\underline{h}\right):=\sum_{h\in\mathbb{Z}_{n}}(\alpha_{h}\beta_{h})\underline{h}\]
Therefore, for $T, T' \subseteq \mathbb{Z}_{n} $ we have $\underline{T} \circ \underline{T'}=\underline{T \cap T'}.$

The $\mathbb{Z}$-submodule of $\mathbb{Z}\lbrack\mathbb{Z}_{n}\rbrack$ generated by elements $\lambda_{1},...,\lambda_{r} \in \mathbb{Z}\lbrack\mathbb{Z}_{n}\rbrack$ will be denoted by
\[\langle \lambda_{1},...,\lambda_{r}\rangle.\]
Therefore the $\mathbb{Z}$-submodule $\langle \lambda_{1},...,\lambda_{r}\rangle$, consists of all linear combinations of $\lambda_{1},...,\lambda_{r}$ and their products.

Assume $T\subseteq\mathbb{Z}_{n}, T=\lbrace t_{1},t_{2},...,t_{r} \rbrace$. Elements of the form
\[\underline{T}:=\sum_{h\in T}\underline{h}\]
are called \emph{simple quantities} of $\mathbb{Z}\lbrack \mathbb{Z}_{n} \rbrack$. One can consider $\underline{T}$ as the formal sum $\sum_{h\in\mathbb{Z}_{n}}\alpha_{h}\underline{h}$ with $\alpha_{h}=1$ if and only if $h \in T$ and $\alpha_{h}=0$ otherwise, that is, a simple quantity is a list in which every entry has multiplicity 1. For $T=\lbrace t_{1},t_{2},...,t_{r} \rbrace$ we use the notation
\[\underline{t_{1},\ldots,t_{r}}\]
instead of $\lbrace t_{1},\ldots,t_{r}\rbrace$.

A subring $\mathcal{S}$ of a group ring $\mathbb{Z}[\mathbb{Z}_n]$ is called a \emph{Schur ring}  $\mathfrak{S}$ or $\mathcal{S}$-ring over $\mathbb{Z}_n$, of rank $r$ if the following conditions hold:
\begin{enumerate}
\item $\mathcal{S}$ is closed under addition and multiplication with elements from $\mathbb{Z}$ (i.e. $\mathcal{S}$ is a $\mathbb{Z}$-module);
\item Simple quantities $\underline{T}_{0},\underline{T}_{1},...,\underline{T}_{r-1}$ exist in $\mathcal{S}$ such that every element $\sigma \in \mathcal{S}$ has a unique representation;
\[\sigma=\sum_{i=0}^{r-1}\sigma_{i}\underline{T}_{i}\]
\item $\underline{T}_{0}=\underline{0}$, $\sum_{i=0}^{r-1}\underline{T}_{i}=\underline{\mathbb{Z}}_n$, that is, $\lbrace T_{0},T_{1},\ldots,T_{r-1}\rbrace$ is a partition of $\mathbb{Z}_n$;
\item For every $i\in \lbrace 0,1,2,...,r-1 \rbrace$ there exists a $j \in \lbrace 0,1,2,\ldots,r-1 \rbrace$ such that $\underline{T}_{j}=\underline{-T}_{i} (=\underline{\lbrace n-x : x \in T_{i}} \rbrace )$ (therefore, $\underline{T_{i}}^{t}=\underline{T_j}$);
\item For $i,j \in \lbrace 1,...,r \rbrace$, there exist non-negative integers $p_{ij}^{k}$ called structure constants, such that
\[\underline{T}_{i}\cdot\underline{T}_{j}=\sum_{k=1}^{r} p_{ij}^{k}\underline{T}_{k}\]
\end{enumerate}

The simple quantities $\underline{T}_{0},\underline{T}_{1},...,\underline{T}_{r-1}$ form a standard basis for $\mathfrak{S}$ and their corresponding sets $T_{i}$ are \emph{basic sets} of the $\mathcal{S}$-ring. The circulant graphs $\Gamma_i=\Cay(\mathbb{Z}_n,T_i)$, where $0\leq i \leq r-1$, are called \emph{basic circulant graphs} \cite{klin&al91}. The following notation will denote a $\mathcal{S}$-ring generated by its basic sets $\underline{T}_{0},\underline{T}_{1},...\underline{T}_{r-1}$:
\[\mathfrak{S}=\langle\underline{T}_{0},\underline{T}_{1},...\underline{T}_{r-1}\rangle.\]
Note that both $\mathbb{Z}(\mathbb{Z}_n)$ and $\langle\underline{0},\underline{\mathbb{Z}_n-\lbrace 0 \rbrace}\rangle$ are Schur rings over $\mathbb{Z}_n$ which we call the trivial Schur rings over $\mathbb{Z}_n$.

A permutation $g:\mathbb{Z}_n\rightarrow\mathbb{Z}_n$ is called an automorphism of an $\mathcal{S}$-ring $\mathfrak{S}$, if it is an automorphism of every graph $\Gamma_{i}$. Equivalently, the intersection of the automorphism groups of the basic circulant graphs of an $\mathcal{S}$-ring $\mathfrak{S}=\langle \underline{T}_0,\underline{T}_1,...,\underline{T}_{r-1} \rangle$, gives the automorphism group of the $\mathcal{S}$-ring.

\begin{equation}
Aut \mathfrak{S}:=\bigcap_{i=0}^{r-1}Aut \Gamma_{i}
\end{equation}

The structural approach to the enumeration of circulants on $n$ vertices is based on the lattice $\mathcal{L}(n)$ of all Schur rings over $\mathbb{Z}_n$ which, together with information on the automorphism groups of the Schur rings, suffices to carry out the enumeration. This enumeration scheme has already been described in \cite{KLP96}, so we give here only a brief summary.

We first use the lattice of Schur rings to count the number of labelled circulant graphs, as follows.
\begin{enumerate}
\item Construct the lattice $\mathcal{L}(n)$ of all Schur rings as a sequence $\mathcal{L}(n)=(\mathfrak{S}_{1},\mathfrak{S}_{2},...\mathfrak{S}_{s})$ such that $\mathfrak{S}_{j}\subseteq\mathfrak{S}_{i}$ implies $j\leq i$;
\item For directed circulants, let $\tilde{d}_{ir}$ be the number of $r$-element basis sets of the $\mathcal{S}$-ring $\mathfrak{S}_{i}$, different from the basis set $T_{0}=\lbrace 0 \rbrace$, that is,
\begin{equation*}
\tilde{d}_{ir}:=\vert \lbrace T_{(x)}\in \mathfrak{S}_{i}\vert \hspace{0.7mm} x\neq 0 \hspace{0.7mm} \mbox{ and }\hspace{0.7mm} \vert T_{(x)}\vert=r \rbrace\vert
\end{equation*}

\item For undirected circulants, let $d_{ir}$ be the number of $r$-element symmetrized (that is closed under taking of inverses) basis sets of $\mathfrak{S}_{i}$, different from $T_{0}$. That is,
\begin{equation*}
d_{ir}:=\vert\lbrace T_{(x)}^{sym} \vert \hspace{0.7mm} x\neq 0 \hspace{0.7mm} \mbox { and } \hspace{0.7mm} \vert T_{(x)}^{sym}\vert=r \rbrace\vert
\end{equation*}
\item Enumeration of all labelled directed and undirected circulant graphs which belong to the Schur ring $\mathfrak{S}_{i}$ may then be carried out by making use of generating functions $\tilde{f}_{i}(t)$ and $f_i(t)$ respectively, given by:
\begin{equation}\label{eq:erl}
\begin{split}
\tilde{f}_{i}(t)&:=\sum_{r=0}^{n-1}\tilde{f}_{ir}t^{r}:=\prod_{r=1}^{n-1}(1+t^{r})^{\tilde{d}_{ir}}\\
f_i(t)&:=\sum_{r=0}^{n-1}f_{ir}t^{r}:=\prod_{r=1}^{n-1}(1+t^{r})^{d_{ir}}
\end{split}
\end{equation}

Substituting $t=1$ in the generating functions, would give us the number of all labelled directed and undirected circulant graphs in $\mathfrak{S}_{i}$. In addition, the graph corresponding to $T \in \mathfrak{S}_{i}$ is of valency $r$ if $T$ has $r$ elements.
\end{enumerate}

The link between the number of labelled and unlabelled circulant graphs is given by this result.

\begin{Lemma}[\cite{KLP96}]
Let $G_{i}= \mbox{ Aut }(\mathfrak{S}_{i})$, let $N(G_{i})=N_{S_{n}}(G_{i})$ be the normalizer of the group $G_{i}$ in $S_{n}$, and let $\Gamma$ be a circulant graph belonging to $\mathfrak{S}_{i}$. Then

\begin{itemize}
\item[\text{(a)}] $\mbox{ Aut }(\Gamma)=G_{i}\Longleftrightarrow \Gamma$ generates $\mathfrak{S}_{i}$.
\item[(b)] If $\mbox{ Aut }(\Gamma)=G_{i}$ then there are exactly $\lbrack N(G_{i}):G_{i} \rbrack$ (that is, equal to the number of cosets of $G_i$ in $N(G_{i})$) distinct circulant graphs which are isomorphic to $\Gamma$.
\end{itemize}
\end{Lemma}

So, let the generating function for the number of non-isomorphic undirected circulant graphs with automorphism group $G_{i}$ be given by
\begin{equation*}
g_i(t)=\sum_{r=0}^{n-1}g_{ir}t^{r}
\end{equation*}
and let the generating function for the number of  non-isomorphic directed circulant graphs with automorphism group $G_{i}$ be given by
\begin{equation*}
\tilde{g}_{i}(t)=\sum_{r=0}^{n-1}\tilde{g}_{ir}t^{r}
\end{equation*}
(In all our generating functions, the coefficient of $t^r$ equals the number of circulants under consideration in which all vertices have valency $r$.

Moreover, let
\begin{equation*}
g(t)=g(n,t) \mbox{ and } \tilde{g}(t)=\tilde{g}(n,t)
\end{equation*}
denote the generating functions for the number of  non-isomorphic undirected and directed circulant graphs, respectively, with $n$ vertices. The values $g(1)$ and $\tilde{g}(1)$ therefore give the numbers of all non-isomorphic undirected and directed circulant graphs, respectively, with $n$ vertices. These generating functions are then given by the following theorem whose proof is based on the inclusion-exclusion principle.

\begin{Theorem}[\cite{KLP96}]
\label{thm:Theorem5}
\begin{equation}\label{eq:erl1}
\begin{split}
g_{i}(t)&=\frac{\vert G_{i} \vert}{\vert N(G_i)\vert} \left( f_{i}(t)-\sum_{\mathfrak{S}_{j} \subseteq \mathfrak{S}_{i}}\frac{\vert N(G_{j})\vert}{\vert G_{j} \vert}g_{j}(t)\right),\\
\tilde{g}_{i}(t)&=\frac{\vert G_{i} \vert}{\vert N(G_i)\vert}\left(\tilde{f}_{i}(t)-\sum_{\mathfrak{S}_{j}\subseteq\mathfrak{S}_{i}}\frac{\vert N(G_{j})\vert}{\vert G_{j} \vert}\tilde{g}_{j}(t)\right),\\
g(t)&=\sum_{i=1}^{s}g_{i}(t), \hspace{3mm} \tilde{g}(t)=\sum_{i=1}^{s}\tilde{g}_{i}(t).
\end{split}
\end{equation}
\end{Theorem}

In Section \ref{sec:structural} we shall give a few simple examples of this approach towards the enumeration of circulant graphs.

\subsection{The multiplier approach: an introduction}

Let $\mathbb{Z}_n^*$ be the multiplicative group consisting of all the units in $\mathbb{Z}_n$ (when $n$ is prime, $\mathbb{Z}_n^*=\mathbb{Z}_n-\{0\}$).
It is clear that if $\Gamma_1=\Cay(\mathbb{Z}_n,S)$ and $\Gamma_2=\Cay(\mathbb{Z}_n,T)$ are circulants such that there exists an $m\in \mathbb{Z}_n^*$ with $mS=\{ms:s\in S\}=T$, then $\Gamma_1$ and $\Gamma_2$ are isomorphic. In this case we say that the connecting sets are \emph{equivalent}. In \cite{adam67} \'Ad\'am conjectured that the converse is also true, that is, two isomorphic circulant graphs have equivalent connecting sets. This conjecture turned out to be false. The following is the smallest counterexample, found by Elspas and Turner \cite{elspas&tur70}. It is a pair of directed circulants. In $\mathbb{Z}_8$, let $S=\{1,2,5\}$ and $T=\{1,5,6\}$, and let $\Gamma_1, \Gamma_2$ be the corresponding circulant graphs $\Cay(\mathbb{Z}_8,S)$ and $\Cay(\mathbb{Z}_8,T)$. Then the sets $S, T$ are not equivalent but $\Gamma_1, \Gamma_2$ are isomorphic via the map 
\[ i \mapsto 4\left\lfloor\frac{i+1}{2}\right\rfloor + i.\] 
Further counterexamples with undirected circulants were also subsequently found.

The principal theorem which gives the most correct version of \adm's Conjecture is the following due to Muzychuk \cite{muzychuk95}.

\begin{Theorem}
Let $\Gamma_1$ and $\Gamma_2$ be two circulant graphs on $n$ vertices, and suppose that $n$ is square-free. Then $\Gamma_1, \Gamma_2$ are isomorphic if and only if their connecting sets are equivalent. 
\end{Theorem}

In any enumeration problem, determining when two objects are ``isomorphic'' is an essential step. Muzychuk's Theorem therefore divides the problem into two classes: when $n$ is square free and when $n$ has repeated prime factors. The easiest square-free case occurs when $n$ is prime, and the fact that, in this case, \adm's Conjecture holds, was first proved by Elspas and Turner \cite{elspas&tur70}. This reduced the problem of enumerating circulant graphs on a prime number of vertices to that of determining the number of subsets of $\mathbb{Z}^*_p$ which are not similar under the regular action of the multiplicative group $\mathbb{Z}^*_p$ on itself. 
Elspas and Turner used this method to count the number of directed and undirected circulants on $p$ vertices by means of a clever use of P\'olya's enumeration theorem.

In view of this result and Muzychuk's Theorem, the natural non-square-free cases to consider would be when the order $n$ is a power $k$ of a prime, that is, $n=p^k$, for $k\geq 2$. But to enumerate circulant graphs of such an order requires some multiplicative relations between the connecting sets of two circulant graphs which are necessary and sufficient for them to be isomorphic, that is, we require the correct version of \adm's Conjecture for $n=p^k$. We call this method of enumerating non-isomorphic circulants the multiplier approach. We shall consider in some detail the multiplier approach for $k=2, 3$ and $p=3,5$ in Section \ref{sec:multiplier}. 


\section{On the automorphism groups of prime cubed circulants}
\label{sec:misha}

We pause here to provide a brief survey of some background concepts and results in algebra graph theory, which may be used for a full justification of the theoretical results at the basis of our computational approach to constructive and analytical enumeration of circulant graphs as developed below. As a rule, rigorous and precise proofs are avoided in this section. Our modest aim in this section is to help the reader achieve a satisfactory intuitive feel for these results and after that the more interested reader, may be able to obtain a full understanding of how these proofs are obtained. Such an understanding is, however, not required in order to follow the arguments presented in the subsequent sections of the paper. At the end of this section, some source references are given. An interested reader with the aid of these texts may reach full lucidity, which however is not required in order to follow the main lines of presentation in the current paper.

We denote the cyclic group of order $n$ by $\Z_n$. Usually, this notation implies the additive group modulo $n$. Simultaneously we denote by $(\Z_n,\Z_n)$ the regular cyclic permutation group  acting on the set ${0,1,\ldots,n-1}$ and generated by the cyclic shift $(0,1,2,\ldots,n-1)$.

The automorphism group $\Z_n^*=\Aut(\Z_n)$ of the group $\Z_n$ is the multiplicative group, modulo $n$, of units of $\Z_n$, denoted by $\Z_n^*$. It has order $\varphi(n)$, where $\varphi$ is the famous Euler function. The group $\Z_n^*$ acts on $\Z_n$ by multiplication modulo n. For our goals it is enough to consider the case $n=p^k$ where $p$ is prime and $k$ is mainly 1, 2 or 3. It is well-known that $\varphi(p^k)=(p-1)p^{k-1}$. For these values of $n$ the group $\Z_n^*$ is also cyclic.

\medskip\noindent
More generally, let $(G,\Omega)$ be a finite permutation group acting on the set $\Omega$. Denote by $2-\mbox{orb}(G,\Omega)$ the set of all 2-orbits of $(G,\Omega)$ (in the sense of H. Wielandt), that is, the set of the orbits of the natural induced action $(G,\Omega^2)$ as follows: for $g\in G$ and $(\alpha,\beta)\in \Omega^2$acting naturally on $\Omega\times\Omega$, $(\alpha,\beta)^g=(\alpha^g,\beta^g)$, where $x^g$ is the image of the element $x\in \Omega$ under the action of $g\in G$.  Let $2-\mbox{orb}(G,\Omega)=\{R_i: 0\leq i \leq r-1\}$. Each pair $\Gamma_i=(\Omega,R_i)$ is regarded as a directed graph with vertex-set $\Omega$. Then the group $(G^{(2)},\Omega)$, where 
\[G^{(2)} = \cap_{i=0}^{r-1} \mbox{Aut}(\Gamma_i),\]
is called the \emph{$2$-closure} of $(G,\Omega)$. The group $(G,\Omega)$ is called \emph{$2$-closed} if $(G^{(2)}, \Omega)= (G,\Omega)$.

A classical (almost trivial) result by Wielandt claims that each regular permutation group is 2-closed. Thus, in particular, $(\Z_n,\Z_n)$ is 2-closed. (Note that later on the notation $(\Z_n,\Z_n)$ may be reduced to $\Z_n$ if is it clear from the context that we mean the regular action of $\Z_n$.)

\medskip\noindent
The main target of interest in this section is the lattice of all 2-closed overgroups of the group $(\Z_n, \Z_n)$. Here, by an overgroup of $\Z_n$ we understand a subgroup $G$ of $S_n$, the symmetric group on $\{0,1,2,3,\ldots,n-1\}$, which contains $(\Z_n,\Z_n)$. Note that here $Z_n$ appears in two different roles: the set of elements $\{0,1,2,\ldots,n-1\}$ and the set of permutations.

It turns out that for the case $n=p^k$ there exists an anti-isomorphism 
between the lattice of all 2-closed overgroups of $\Z_n$ in $S_n$ and the lattice  $\mathfrak{S}_n$ of all $S$-rings over $\Z_n$.

For arbitrary values of $n$, establishment of such a bijection between overgroups and $S$-rings turns out to be more sophisticated: one has to consider only so called Schurian $S$-rings. The good news for the case $n=p^k$ is that here all $S$-rings are Schurian, a result due to R. P\"{o}schel \cite{pos74}.

For all values of $n$ there are two trivial overgroups, the minimal $(\Z_n, \Z_n)$ and the maximal $(S_n, \Z_n)$. All other overgroups appear between these two extremal objects. 

\subsection{Wreath products}

Let $(G_1, M_1)$ and $(G_2, M_2)$ be two permutation groups and let $M=M_1\times M_2$. Let $G$ be the set of all mappings $g:M\to M$ such that, for $x=(x_1,x_2)\in M$, $x^g=(x_1,x_2)^g=(f_1(x_1,x_2), f_2(x_1,x_2))$, the following conditions hold:
\begin{itemize}
\item[a] $f_1$ depends only on coordinate $x_1$;
\item[b] the mapping $x_1 \mapsto f_1(x_1, x_2)$ is a permutation in $G_1$;
\item[c] for every $x_1\in M_1$, the mappings $x_2 \mapsto f_2(x_1,x_2)$ are permutations which belong to $G_2$.
\end{itemize}
In this case, for brevity, the notation $g=[g_1, g_2(x_1)]$ is used and is called the \emph{table} of $g$. By definition, we we have 
\[x^g = (x_1,x_2)^g = (x_1^{g_1}, x_2^{g_2(x_1)}).\]

It is easy to check that here $G$ is a permutation group $G\leq Sym(M)$. The group $(G,M)$ is called the \emph{wreath product} of $(G_1,M_1)$ and $(G_2, M_2)$ and will be denoted by $(G_1\wr G_2, M_1\times M_2)$  or $(G_1,M_1)\wr(G_2,M_2)$.

The wreath product is a group of order $|G_i|\cdot|G_2|^{M_1}$ ($=$ the number of all tables); sometimes, $G_1$ and $G_2$ are called active and passive factors of $G$, respectively. Note also that we are using here so-called orthodixal notation for the wreath product, due to L.A. Kalu\v{z}nin. 

Let $\Gamma_1=(V_1,E_1), \Gamma_2=(V_2, E_2)$ be two graphs (directed or undirected). Then the graph $\Gamma=(V,E)$ defined by $V=V_1\times V_2$, $E=\{\left((x_1,x_2), (x_1',x_2')\right): (x_1,x_2')\in E_1 \lor (x_1=x_1' \land (x_2,x_2')\in E_2)\}$ is called the \emph{composition} of $\Gamma_1$ and $\Gamma_2$ and is usually denoted by $\Gamma_1[\Gamma_2]$ (substitute $\Gamma_2$ for each vertex of $\Gamma_1$ and connect all vertices of the corresponding components according to the connections in $\Gamma_1$).

One of the traditional questions in graph theory was to study conditions under which the automorphism group of the composed graph $\Gamma_1[\Gamma_2]$ equals the wreath product of $\Aut(\Gamma_1)$ with $\Aut(\Gamma_2)$. 

In cases where we agree to consider transitive permutation groups only, this question, in the context of the current presentation, finds a very suitable solution, namely it turns out that 
\[(G_1\wr G_2, M)^{(2)}= (G_1,M_1)^{(2)}\wr (G_2, M_2)^{(2)} \]
where $M= M_1 \times M_2$.

In order to explain the significance of this equality in the context of Schur rings, it is very natural to use one more definition. An $S$-ring $\mathcal{S} = \langle \underline{T_1}, \underline{T_2}, \ldots, \underline{T_d}\rangle$ over $\Z_{p^k}$ is called \emph{wreath decomposable} (or briefly decomposable) if there exists a non-trivial proper subgroup $K\leq \Z_{p^k}$ such that for each basic set either $T_i\subseteq K$ or $T_i$ is a union of suitable cosets of $\Z_{p^k}/K$ (that is, $T_i = \cup _{x\in T_i}K+x)$. In particular, one has $\underline{K}\in \mathcal{S}$. The S-ring $\mathcal{S}$ is called \emph{wreath indecomposable} (briefly \emph{indecomopsable}) if it is not wreath decomposable.

Using these concepts and facts, one can prove that for an arbitrary $S$-ring $\mathcal{S}$ over $\Z_{p^k}$ is true that $\mathcal{S}$ is wreath decomposable if and only if $\Aut(\mathcal{S})$ can be represented as the wreath product of the automorphism groups of $S$-ring over $\Z_{p^{k-i}}$
and $S$-ring over $\Z_{p^i}$, where $1\leq i < k$. 

\begin{figure}[ht]  
	\centering
	\includegraphics[width=7cm, height=5cm]{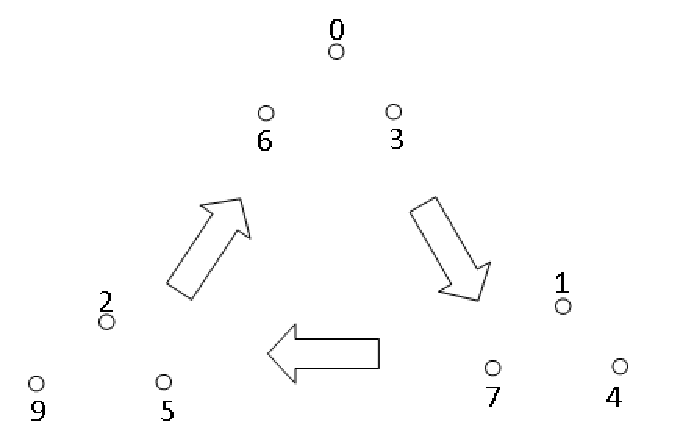}
	\caption{Circulant graph $\Cay(\Z_9, \{1,4,7\})$} \label{fig:misha1}
\end{figure}

\begin{Example}
	Let $n=9$, $\mathcal{S}_1=\langle \underline{0}, \underline{3}, \underline{6}, \underline{1,4,7}, \underline{2,5,8}\rangle$. Then $\Aut(\mathcal{S}_1)=\Z_3\wr S_3 = \Aut(\Gamma_1)$ is a group of order $3\cdot(3!)^3=2^3\cdot3^4=648$. Here $\Gamma_1= \overrightarrow{C}_3[E_3]$ is the composition of directed cycle $\overrightarrow{C}_3$ with the empty 3-vertex graph $E_3$. The graph $\Gamma_1$ is a Cayley graph $\Cay(\Z_9, \{1,4,7\})$, that is, a 9-vertex circulant. Its mnemonic diagram is depicted in Figure \ref{fig:misha1}. Here a big directed arrow substitutes nine arcs from each vertex of initial 3-vertex set to each vertex of targeted 3-vertex set. 
\end{Example}

In what follows, indecomposable $S$-rings over groups $\Z_{p^i}$ will be called atoms. It turns out that if we understand the structure of all atoms then we understand the automorphism groups of all $S$-rings over $\Z_{p^k}$.

\subsection{Affine overgroups of $(\Z_n, \Z_n)$}

Recall that there are two commutative operations of addition and multiplication on the set $\Z_n$. Thus the structure $(\Z_n, +, \cdot)$ forms the classical prototype for all finite commutative rings. The elements of $\Z_n^*$ form the group of invertible elements of the ring $\Z_n$. 
  
Let $\Aff(1,n)$ be the group of all one-dimensional linear transformations over $\Z_n$, that is,
\[ \Aff(1,n):= \{\mathcal{M}_{a,b}: a\in \Z_n^*, b\in \Z_n\},\] 
where the \emph{affine transformation} $\mathcal{M}_{a,b}$ of $\Z_n$ is given by 
\[ \mathcal{M}_{a,b}:x\mapsto ax+b \quad (x\in \Z_n). \]
  
The following well-known facts prove to be very helpful in our context.
\begin{enumerate}
\item Every $\M_{a,b}\in\Aff(1,n)$ is a permutation on $\Z_n$;
\item $(\Aff(1,n),\Z_n)$ is a permutation group;
\item[c)] $|\Aff(1,n)|=n\cdot\varphi(n)$;
\item $\Z_n\cong\{\M_{1,b}: b\in \Z_n\}$ is a normal subgroup
of Aff$(1,n)$.
\end{enumerate}

Here we restrict our considerations to the case $n=p^k$; moreover, only small values of $k$ will actually be required. Therefore
$|\Z_{p^m}^{\ast}|=\varphi(p^m)=p^{m-1}(p-1)$ and
$|\Aff(1,p^m)|=p^{2m-1}(p-1)$.  Note that the permutation $\M_{1,1}$
is nothing else than the standard cycle $(0,1,\dots,n-1)$. Every
subgroup $(G,\Z_n)$ of Aff$(1,n)$ which contains the cycle
$\M_{1,1}$ will be called an \emph{affine overgroup} of $(\Z_n,\Z_n)$ or simply an affine group. It is convenient to call Aff$(1,n)$ the \emph{complete affine group}.

Clearly, each affine group $G$ can be represented as a semidirect product 
\[ G=\Z_n\rtimes L, \]
where $L$ is a subgroup of $\Z_n^{\ast}$, while in the group $G$, the subgroup $L$ is the stabiliser $G_0$ of the element $0\in \Z$. 

This is why orbits of $L$ on the set $\Z_n$ form an $S$-ring  
over $\Z_n$. Such an $S$-ring which stems from a suitable affine group $G$ is called an \emph{affine $S$-ring} over $\Z_n$.

By the given definitions for an affine $S$-ring $\mathcal{S}$, which is obtained from $(G,\Z_n)$, the group $G$ is a subgroup of $\Aut(\mathcal{S})$. A significant issue is to understand the full group $\Aut(\mathcal{S})$, or, in other words, the 2-closure of $(G,\Z_n)$.

For the case $n=p$, the full affine group of order $p(p-1)$ is 2-transitive and its 2-closure is the symmetric group $S_p$. All other affine groups are \emph{uniprimitive}, that is, primitive but not 2-transitive. 

One of the classical results in the theory of permutation groups (due to Burnside and Schur) ia that each uniprimitive permutation group of prime degree $p$ is affine and 2-closed.

For $n=p^k, k>1$ all affine groups are uniprimitive! It turns out that here the description of the 2-closure is becoming a more involved task. Recall that a transitive permutation group $(G,M)$ is called a \emph{Frobenius group} if each non-identical permutation $g\in G$ has at most one fixed point in $M$.

\begin{Proposition} \label{prop:frobenius} Every imprimitive Frobenius group is 2-closed.
\end{Proposition}
	
Note that, in general, the proposition fails for primitive permutation groups, though it remains valid for some restricted classes, like the above-mentioned uniprimitive permutation groups of degree $p$. 

For $k=2$, Proposition \ref{prop:frobenius} allows us to detect a one-parameter family of affine 2-closed groups of order $t\cdot p^k$, where $t$ is any divisor of $p-1$. These imprimitive Frobenius groups form one class of $p^k$-atoms in the process of classification of the automorphism groups of circulants.

There exists also an efficient criterion for the indecomposability of an affine $S$-ring over $\Z_{p^k}$. To avoid technical complications, this criterion will not be exploited here, however it will be exploited implicitly in a further presentation.

At this stage we are sufficiently prepared to go ahead towards formulating the main results in this section. 

\subsection{Main results}

Let us denote by $u_n$ the number of subgroups of the multiplicative group $\Z^*_n$. It is easy to understand that for $n=p^k$, $p$ odd prime, there is the equality $u_n=k\cdot d$, where $d$ is the number of all natural divisors of $p-1$. 

Thus, in the context of the current paper, the following values are mostly significant:
\[ u_3 = 2,\quad u_9=4,\quad u_{27}=6,\quad u_5=3\quad u_{25}=6,\quad u_{125}=9. \]

\begin{Proposition} \label{prop:numberofovergroups}
	There are exactly $u_p$ $2$-closed overgroups of $\Z_p$:
	\begin{itemize}
		\item[(a)] The symmetric group $S_p$; and
		\item[(b)] Frobenius uniprimitive groups $\mathcal{F}_p^s$ of order $sp$, where $s$ is a proper divisor of $p-1$, (that is, $s<p-1$).
	\end{itemize}
\end{Proposition}
 
All the groups which appear in the formulation of Proposition \ref{prop:numberofovergroups} play the role of $p$-atoms in the recursive description of the 2-closed overgroups of $\Z{p^k}, k>1$.
 
\begin{Proposition}
	Every 2-closed overgroup of $\Z_{p^2}$ is of one of the following types:\begin{itemize}
		\item[(a)] wreath product of $p$-atoms;
		\item[(b)] $S_{p^2}$; 
		\item[(c)] Frobenius group $\mathcal{F}_{p^2}^s$, where $s$ is any divisor of $p-1$. 
	\end{itemize}
\end{Proposition}

Groups of types (b) or (c) in the above proposition will be called $p^2$-atoms.

\begin{Corollary}
	There are exactly $1+u_p+(u_p)^2$ 2-closed overgroups of $\Z_{p^2}$.
\end{Corollary}
 
The first difficulty in the description of the 2-closure of affine overgroups of $\Z_{p^k}$ appears in the case when $k=3$. Each affine group $G$ over $\Z_{p^3}$ can be presented in the form $G=\Z{p^3} \rtimes L$, where $L\leq \Z_{p^3}^*$.

It turns out that we have to classify all affine groups into three classes; according to the appearance of elements $(p+1)$ and $(p^2+1)$ in $L$.

If $(p+1)\not\in L$ and $(p^2+1)\not\in L$, then $L$ is a Frobenius (imprimitive) group, which is 2-closed.

If $(p+1)\in L$ and $(p^2+1)\in L$, then the orbits of $L$ define a decomposable $S$-ring. Thus $G$ is not 2-closed, however, $G^{(2)}$ can be described with the aid of the iterated wreath product.

\begin{Example}
	Here $n=27$, $L=\{1,4,7,10,13,16,19,22,25\}$. Clearly both $1+3$ and $1+3^2$ are in $L$. The group $G=\Z_{27}\rtimes L$ defines an $S$-ring $\mathfrak{S}_2$ with the basic elements $T_1=\underline{1,4,7,10,13,16,19,22,25},\ 2T_1,\  T_3=\underline{3,12,21}\cap 2T_3, T_5=\underline{9},\ 2T_5,\  \underline{0}$. 
\end{Example}

Analysing the lattice of $S$-rings over $\Z_{27}$, it is possible to observe that $\Aut(\mathfrak{S}_2)= \Aut(\Gamma_1)\cap \Aut(\Gamma_3)\cap \Aut(\Gamma_5)$, where the $\Gamma_i$ are the circulants defined by $T_i$.

\begin{figure}[ht]  
	\centering
	\includegraphics[width=7cm, height=5cm]{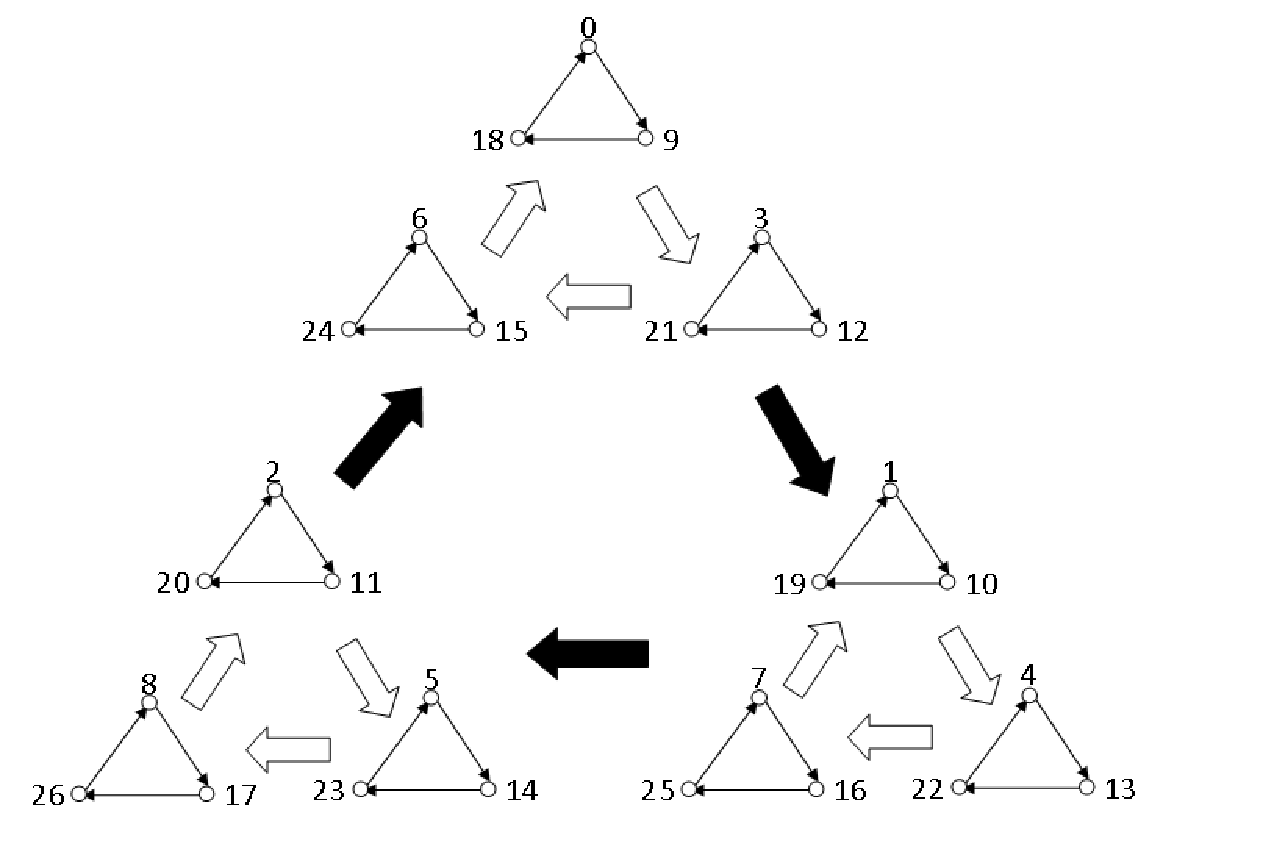}
	\caption{Graphs $\Gamma_1$ $\Gamma_3$ and $\Gamma_5$ depicted together} \label{fig:misha2}
\end{figure}

The reader is welcome to conceptualise the structure of $\Aut(\mathfrak{S}_2)$ from Figure \ref{fig:misha2}. Here, usual arrows depict arcs, large unfilled arrows have the same sense as they had before in Figure \ref{fig:misha1}, while solid black arrows substitute 81 arrows from each of nine vertices of the starting block to each of nine vertices in the target block. 
 
Clearly there is full freedom to rotate cyclically each usual triangle, to rotate cyclically each triangle consisting of three usual ones, and finally to rotate cyclically the global triangle, consisiting of three $9$-vertex ingredients. This means that $\Aut(\mathfrak{S}_2)= \Z_3\wr\Z_3\wr\Z_3$. (Note that the operation of wreath product is associative.) Thus $|\Aut(\mathfrak{S}_2|= 3\cdot(3\cdot3^3)^3= 3^{13}$.

Finally, the third most sophisticated case of affine overgroups of $\Z_{p^3}$ is when $(p+1)\not\in L$, but $(p^2+1)\in L$. In this case the affine group $G\rtimes L$ defines an indecomposable $S$-ring, however the group $G$ is not 2-closed.

Again we will try to understand this more sophisticated situation with the aid of an example. 

\begin{Example}
	Here again $n=27$. Consider $L=\{1,10,19\}$. Note that $(3+1)\not\in L$, however $(3^2+1)\in L$, thus we indeed face the third case which we declared to be the most difficult one. The group $G=\Z_{27}\rtimes L$ defines the $\mathfrak{S}_3$ with the basic quantities of lengths $3$ and $1$ as follows: 
	\[\mathfrak{S}_3= \langle \underline{0},\underline{1,10,19}, \underline{2,11,20}, \underline{8,17,26}, \underline{7,16,25}, \underline{4,13,22}, \underline{5,14,23}, \underline{3}, \underline{6}, \underline{9}, \underline{12}, \underline{15}, \underline{18}, \underline{21}, \underline{24}\rangle.\] 
\end{Example} 
 
\begin{figure}[ht]  
	\centering
	\includegraphics[width=6cm, height=5cm]{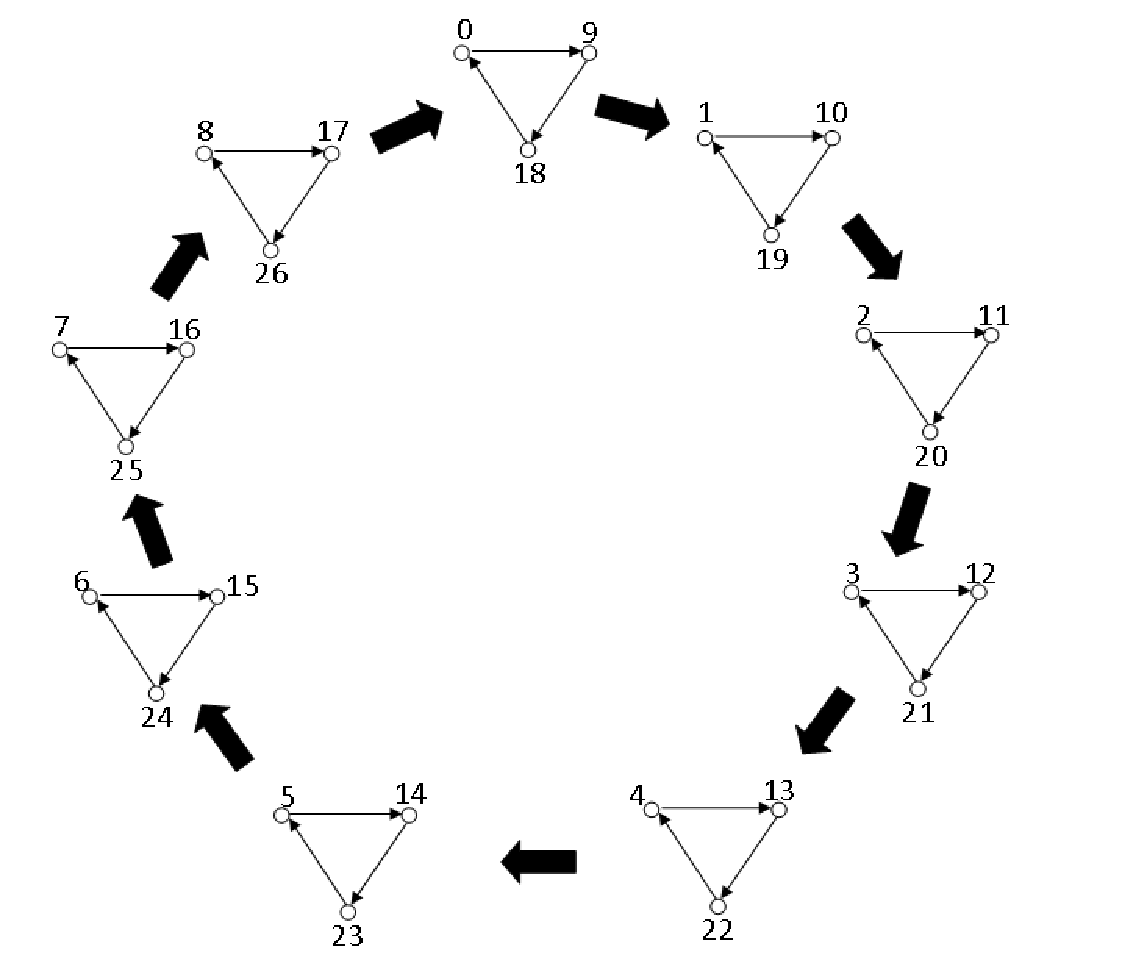}
	\caption{The circulant graphs $\Gamma_1$ (thin arcs) and $\Gamma_2$ (solid arcs)} \label{fig:misha3}
\end{figure}

It is possible to prove that in this case $\Aut(\mathfrak{S}_3)= \cap_{i=1}^3\Aut(\Gamma_i)$, where $\Gamma_1=\Cay(\Z_{27}, \{1,10,19\})$, $\Gamma_2=\Cay(\Z_{27}, \{9\})$ and $\Gamma_3=\Cay(\Z_{27},\{3\})$. These circulant graphs $\Gamma_1$ and $\Gamma_2$ are again depicted together in Figure \ref{fig:misha3} exactly in the same fashion as was done before. 

The situation with graph $\Gamma_3$ is more complicated, provided we want to arrange its vertices to afford the viewer a ``correct'' visualisation: that is to consider $\Gamma_3$ together with the two previous graphs in order to visualisen the action of the group $\Aut(\mathfrak{S}_3)$.

\begin{figure}[ht]  
	\centering
	\includegraphics[width=8cm, height=7cm]{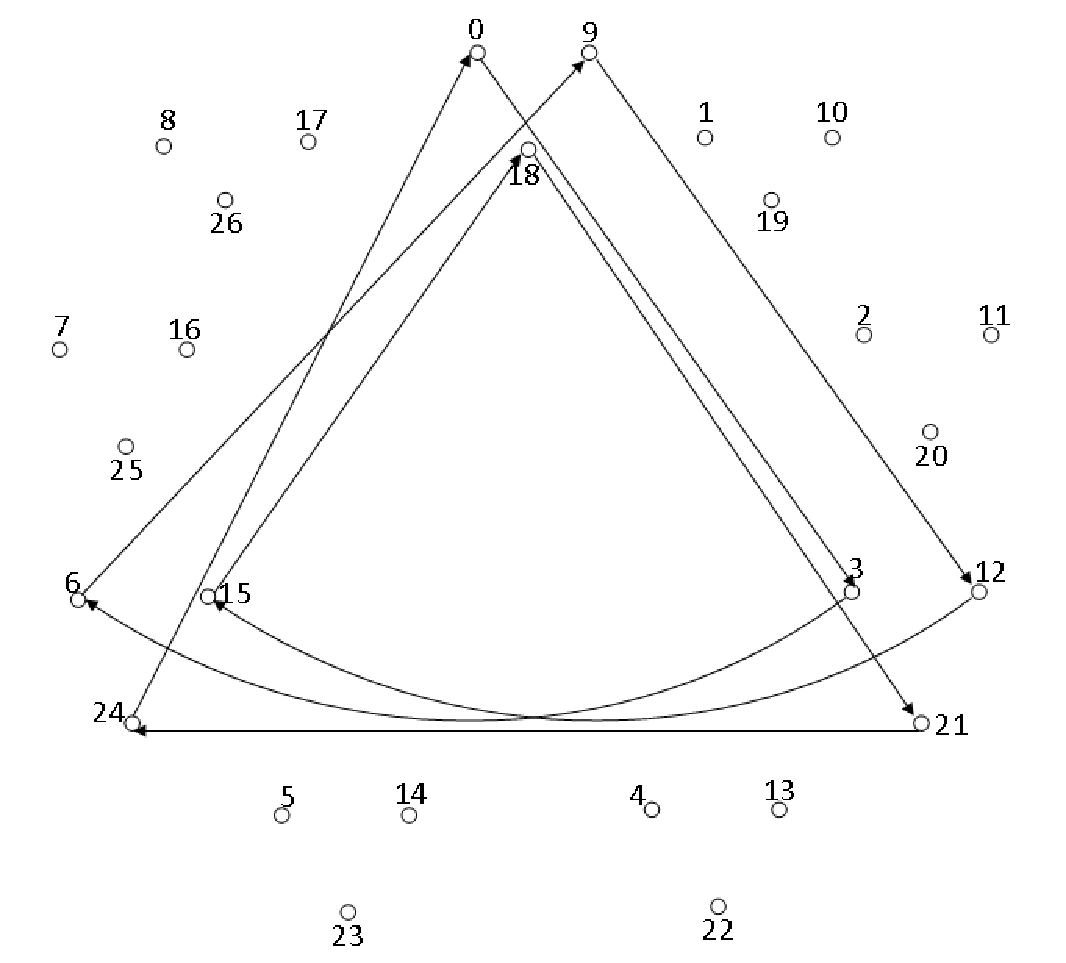}
	\caption{One of the three connected components of the graph $\Gamma_3$} \label{fig:misha4}
\end{figure}

In principle, the graph $\Gamma_3$ has a very simple structure of the form $3\overrightarrow{C}_9$, that is, the disjoint union of three directed cycles of length 9. Nevertheless, we intentionally prefer to depict it in a more sophisticated ``skew'' manner, as it appears in Figure \ref{fig:misha4}. 
 
In this figure we have deliberately shown only one of the three connected components of $\Gamma_3$ (it is the one spanning  those vertices which are multiples of 3 in $\Z_{27}$). The two other isomorphic connected components are omitted.

From Figure \ref{fig:misha4} it immediately follows that 
$\Aut(\Gamma_1)\cap \Aut(\Gamma_2)= \Z_9\wr\Z_3$. Clearly, 
$\Aut(\Gamma_3)=\Aut(3\overrightarrow{C}_9)=S_3\wr\Z_3$. Thus, in principle, $G^{(2)}= \Aut(\mathfrak{S}_3)=(\Z_9\wr\Z_3)\cap(S_3\wr\Z_3)$. However, this abstract formula is useful only if we take into account the actual form of the graph $\Gamma_3$ relative to the other two graphs. The reader is welcome to check that the group $G^{(2)}$ has the form $(\Z_3)^3\cdot\Z_9$ (note that this is not a semidirect product!) and
\[ G^{(2)} = \langle h, h_0, h_1, h_2\rangle = \mathcal{M}_{1,1}. \]
Here $\mathcal{M}_{1,1}$ is a generating element of $\Z_{27}$, while
\begin{eqnarray}
h_0 & = & (0\ 9\ 18)(3\ 12\ 21)(6\ 15\ 24)\\
h_1 & = & (1\ 10\ 19)(4\ 13\ 22)(7\ 16\ 25)\\
h_2 & = & (2\ 11\ 20)(5\ 14\ 23)(8\ 17\ 26).
\end{eqnarray}
 
Indeed, $h_0$ preserves the copy of $\overrightarrow{C}_9$ depicted in Figure \ref{fig:misha4}, fixing the two other copies of $\overrightarrow{C}_9$, preserving simultaneously graphs $\Gamma_1$ and $\Gamma_2$. The permutations $h_1$ and $h_2$ have similar interpretations.

Finally we wish to stress that while $|G|=|\Z_{27}|\cdot|L|=27\cdot3=81$, for the 2-closure we get $|G^{(2)}|=3^3\cdot9=243$. In other words, the order of the 2-closure $G^{(2)}$ is $p$ times larger in comparison with the order of the affine group $G=\Aut(\mathfrak{S}_3)$. Note also that the group $G^{(2)}$ is not wreath decomposable, that is,  it cannot be represented as a wreath product at suitable smaller atoms. According to Klin and P\"{o}schel, the group $G^{(2)}$ is called a subwreath product (see the discussion below). 

Now we are prepared for the consideration of the next theoretical claim, whuich again will be presented here without its formal justification.

\begin{Theorem}
	Every $2$-closed overgroup of $\Z_{p^3}$ belongs to one of the following types:
	\begin{itemize}
		\item[(a)] Wreath product of $p$ atoms and $p^2$-atoms or vice-versa;
		\item[(b)] wreath product of three $p$-atoms;
		\item[(c)] $S_{p^3}$;
		\item[(d)] A Frobenius group $\Z_{p^3}\rtimes L$ such that $(p+1)\not\in L$ and $(p^2+1)\not\in L$;
		\item[(e)] The permutation group $(G^{(2)}, \Z_{p^3})$, where $G^{(2)}$ is the 2-closure of $G=\Z_{p^3}\rtimes L$, such that $(p^2+1)\in L$, but $(p+1)\not\in L$. 	\end{itemize}
\end{Theorem}
 
\begin{Corollary}
	There exist exactly $1+4u_p+ 2(u_p)^2+(u_p)^3$ different $2$-closed overgroups of $\Z_{p^3}$.
\end{Corollary}

As it will be explained in the next section, in order to install the structural approach for the enumeration of circulants, besides knowledge of all 2-closed overgroups of $\Z_{p^k}$, the orders of their normalisers (in $S_{p^k}$) is also required. In general, we were not trying to obtain this information on a theoretical level since, for small values of $p$ this information may be obtained with the aid of computer algebra packages such as {\tt GAP}.

\medskip\noindent
The alternative, so-called multiplier approach, relies on some implicit information regarding the behaviour of the normalisers of 2-closed overgroups of $\Z_{p^k}$, namely, the knowledge of necessary and sufficient conditions of isomorphism of circulants is enough for this approach. For the case $n=p^k$, $p$ odd, these conditions were formulated in terms very suitable for the purposes of analytical information. For general $n$, the problem of analytical enumeration of the $n$-vertex circulants might become a subject of future special attention. 

\subsection{Brief historical summary} 
 
 The crucial concept of Schur ring, which is used in this paper, goes back to the seminal paper \cite{sch33,sch03}. For more than two decades this work was known only to a few colleagues and followers of Schur. In a sense, this concept was revitalised by H. Wielandt, in particular, due to a special chapter on $S$-rings in his classic book \cite{wie64}.
  Nowadays, $S$-rungs are discussed in many modern textbooks on group theory, among them we should mention \cite{sco87} and \cite{dixM96}.
 
 The concept of $k$-closure, an in particular of 2-closure, of a given permutation group belongs to H. Wielandt \cite{wie69}. Some of its roots are attributed by Wielandt to M. Krasver, a collaborator of L.A. Kalu\v{z}{n}in, the scientific advisor of the author M.K. and R. P\"{o}schel. Note also that in his early stages of education, Kalu\v{z}{n}in was strongly influenced by Schur (attending, for a couple of years, the home seminar of Schur at the time when Schur was pushed out by the Nazis from his professorship at the University of Berlin). We refer, for more details, to the biographies of Schur and Kalu\v{z}{n}in on the famous online source ``The MacTutor History of Mathematics''. 
 
 It is worthy to mention that the methodology of invariant relations of permutation groups (conceptually quite close to Wielandt's $k$-closures) was developed in the school of Kalu\v{z}{n}in. First, very successful, stages of this development are reflected in the monograph \cite{posK79}. 
 
 For about 40 years since their creation, $S$-rings were used and considered quite sporadically. A significant step was taken by P\"{o}schel \cite{pos74}, who achieved full classification of $S$-rings over $\Z{p^k}$, $p$ an odd prime. The results of P\"{o}schel strongly influenced M.K. who, quite soon, suggested to apply P\"{o}schel's results to graph theory in order to recognise the structure of sutomorphism groups of circulants and to elaborate criteria for their isomorphism. The abstract \cite{klinP75} of the lecture, presented at the famous Zykov seminar, is the first document reflecting the start of the use of $S$-rings in combinatorics. 
 The preprint \cite{kliP78} and the 
 the paper \cite{klinP81} contain more systematic presentations. We also mention that Chapter 8 in \cite{posK79} was written in collaboration with M.K. This chapter was definitely the first attempt to consider together at the level of a monograph, $k$-closures, invariant relations, relational algebras and Schur rings, especially over cyclic groups. 
 
A paper \cite{vysKCh78} is a direct predecessor of the current text. It relies on \cite{klinP81} and presents first attempts to  create (with a computer) a full catalogue of $S$-rings over the cyclic group $\Z_{125}$. This catalogue contains exactly 58 $S$-rings. Note that $u_5=3$ and $1+4\cdot3 + 2\cdot3^2+3^3=58$. At the same time, S.P. Yushchenko prepared at Kiev State University a masters thesis (unpublished) where the results of some attempts of enumeration of circulants on 27 vertices were presented. The background of this thesis was also a full catalogue of $S$-rings over $\Z_{27}$.

This successful experience of using Schur rings in applied combinatorics served as a strong motivation for further theoretical activities as well as for more involved use of computers. In the preprint \cite{kliP80}, necessary and sufficient conditions for the isomorphism of circulant graphs with $p^k$ vertices were formulated and justified, based on \cite{pos74}. These conditions, as was mentioned, created a background for what is now called the multiplier approach.

Already, at that stage, it became clear that in the problem of description of automorphism groups of circulants, a crucial step was to move from $n=p^2$ to $n=p^3$. Main ideas in this direction were announced in the abstract \cite{kli94} and also in \cite{VL17} and \cite{LP2000}. 

Another pioneering text \cite{KLP96} dealt with analytical enumeration of prime-squared circulants. It inspired a few other papers in this direction, co-authored by V.L., to be mentioned later on in the current paper.

It is worth mentioning that the automorphism groups of $p^k$-vertex circulants were considered by some other authors, among them E. Dobson. His unpublished preprint \cite{dob95} presented a successful attempt towards classification of such groups. Becoming acquainted at that time with the results of M.K. et al., based on $S$-rings, Dobson decided to postpone publication of his preprint, aiming to develop his own self-contained approach. A number of related results were published later on, (see, for example, \cite{dob00,dobW02,dob05,dobM09,dob10}). Some of these papers consider diverse natural extensions of the original problem, related to general Cayley graphs.

The same problem was also considered by I. Kovacs \cite{kov02}. This paper was influenced by M.K.: originally Kovacs was not familar at all with the approach based on the use of $S$-rings. Finally the author elaborated his own, quite original way, using essentially spectral techniques. As a result, recursive description of automorphism groups was reached till $k=4$, hopefully convincing the reader that, in principle, in this way the problem may be resolved for arbitrary values of $k$.

Circulants with $2^k$ vertices provide another interesting (slightly more sophisticated) line of investigation, which is not touched upon in this paper, although it also has a reasonably striking history.

First, a computer algorithm was elaborated in order to describe all $S$-rings over $\Z_{2^k}$. It was implemented on the computer EC1020, and results were obtained in 1981 for $k\leq 6$ in \cite{kliNajP81}. The number of all $S$-rings for $k=3, 4, 5, 6$ was equal to $10, 37, 151$ and $657$, respectively. After that, a purely theoretical generalisation was reached, full description of $S$-rings over $\Z_{2^k}$ was announced and justified in 
\cite{kliNajP81}. Two decades later, Kovacs, basing himself on the results in \cite{klinP81,golNP85} and acting in the spirit of \cite{LP2000}, described the total number of indecomposable Schur rings over $\Z_{2^k}$, using Catalan and Schr\"{o}der numbers. Analytical enumeration of of $2^k$-vertex circulants has still not been tackled.

Finally, we describe very briefly the general problem of isomorphism of circulant graphs. Its consideration goes back to \cite{muzKP01}. A significant background was provided by a few papers by K.H. Leung and L.L. Ma, as well as in further papers by Muzychuk. Final dots were put in a significant paper  \cite{muz04} where a full solution of the isomorphism problem for circulant graphs with an arbitrary number of vertices was presented. The description of the automorphism groups of so-called rational circulant graphs can be found in \cite{kliK12}. This text contains also a historical digest together with quite a rich bibliography.

The recent paper \cite{mis15} describes a new computational approach to the counting of the number of $S$-rings over cyclic groups of order $n$, with special attention to the prime-power case. This approach, which is also based on $S$-rings, correlates with classical results in the enumeration of $p$-groups, stemming from a conjecture of Graham Hugman. 

\nocite{faradzev&l90}

\section{The structural approach} \label{sec:structural}

\subsection{The case $p^2$ for $p=3$}

We shall first describe the structural approach for $p^2$ with $p=3$. This work has already been shown in \cite{KLP96} but we present it here in order to illustrate the method. Using the techniques of wreath products and wreath decomposition of Schur rings as described in \cite{KLP96} one obtains that the following is the list of all Schur rings over $\mathbb{Z}_9$

\begin{equation*}
\begin{split}
\mathfrak{S}_{1}&=\langle \underline{0},\underline{1,2,3,4,5,6,7,8}\rangle,\\
\mathfrak{S}_{2}&=\langle \underline{0},\underline{1,2,4,5,7,8},\underline{3,6}\rangle,\\
\mathfrak{S}_{3}&=\langle \underline{0},\underline{1,4,7},\underline{2,5,8},\underline{3,6}\rangle,\\
\mathfrak{S}_{4}&=\langle \underline{0},\underline{1,2,4,5,7,8},\underline{3},\underline{6}\rangle,\\
\mathfrak{S}_{5}&=\langle \underline{0},\underline{1,4,7},\underline{2,5,8},\underline{3},\underline{6}\rangle,\\
\mathfrak{S}_{6}&=\langle \underline{0},\underline{1,8},\underline{2,7},\underline{3,6},\underline{4,5}\rangle,\\
\mathfrak{S}_{7}&=\langle \underline{0},\underline{1},\underline{2},\underline{3},\underline{4},\underline{5},\underline{6},\underline{7},\underline{8}\rangle,\\
\end{split}
\end{equation*}

We now show how this list can be used to enumerate all circulant graphs of order 9. 

We give only the briefest necessary information about the automorphism groups of all S-rings in $\mathcal{L}$. These were obtained using GAP (see \cite{gap99}), the automorphism groups being the intersection of the automorphism groups of the basic Cayley graphs associated with each Schur-ring. 

\begin{quote}
\centering
\begin{tabular}{l l}
Automorphism Group&Normalizer\\
$G_{1}=S_{9}$,&$\lbrack N(G_{1}):G_{1}\rbrack=1$,\\
$G_{2}=S_{3}\wr S_{3}$,&$\lbrack N(G_{2}):G_{2}\rbrack=1$,\\
$G_{3}=\mathbb{Z}_{3}\wr S_{3}$,&$\lbrack N(G_{3}):G_{3}\rbrack=2$,\\
$G_{4}=S_{3}\wr \mathbb{Z}_{3}$,&$\lbrack N(G_{4}):G_{4}\rbrack=2$,\\
$G_{5}=\mathbb{Z}_{3}\wr\mathbb{Z}_{3}$,&$\lbrack N(G_{5}):G_{5}\rbrack=4$,\\
$G_{6}=D_{9}$,&$\lbrack N(G_{6}):G_{6}\rbrack=3$,\\
$G_{7}=\mathbb{Z}_{9}$,&$\lbrack N(G_{7}):G_{7}\rbrack=6$.
\end {tabular}
\end{quote}

Now we are able to use the structural approach in order to count the number of undirected and directed circulant graphs of order 9.

\begin{quote}
\centering
\begin{tabular}{ll}
$f_{1}(1)=$&$\tilde{f}_{1}(1)=2$,\\
$f_{2}(1)=$&$\tilde{f}_{2}(1)=2^2$,\\
$f_{3}(1)=2^2$,&$\tilde{f}_{3}(1)=2^3$,\\
$f_{4}(1)=2^2$,&$\tilde{f}_{4}(1)=2^3$,\\
$f_{5}(1)=2^2$,&$\tilde{f}_{5}(1)=2^4$,\\
$f_{6}(1)=$&$\tilde{f}_{6}(1)=2^4$,\\
$f_{7}(1)=2^4$,&$\tilde{f}_{7}(1)=2^8$.\\
\end{tabular}
\end{quote}

Therefore
\begin{quote}
\centering
\begin{tabular}{l l}
$g_{1}(1)=2$,&$\tilde{g}_{1}(1)=2$,\\
$g_{2}(1)=2^{2}-2=2$,&$\tilde{g}_{2}(1)=2$,\\
$g_{3}(1)=\frac{1}{2}(2^2-2-2)0$,&$\tilde{g}_{3}(1)=\frac{1}{2}(2^3-2-2)=2$,\\
$g_{4}(1)=\frac{1}{2}(2^2-2-2)=0$,&$\tilde{g}_{4}(1)=\frac{1}{2}(2^3-2-2)=2$,\\
$g_{5}(1)=\frac{1}{4}(2^2-2-2)=0$,&$\tilde{g}_{5}(1)=\frac{1}{4}(2^4-2-2-4-4)=1$,\\
$g_{6}(1)=\frac{1}{3}(2^4-2-2)=4$,&$\tilde{g}_{6}(1)=4$,\\
$g_{7}(1)=\frac{1}{6}(2^4-2-2-3.4)=0$,&$\tilde{g}_{7}(1)=\frac{1}{6}(2^8-2-2-4-4-4-12)=38$.\\
\end{tabular}
\end{quote}

\begin{quote}
\centering
\begin{tabular}{l l}
$g(1)=g(9,1)=8$,&$\tilde{g}(1)=\tilde{g}(9,1)=51$,\\
\end{tabular}
\end{quote}

\subsection{The case $p^3$ for $p=3$}

We shall now use the structural approach to enumerate the \emph{undirected} circulant graphs of order 27. The number of such circulant graphs has already been determined by Brendan McKay and listed in \cite{KLP03}, but we shall here also obtain the generating function by degree for the number of these circulant graphs.

A list of symmetric Schur rings over $\mathbb{Z}_{27}$ was obtained using the package COCO (see \cite{faradzev&klin91}). (By a ``symmetric Schur-Ring", we mean one in which every basic set $T$ satisfies $-T=T$. This is sufficient for our purpose of enumerating undirected circulant graphs.) The following is the list.

\begin{equation*}
\begin{split}
\mathfrak{S}_{1}&=\langle \underline{0},\underline{1,26},\underline{2,25},\underline{3,24},\underline{4,23},\underline{5,22},\underline{6,21},\underline{7,20},\underline{8,19},\underline{9,18},\underline{10,17},\underline{11,16},\underline{12,15},\underline{13,14}\rangle,\\
\mathfrak{S}_{2}&=\langle \underline{0},\underline{1,26,8,19,10,17},\underline{2,25,7,20,11,16},\underline{3,24},\underline{4,23,5,22,13,14},\underline{6,21},\underline{9,18},\underline{12,15}\rangle,\\
\mathfrak{S}_{3}&=\langle \underline{0},\underline{1,26,2,25,4,23,5,22,7,20,8,19,10,17,11,16,13,14},\underline{3,24},\underline{6,21},\underline{9,18},\underline{12,15}\rangle,\\
\mathfrak{S}_{4}&=\langle \underline{0},\underline{1,26,8,19,10,17},\underline{2,25,7,20,11,16},\underline{3,24,6,21,12,15},\underline{4,23,5,22,13,14},\underline{9,18}\rangle,\\
\mathfrak{S}_{5}&=\langle \underline{0},\underline{1,26,2,25,4,23,5,22,7,20,8,19,10,17,11,16,13,14},\underline{3,24,6,21,12,15},\underline{9,18}\rangle,\\
\mathfrak{S}_{6}&=\langle \underline{0},\underline{1,26,2,25,4,23,5,22,7,20,8,19,10,17,11,16,13,14},\underline{3,24,6,21,9,18,12,15}\rangle,\\
\mathfrak{S}_{7}&=\langle \underline{0},\underline{1,26,2,25,3,24,4,23,5,22,6,21,7,20,8,19,10,17,11,16,12,15,13,14},\underline{9,18}\rangle,\\
\mathfrak{S}_{8}&=\langle \underline{0},\underline{1,2,3,4,5,6,7,8,9,10,11,12,13,14,15,16,17,18,19,20,21,22,23,24,25,26}\rangle\\
\end{split}
\end{equation*}

\begin{figure}[h!]
\centering
\includegraphics[width=0.3\linewidth]{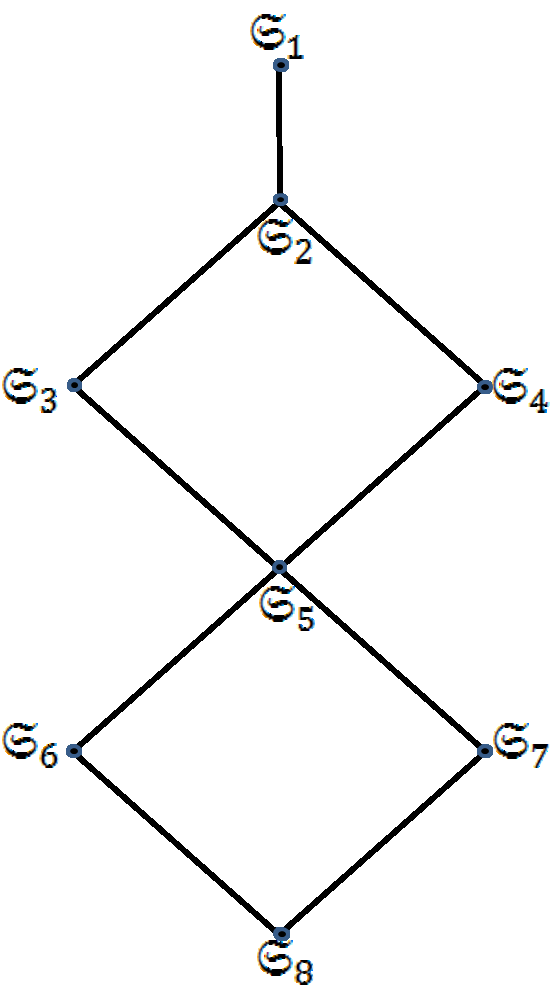}
\caption{Lattice of all $\mathcal{S}$-rings for $n=27$ Undirected}
\label{Fig7}
\end{figure}

In this list, we can observe that $\mathfrak{S}_{1}$ is the finest with the smallest automorphism group, while $\mathfrak{S}_{8}$ has the largest automorphism group. Therefore $\mathfrak{S}_{1}$ contains all the other Schur rings. We may now construct the lattice of Schur rings. This is given in Figure~\ref{Fig7}.
\\
\vspace{5mm}
Using Equation \eqref{eq:erl} 
we can obtain the generating functions $f_{i}(t)$. These are as follows
\begin{equation*}
\begin{split}
f_1(t)&=(1+t^2)^{13}\\
f_2(t)&=(1+t^6)^3(1+t^2)^4\\
f_3(t)&=(1+t^{18})(1+t^2)^4\\
f_4(t)&=(1+t^6)^4(1+t^2)\\
f_5(t)&=(1+t^{18})(1+t^6)(1+t^2)\\
f_6(t)&=(1+t^{18})(1+t^8)\\
f_7(t)&=(1+t^{24})(1+t^2)\\
f_8(t)&=(1+t^{26})
\end{split}
\end{equation*}
Table ~\ref{table:Table4} gives a list of the sizes of the automorphism groups and their normalizers. These were again obtained using GAP.
 
\begin{table}[h!]
\caption{Sizes of Automorphism Groups and their Normalizers for the Case n=27}
\begin{adjustwidth}{-1.2cm}{}
\begin{quote}
\centering
\begin{tabular}{|c|c|c|c|}
\hline
$G_i$&$|G_i|$&$|N(G_i)|$&$\frac{|G_i|}{|N(G_i)|}$\\
\hline
$G_1$&54&486&$1/9$\\
$G_2$&486&4374&$1/9$\\
$G_3$&34992&104976&$1/3$\\
$G_4$&181398528&544195584&$1/3$\\
$G_5$&13060694016&13060694016&1\\
$G_6$&286708355039232000&286708355039232000&1\\
$G_7$&3656994324480&3656994324480&1\\
$G_8$&10888869450418352160768000000&10888869450418352160768000000&1\\
\hline
\end{tabular}
\label{table:Table4}
\end{quote}
 \end{adjustwidth}
\end{table}
\vspace{8mm}

We may now determine $g_i=g_{i}(t)$ for $i=1,2,...,8$, using \eqref{eq:erl1} and Figure~\ref{Fig7}.
\begin{equation*}
\begin{split}
g_8&=f_8=1+t^{26}\\
g_7&=f_7-g_8=t^{24}+t^2\\
g_6&=f_6-g_8=t^{18}+t^8\\
g_5&=f_5-(g_8+g_7+g_6)=t^{20}+t^6\\
g_4&=\frac{1}{3}(f_4-g_8-g_7-g_6-g_5)=t^{20}+t^{18}+2t^{14}+2t^{12}+t^8+t^6\\
g_3&=\frac{1}{3}(f_3-g_8-g_7-g_6-g_5)=t^{24}+2t^{22}+t^{20}+t^6+2t^4+t^2\\
g_2&=\frac{1}{9}(f_2-g_8-g_7-g_6-g_5-3g_4-3g_3)=t^{18}+2t^{16}+t^{14}+t^{12}+2t^{10}+t^8\\
g_1&=\frac{1}{9}(f_1-g_8-g_7-g_6-g_5-3g_4-3g_3-9g_2)\\
&=t^{24}+8t^{22}+31t^{20}+78t^{18}+141t^{16}+189t^{14}+189t^{12}+141t^{10}+78t^8+31t^6+\\
& 8t^4+t^2\\
\end{split}
\end{equation*}
Therefore
\begin{equation*}
\begin{split}
g(t)=g_1+g_2+...+g_8=t^{26}+&3t^{24}+10t^{22}+34t^{20}+81t^{18}+143t^{16}+192t^{14}+\\
&192t^{12}+143t^{10}+81t^8+34t^6+10t^4+3t^2+1\\
\end{split}
\end{equation*}

This gives the same generating function as that obtained below using the multiplier method below. It confirms McKay's old result [unpublished, 1995] that there are 928 non-isomorphic, undirected circulant graphs on 27 vertices.


\section{The multiplier approach for $n=p^2$ when $p=3$} \label{sec:multiplier}


Since \adm's Conjecture does not hold for $n=p^2$ we need the next result which tells us, in terms of their connecting sets, when two circulant graphs of this order are isomorphic. This isomorphism criterion will require us to partition the elements of the connecting sets into \emph{layers}. This is done in the following way: We will first consider the set $\mathbb{Z}'_{p^2}= \mathbb{Z}'_{p^2}-\{0\}$ and divide its elements into two layers, namely $Y_0$ and $Y_1$, where $Y_0$ will contain those elements which do not have $p$ as a factor and $Y_1$ will contain those elements which do have $p$ as a factor. A connecting set $X$ is then given by
\[X=X_{0}\dot{\cup} X_{1}\]
where $X_{0}=X \cap Y_0$ and $X_{1}=X \cap Y_1$. The layer $X_{0}$ is a subset of $\mathbb{Z}^\ast_{p^2}$ while the layer $X_{1}$ is a subset of $p\mathbb{Z}^\ast_p$. In addition, when these layers are acted upon (multiplicatively) by elements of $\mathbb{Z}^\ast_{n}$, where in this case $n=p^2$, these layers are invariant. In \cite{KLP96} the following isomorphism criterion for circulant graphs of order $p^2$ was presented. 

\begin{Theorem}[\cite{KLP96}]
\label{thm:theorem2}
Two circulant graphs $\Gamma(\mathbb{Z}_n,X)$ and $\Gamma'=\Gamma(\mathbb{Z}_n,X')$ with $n=p^2$ vertices, are isomorphic if and only if their respective layers are multiplicatively equivalent, i.e.
\begin{equation}\tag{$M_{2}$}
X'_{0}=m_{0}X_{0},   X'_{1}=m_{1}X_{1},
\end{equation}
for a pair of multipliers $m_{0}, m_{1} \in \mathbb{Z}^\ast_{p^{2}}$. Moreover, in the above, one must have
\begin{equation}\tag{E}
m_{0}=m_{1}
\end{equation}
whenever
\begin{equation}\tag{R}
(1+p)X_{0}\not= X_{0}
\end{equation}
\end{Theorem}

We shall illustrate in some detail the use of this result for counting circulants of order 9, based on the treatment given in \cite{KLP96}, in order to introduce the use of the inclusion-exclusion principle and also some techniques and notation which will be expanded upon in the next section. Our detailed treatment should help to make the more difficult case of $n=p^3$ clearer.

First of all, in practice, it is easier to count orbits under invariance conditions
\begin{equation}\tag{$\neg R$}
(1+p)X_{0}=X_{0},
\end{equation}
that is, when the restrictions have an equality, rather than under the non-invariance condition $(R)$. Therefore, when the problem under consideration includes the non-invariance condition $(R)$, this is changed to the invariance condition $(\neg R)$ and then the result is subtracted from the total amount. 

Let us consider the case when $n=9$. In this case we have
\begin{equation*}
\begin{split}
\mathbb{Z}^\ast_{9}&=\lbrace 1,2,4,5,7,8\rbrace \mbox{ and }\\
\mathbb{Z}'_{9}&=\lbrace 1,2,3,4,5,6,7,8\rbrace
\end{split}
\end{equation*}
that is, the connecting set $X$ is a subset of $\mathbb{Z}'_{9}$  and the multipliers $m_0$ and $m_1$ come from $\mathbb{Z}^\ast_9$.
Let $Y_0$ and $Y_1$ be the two layers of $\mathbb{Z}'_{9}$.  Therefore $\mathbb{Z}'_9=Y_0 \dot{\cup}Y_1$ where
\begin{equation*}
\begin{split}
Y_0&=\lbrace 1,2,4,5,7,8\rbrace\\
Y_{1}&=\lbrace3,6\rbrace
\end{split}
\end{equation*}
The connecting sets are then given by the layers
\begin{equation*}
\begin{split}
X_{0}&=X\cap Y_0\\
X_{1}&=X\cap Y_1
\end{split}
\end{equation*}

Now two circulant graphs may be isomorphic either under one multiplier, that is, $m_0=m_1$, or two distinct multipliers, that is, when $m_0\neq m_1$. From Theorem ~\ref{thm:theorem2}, we have that the non-invariance condition $(R)$, holds only when the multipliers are equal. Therefore in order to count those circulant graphs which are equivalent via two different multipliers, we need to consider the invariance relation $(\neg R)$ given by $4X_{0}=X_{0}$. One must note however, that this relation may still hold when the multipliers are equal.

When enumerating under this invariance condition, the set $X_{0}$ must be taken from whole subsets of $Y_0$ which are invariant under $4Y_0=Y_0$ such that $4X_{0}=X_{0}$. This will give the partition of $Y_0$ as $Y^\ast_0$. In this case $Y^\ast_0=\lbrace \lbrace 1,4,7 \rbrace, \lbrace 2,8,5 \rbrace \rbrace$. Therefore under the condition $4X_{0}=X_{0}$, the set $X_{0}$ must be a union of these parts and the multiplicative action is on the sets $\lbrace 1,4,7 \rbrace$ and $\lbrace 2,8,5\rbrace$. That is, $X_{0}$ must either contain all of the set $\lbrace 1,4,7 \rbrace$, or none of it and similarly all of the set $\lbrace 2,8,5 \rbrace$ or none of it. 

In order to count the number of non-isomorphic directed circulant graphs on 9 vertices, we will divide the counting problem into two subproblems: determining the orders of $A_1(9)$ and $A_2(9)$, in which

\begin{itemize}
\item
$A_1(9)$ is the set of all those circulant digraphs which are distinct under the invariance
condition $(\neg R)$ with no restriction on the multipliers;
\item
$A_2(9)$ is the set of all those circulant digraphs which are distinct under the non-invariance condition $(R)$  and $m_0=m_1$. 
\end{itemize}

Let us first consider $A_1(9)$. In this case we have the invariance condition $4X_{0}=X_{0}$ and no restriction on the multipliers, that is, two circulant graphs may be equivalent under one multiplier or two distinct multipliers.

Here we need the action of $\mathbb{Z}^\ast_{9}\times \mathbb{Z}^\ast_{9}$ on $\lbrace Y_1 \cup Y^\ast_0 \rbrace$. The multiplier on $Y_1$ can be the same or different from that on $Y^\ast_0$. Therefore we have to consider the action of all $(i,j)\in \mathbb{Z}^\ast_{9} \times \mathbb{Z}^\ast_{9}$ on $\lbrace 3,6\rbrace \cup \lbrace F,H \rbrace$ where $F=\lbrace1,4,7\rbrace$ and $H=\lbrace2,8,5\rbrace$. For example $(2,5)$ has the action (3\hspace{2mm}6)(F\hspace{2mm}H), where 2 acts on $\lbrace3,6\rbrace$ and 5 acts on $\lbrace F,H\rbrace$. This gives the monomial $x_{1}^{2}x_{2}$. There are $6^2$ such actions since $|\mathbb{Z}^\ast_{9}|=6$.
We may however, determine $A_{1}(9)$ more simply, by finding the cycle index of $\mathbb{Z}^\ast_{9}$ on $\lbrace 3,6\rbrace$ and $\mathbb{Z}^\ast_{9}$ on $\lbrace F,H\rbrace$ and take the product. Both of these are equivalent to the action of $\mathbb{Z}^\ast_{9}$ on $\lbrace 1,2\rbrace$ mod 3, so we may simply find the latter and square.

\begin{quote}
\centering
\begin{tabular}{|c|c|c|}
\hline
Action of&Action&Cycle Structure\\
\hline
1&(1)(2)&$x_{1}^{2}$\\
2&(1,2)&$x_{2}$\\
4&(1)(2)&$x_{1}^{2}$\\
5&(1,2)&$x_{2}$\\
7&(1)(2)&$x_{1}^{2}$\\
8&(1,2)&$x_{2}$\\
\hline
\end{tabular}
\end{quote}
Therefore the cycle index is:

\[\frac{1}{6}(3x_{1}^{2}+3x_{2})\]
Squaring and simplifying gives the following generating function for the non-isomorphic circulant graphs in $A_2(9)$: 
\[\widetilde{A_2(9)}(\mathbf{x})=\frac{1}{4}(x_{1}^{2}+x_{2})^2,\]
where, as is usual in graph enumeration, substituting $x_{i}=2$ for all $i$ gives the size of the set in question. In this case, $|A_1(9)|=9$.

Let us now consider $A_2(9)$. Since in this case we have the non-invariance condition $(R)$, we shall make use of the counting procedure described previously. Let
\begin{equation*}
\begin{split}
A_{21}& \mbox{ be the set of all the circulants which are distinct under } m_0=m_1, \mbox{ that is, when }\\
&\mbox{ \'{A}d\'{a}m's condition holds and }\\
A_{22}& \mbox{ the set of of circulants which are distinct under } (\neg R) \mbox{ and having }\\
&m_0=m_1\\
\end{split}
\end{equation*}
Then
\[|A_2(9)|=|A_{21}|-|A_{22}|.\]

Let us first determine $A_{21}$. Here we need to count all circulant graphs assuming \'{A}d\'{a}m's conjecture holds, that is the number of one-multiplier equivalent directed circulant graphs. This is the number of orbits under the action $(\mathbb{Z}^\ast_{9},\mathbb{Z'}_{9})$.

\begin{quote}
\centering
\begin{tabular}{|c|c|c|}
\hline
Action of&Action&Cycle Structure\\
\hline
1&(1)(2)(3)(4)(5)(6)(7)(8)&$x_{1}^{8}$\\
2&(1,2,4,8,7,5)(3,6)&$x_{2}x_{6}$\\
4&(1,4,7)(2,8,5)(3)(6)&$x_{1}^{2}x_{3}^{2}$\\
5&(1,5,7,8,4,2)(3,6)&$x_{2}x_{6}$\\
7&(1,7,4)(2,5,8)(3)(6)&$x_{1}^{2}x_{3}^{2}$\\
8&(1,8)(2,7)(4,5)(3,6)&$x_{2}^{4}$\\
\hline
\end{tabular}
\end{quote}
Therefore the cycle index corresponding to $A_{21}$ is given by:
\[\widetilde{A_{21}}(\mathbf{x})=\frac{1}{6}(x_{1}^{8}+2x_{2}x_{6}+2x_{1}^{2}x_{3}^{2}+x_{2}^{4})\]
and substituting $x_{i}=2$ for all $i$ we obtain
\[|A_{21}|=\frac{1}{6}(2^{8}+2(2)(2)+2(2)^{2}(2)^{2}+(2)^4)=52\]

Let us now determine $|A_{22}|$. Since in $A_{22}$ the multipliers are equal and we have that $4X_0=X_0$, we need to consider the orbits of the action $(\mathbb{Z}^\ast_9,Y^\ast_0 \cup Y_1)$, that is, $(\mathbb{Z}^\ast_9,\lbrace\lbrace 1,4,7\rbrace,\lbrace 2,8,5 \rbrace, 3,6 \rbrace)$. The sets $\lbrace 1,4,7\rbrace$ and $\lbrace 2,8,5 \rbrace$ are blocks, that is, each block must appear whole as a neighbour or not at all. The contents of $X_1$ do not influence whether or not $X_0$ is 4-invariant, therefore 3 and 6 are acted upon separately. Now the action of $\mathbb{Z}^\ast_9$ on $\lbrace\lbrace 1,4,7\rbrace, \lbrace2,8,5\rbrace\rbrace$ as two points is equivalent to the action of $\mathbb{Z}^\ast_9$ on $\lbrace1,2\rbrace \bmod\ 3$. Therefore, this action may be expressed as $(\mathbb{Z}^\ast_{9},\lbrace1',2',1,2\rbrace)$ mod3 where 1', 2' represent $\lbrace1,4,7\rbrace$ and $\lbrace 2,8,5\rbrace$ respectively.

\begin{quote}
\centering
\begin{tabular}{|c|c|c|}
\hline
Action of&Action&Cycle Structure\\
\hline
1&$(1)(2)(1')(2')$&$x_{1}^{4}$\\
2&$(1,2)(1',2')$&$x_{2}^{2}$\\
4&$(1)(2)(1')(2')$&$x_{1}^{4}$\\
5&$(1,2)(1',2')$&$x_{2}^{2}$\\
7&$(1)(2)(1')(2')$&$x_{1}^{4}$\\
8&$(1,2)(1',2')$&$x_{2}^{2}$\\
\hline
\end{tabular}
\end{quote}

Therefore the cycle index corresponding to $A_{22}$ is given by:

\[\widetilde{A_{22}}(\mathbf{x})=\frac{1}{6}(3x_{1}^{4}+3x_{2}^{2})=\frac{1}{2}(x_{1}^{4}+x_{2}^{2})\]
and substituting $x_{i}=2$ for all $i$ we obtain

\[|A_{22}|=\frac{1}{2}(2^{4}+2^{2})=10\]
Therefore we have $|A_2(9)|=|A_{21}|-|A_{22}|=52-10=42$.

Therefore combining our results we obtain:

\[|A_1(9)|+|A_2(9)|=9+42=51.\]
This means that 51 directed, non-isomorphic circulant graphs on 9 vertices exist.

\vspace{3mm}
Summarizing the above in order to see more clearly the role of inclusion-exclusion, what we have essentially, is the set $A_{21}$ which contains the distinct circulant graphs under the conditions $(R)$ and $(\neg R)$ and $m_0=m_1$. The set $A_1$ contains distinct circulants under the condition $(\neg R)$ with no restriction on the multipliers (that is the multipliers could be the same or different). Now what we require is $|A_{21} \cup A_1|$. We know that $|A_{21} \cup A_1|=|A_{21}|+|A_1|-|A_{21} \cap A_1|$, where $|A_{21} \cap A_1|$ counts all those circulant graphs with $(\neg R)$ and $m_1=m_0$. This is simply $|A_{22}|$ as mentioned above.

\vspace{3mm}
The same procedure may be repeated for undirected circulants. However, as previously stated, we now need a slight modification on the connecting set,
that is, the connecting set is a subset $X$ of $\mathbb{Z}'_9$ which must have the property $X=-X$, and the multipliers $m_0$ and $m_1$ come from $\mathbb{Z}^\ast_9$.

Since the elements in the connecting sets are paired by inversion,  we partition $\mathbb{Z}'_9$ as
\[\mathbb{Z}'_9=\lbrace \lbrace1,8\rbrace,\lbrace2,7\rbrace,\lbrace3,6\rbrace,\lbrace4,5\rbrace\rbrace\]

One must note that any connecting set we shall work with must have either both elements of a given pair or none. The multiplicative action must therefore be taken on these pairs.

In this case we have
\begin{equation*}
\begin{split}
Y_0&=\lbrace\lbrace1,8\rbrace,\lbrace2,7\rbrace,\lbrace4,5\rbrace\rbrace\\
Y_1&=\lbrace{3,6}\rbrace
\end{split}
\end{equation*}
still partitioned into inverse pairs and
\[X_{0}=X \cap Y_0 \mbox{ and } X_{1}=X \cap Y_1.\]
We note here that $Y^\ast_0=\lbrace 1,2,4,5,7,8\rbrace$, that is the two separate blocks obtained previously in $Y^\ast_0$ are now merged together so that every element and its inverse are in the same set.

In order to determine $|A_{21}|$, again we require the orbits of the action $(\mathbb{Z}^\ast_{9},\mathbb{Z}'_{9})$ and therefore we will consider
\[(\mathbb{Z}^\ast_{9},\lbrace \lbrace1,8\rbrace,\lbrace2,7\rbrace,\lbrace3,6\rbrace,\lbrace4,5\rbrace\rbrace).\]
Again $\lbrace1,8\rbrace,\lbrace2,7\rbrace,\lbrace3,6\rbrace,\lbrace4,5\rbrace$ are considered as blocks, that is, we may consider $(\mathbb{Z}^\ast_{9},\lbrace K,L,M,N\rbrace)$ where $K=\lbrace1,8\rbrace$, $L=\lbrace2,7\rbrace$, $M=\lbrace3,6\rbrace$, $N=\lbrace4,5\rbrace$. The cycle index corresponding to this action is given by
\[\widetilde{A_{21}}(\mathbf{x})=\frac{1}{6}(2x_{1}^{4}+4x_{3}x_{1})\]
and substituting $x_{i}=2$ for all $i$, we obtain $A_{11}=8$.

To determine $|A_{22}|$  we again require the action $(\mathbb{Z}^\ast_{9},Y^\ast_0\cup Y_1)$. Therefore we have
\[(\mathbb{Z}^\ast_{9},\lbrace \lbrace1,2,4,5,7,8\rbrace,\lbrace3,6\rbrace \rbrace)\]
This action may be seen as $(\mathbb{Z}^\ast_{9},\lbrace M,O \rbrace)$ where $O=\lbrace 1,2,4,5,7,8\rbrace$. The cycle index corresponding to this action is
\[\widetilde{A_{22}}(\mathbf{x})=\frac{1}{6}(6)x_{1}^{2}=x_{1}^{2}\]
and substituting $x_{i}=2$ we obtain $A_{22}=4$.

Therefore
\[|A_2(9)|=|A_{21}|-|A_{22}|=4.\]

For $A_1(9)$ we need to multiply the cycle indices of the actions $(\mathbb{Z}^\ast_{9},\lbrace3,6\rbrace)=(\mathbb{Z}^\ast_{9},\lbrace M\rbrace)$ and $(\mathbb{Z}^\ast_{9},\lbrace 1,2,4,5,7,8\rbrace)=(\mathbb{Z}^\ast_{9},\lbrace O\rbrace)$. Both these cycle indices are equivalent to $x_{1}$. Therefore

\[\widetilde{A_1(9)}(\mathbf{x})=x_{1}^{2}.\]
Substituting, we obtain $|A_1(9)|=4$.

Therefore the number of non-isomorphic, undirected circulant graphs on 9 vertices is

\[|A_1(9)|+|A_2(9)|=4+4=8.\]

\section{The multiplier approach for $n=p^3$}

\subsection{The Main Isomorphism Theorem}

We shall first state a general isomorphism theorem for circulant graphs which was proved by Klin and P\"{o}schel in \cite{KP80}.  
Theorem \ref{thm:theorem2} which we used for the enumeration of circulants of order $p^2$ is a special case of this. We then state the special case of the result of Klin and P\"{o}schel for order $p^3$, which will be our main tool. 

\begin{Theorem}[\cite{KP80},\cite{LP2000}]
\label{thm:theorem3}
Let $n=p^k$ ($p$ an odd prime) and let $\Gamma$ and $\Gamma'$ be two $p^k$-circulants with the connecting sets $X$ and $X'$, respectively. Then $\Gamma$ and $\Gamma'$ are isomorphic if and only if their respective layers are multiplicatively equivalent, that is,
\begin{equation}\tag{$M_{k}$}
X'_{i}=m_{i}X_{i}, \hspace{5mm} i=0,1,...,k-1,
\end{equation}
for an  arbitrary set of multipliers $m_0, m_1,...m_{k-1}\in \mathbb{Z}^\ast_p$ which satisfy the following constraints: whenever the layer $X_{i}$ satisfies the non-invariance condition

\begin{equation}\tag{$R_{ij}$}
(1+p^{k-i-j-1})X_{i}\neq X_{i}
\end{equation}
for some $i \in \lbrace 0,1,...,k-2\rbrace$ and $j \in \lbrace 0,1,...,k-2-i\rbrace$, the successive multipliers $m_{i},...,m_{k-j-1}$ meet the system of congruences
\begin{equation}\tag{$E_{ij}$}
\begin{split}
m_{i+1}&\equiv m_i(\bmod\ \ p^{k-i-j-1}),\\
m_{i+2}&\equiv m_{i+1}(\bmod\  p^{k-i-j-2}),\\
&\vdots\\
m_{k-j-1}&\equiv m_{k-j-2}(\bmod\  p).
\end{split}
\end{equation}
\end{Theorem}

For the case when $k=3$, Theorem~\ref{thm:theorem3} translates to the following theorem which we shall refer to as the Main Isomorphism Theorem

\begin{Theorem}[Main Isomorphism Theorem]
\label{thm:theorem4}
Let $n=p^3$ ($p$ an odd prime) and let $\Gamma$ and $\Gamma'$ be two $p^3$-circulants with the connecting sets $X$ and $X'$, respectively. Then $\Gamma$ and $\Gamma'$ are isomorphic if and only if their respective layers are multiplicatively equivalent, that is,
\begin{equation}\tag{$M_{3}$}
X'_{0}=m_{0}X_{0},\hspace{3mm}X'_{1}=m_{1}X_{1}\hspace{3mm}X'_{2}=m_{2}X_{2},
\end{equation}
for an  arbitrary set of multipliers $m_0, m_1, m_2 \in \mathbb{Z}^\ast_{p^3}$.  Moreover, in the above, one must have

\begin{itemize}
\item[(i)] \hfill \makebox[0pt][r]{%
            \begin{minipage}[b]{\textwidth}
              \begin{equation}
                 m_1 \equiv m_0(\bmod\  p^2) \mbox{ and } m_2 \equiv m_1(\bmod\  p)\tag{$E_{00}$}
              \end{equation}
          \end{minipage}}
          \end{itemize}
whenever
\begin{equation}\tag{$R_{00}$}
(1+p^2)X_{0}\neq X_{0},
\end{equation}

\begin{itemize}
\item[(ii)] \hfill \makebox[0pt][r]{%
            \begin{minipage}[b]{\textwidth}
              \begin{equation}
                 m_1 \equiv m_0(\bmod\  p)\tag{$E_{01}$}
              \end{equation}
          \end{minipage}}
          \end{itemize}
whenever
\begin{equation}\tag{$R_{01}$}
(1+p)X_{0}\neq X_{0},
\end{equation}

\begin{itemize}
\item[(iii)] \hfill \makebox[0pt][r]{%
            \begin{minipage}[b]{\textwidth}
              \begin{equation}
                 m_2 \equiv m_1(\bmod\  p)\tag{$E_{10}$}
              \end{equation}
          \end{minipage}}
          \end{itemize}
whenever
\begin{equation}\tag{$R_{10}$}
(1+p)X_{1}\neq X_{1}.
\end{equation}
\end{Theorem}

Whereas Theorem~\ref{thm:theorem2} for $p^2$ involved two multipliers and two subcases depending on non-invariance conditions on the layers, the isomorphism theorem for $p^k$ involves $k$ multipliers and $\binom{k}{2}$ cases coming from the non-invariance conditions $R_{ij}$, making it more difficult to apply in practice for enumeration purposes. And what makes the enumeration problem particularly difficult is not only that there are multipliers for the separate layers of the connecting sets, but that, depending on non-invariance conditions, some multipliers must be equal in certain cases. Moreover, the intersection between the conditions makes this case even more difficult. The two cases for $p^2$ involved three different enumeration problems, as we have seen, and for $p^3$, the three non-invariance relations $R_{00}, R_{01}, R_{10}$ below, will break up into five cases which will eventually give eleven enumeration subproblems, as we shall see below.


\subsection{Representation and computational implementation of the Main Isomorphism Theorem, $p=3,5$}

 Liskovets and P\"{o}schel in \cite{LP2000} manage, for $n=p^3$, to partition the conditions of the Main Theorem into five parts which makes their use in enumeration much easier. These authors take into consideration all combinations of non-invariance conditions $(R_{ij})$, together with the remaining invariance conditions
 
 \begin{equation}\tag{$\neg R_{ij}$}
 (1+p^{k-i-j-1})X_{i}=X_{i}
 \end{equation}
 and make use of a number of results, in order to obtain the subproblem list for counting circulants of order $p^k$, $k\leq 4$. For details of how this list has been generated from the Main Theorem using results from number theory and walks through a rectangular lattice, the reader is referred to \cite{LP2000}. The necessary information required for the enumeration of $p^3$ circulants is listed in Table~\ref{table:Table3}, which we therefore take to be a rewording of Theorem~\ref{thm:theorem4}. This has been obtained from Table 1 in \cite{LP2000}.
 
 \begin{table}[ht]
 \begin{adjustwidth}{-1cm}{}
 \caption{The conditions for isomorphism of circulants of order $p^3$}
 \begin{quote}
 \centering
 \begin{tabular}{|c|c|c|c|}
 \hline
 Subproblem&Non-Invariance&Invariance&Conditions on\\
 &Conditions&Conditions&Multipliers\\
 \hline
 $A_1$&$\emptyset$&$\neg R_{01}$, $\neg R_{10}$&no restriction\\
 $A_2$&$R_{00}$&$\emptyset$&$m_2=m_1=m_0$\\
 $A_3$&$R_{01}$&$\neg R_{00}$, $\neg R_{10}$&$m_1=m_0$\\
 $A_4$&$R_{10}$&$\neg R_{01}$&$m_2=m_1$\\
 $A_5$&$R_{01},R_{10}$&$\neg R_{00}$&$m_2=m_1$ and $m_1 \equiv m_0(\bmod\  p)$\\
 \hline
 \end{tabular}
 \label{table:Table3}
 \end{quote}
 \end{adjustwidth}
 \end{table}

The five subcases $A_1$ to $A_5$ shown in Table~\ref{table:Table3}, give conditions on the three multipliers $m_0$, $m_1$ and $m_2$ for two circulants of order $p^3$ to be isomorphic. The re-interpretation of the Main Isomorphism Theorem by Liskovets and P\"{o}schel says that the three multipliers must satisfy at least one of the five sets of conditions for two ciruclants of order $p^3$ to be isomorphic. So, for example, condition  $A_3$ means that if the non-invariance condition $R_{01}$ holds, together with the invariance conditions $\neg R_{00}$ and $\neg R_{10}$, then $m_1=m_0$ but $m_2$ can be independent. 

If we let $A$ denote set of all non-isomorphic circulants of order $p^3$  and let $A_i$, for $i=1,\ldots,5$, also denote the set of circulants which are non-isomorphic under the respective condition of Table \ref{table:Table3}, then, the result of Liskovets and P\"{o}schel says that   
\[ A = A_1\cup A_2\cup A_3\cup A_4\cup A_5. \]

In addition to this information given in Table \ref{table:Table3}, we shall also use, without explicit mention, the following observations \cite{LP2000}.
\begin{equation*}
(R_{ij})\Rightarrow(R_{ij'}) \mbox{ whenever } j'\geq j
\end{equation*}
Therefore
\begin{equation*}
\neg(R_{ij'})\Rightarrow \neg(R_{ij})
\end{equation*}
As a result we have that
\begin{equation}
\neg(R_{01})\Rightarrow \neg(R_{00})
\end{equation}
In addition,
\begin{equation*}
(E_{ij})\Rightarrow(E_{i'j}) \mbox{ whenever } i'\geq i
\end{equation*}
Therefore
\begin{equation*}
\neg(E_{i'j})\Rightarrow \neg(E_{ij})
\end{equation*}
that is
\begin{equation}
\neg(E_{10})\Rightarrow \neg(E_{00})
\end{equation}
 

As we explained before, when the subproblem in question includes one or more non-invariance conditions, these are changed to invariance conditions and then the result is subtracted from the total number. Therefore, in order to count under a given non-invariance relation $R_{ij}$, we first
\begin{itemize}
\item[(i)]Determine the count under the action assuming the invariance relation $\neg (R_{ij})$,
\item[(ii)] Determine the count under the action without any (non)-invariance relations,
\item[(iii)] Subtract the result of (i) from (ii).
\end{itemize}
This procedure is often complicated by having both invariance and non-invariance conditions. For example to count the number of non-isomorphic circulants in case $A_3$ of Table~\ref{table:Table3}, we
\begin{itemize}
\item[(i)]First count under the conditions $m_1=m_0$, $\neg (R_{00}), \neg (R_{10})$.
\item[(ii)]Then count under the conditions $m_1=m_0$, $\neg (R_{01}),\neg (R_{00}),\neg (R_{10})$.
\item[(iii)] Then subtract the result of (ii) from (i).
\end{itemize}

Having counted the number of non-isomorphic circulants under each of the five isomorphism conditions we then need to calculate $|A|$ and therefore we would need to consider the intersections between the $A_i$. It however transpires that these intersections are empty. This can be seen by considering the invariance and non-invariance relations. It therefore follows that 
\[ |A| = |A_1|+ |A_2| + |A_3| + |A_4| +|A_5| \qquad (*) \]
This is the main reason why it is easier to do enumeration using the formulation of Theorem~\ref{thm:theorem4} as in Table~\ref{table:Table3}.

We shall now explain how these five problems lead to eleven subproblems using the case of undirected circulant graphs of order $p^3$ for $p=3$. 
 Roughly speaking, our goal is to solve each concrete subproblem by considering a suitable group of multipliers acting on a suitable combination of layers and applying the standard enumeration technique of P\'{o}lya-Redfield.
The actual values given by these subproblems and their generating functions will be given in the next two subsections for $p=3$ and $p=5$, respectively. 

So, let $\mathbb{Z}^*_{27}$ denote the set of units in $\mathbb{Z}_{27}$ and $\mathbb{Z}'_{27}$ the set $\mathbb{Z}_{27}-\{0\}$. Therefore the connecting set of the circulant graph is a subset $X$ of 
$\mathbb{Z}'_{27}$ and the multipliers $m_0,m_1,m_2$ come from $\mathbb{Z}^\ast_{27}$.

Let $Y_0$, $Y_1$, $Y_2$, be the three layers of $\mathbb{Z}'_{27}$ where $Y_0$ contains all elements of $\mathbb{Z}'_{27}$ which are relatively prime to 27, $Y_1$ contains those elements which are divisible by 3 and not by 9 and $Y_2$ contains those divisible by 9. 

Now by Theorem~\ref{thm:theorem4}, the non-invariance conditions in this case are:
\begin{equation*}
\begin{split}
R_{00}&: \hspace{3mm} 10X_{0}\neq X_{0}\\
R_{01}&: \hspace{3mm} 4X_{0}\neq X_{0}\\
R_{10}&: \hspace{3mm} 4X_{1}\neq X_{1},
\end{split}
\end{equation*}
where $X_{0}=X \cap \mathbb{Z}^\ast_{27}$ that is, $X_{0}=X \cap Y_0$ and $X_{1}=X \cap 3\mathbb{Z}^\ast_{9}$, that is, $X_{1}=X \cap Y_1$. Recall that $X_{2}=X\cap Y_2$.

As described in the $p^2$ case, when we enumerate under an invariance condition, such as $10X_{0}=X_{0}$, we must take $X_{0}$ from whole subsets of $Y_0$ which are invariant under $10Y_0=Y_0$. These subsets partition $Y_0$, therefore under the condition $10X_{0}=X_{0}$, the set $X_{0}$ must be a union of these parts. Therefore the multiplicative action is taken on these parts or blocks.

We shall denote the partitioned set corresponding to the invariance condition $10Y_{0}=Y_{0}$ by $Y^\ast_0$, that corresponding to $4Y_{0}=Y_{0}$ by $Y^{\ast\ast}_0$, and that corresponding to the invariance condition $4Y_{1}=Y_{1}$ by $Y^\ast_1$.
We have that
\begin{equation*}
\begin{split}
Y^\ast_0&=\lbrace\lbrace 1,10,19\rbrace,\lbrace 2,20,11\rbrace,\lbrace 4,13,22\rbrace,\lbrace 5,23,14\rbrace,\lbrace 7,16,25\rbrace,\lbrace 8,26,17\rbrace\rbrace\\
Y^{\ast\ast}_0&=\lbrace\lbrace 1,4,16,10,13,25,19,22,7\rbrace,\lbrace2,8,5,20,26,23,11,17,14\rbrace\rbrace\\
Y^\ast_1&=\lbrace\lbrace3,12,21\rbrace,\lbrace6,24,15\rbrace\rbrace
\end{split}
\end{equation*}

\medskip\noindent
Now consider first the subproblem $A_1$. In this case we have that, when the invariance conditions $\neg R_{01}$ and $\neg R_{10}$ hold, there is no restriction on the multipliers. Also, $A_1$ does not include any non-invariance condition, making this subproblem easier because it does not split into further subproblems. We have here an action on the layers arising from
\begin{equation*}
\begin{split}
4X_{0}&=X_{0} \mbox{ and }\\
4X_{1}&=X_{1}
\end{split}
\end{equation*}
Under the invariance condition $4X_{0}=X_{0}$, we must take $X_{0}$ from whole subsets of $Y_0$ which are invariant under $4Y_0=Y_0$. Therefore the set $X_{0}$ must be a union of the layers in $Y^{\ast\ast}_0$. Consequently, instead of $Y_0$ we shall make use of $Y^{\ast\ast}_0$. Similarly, under the invariance condition $4X_{1}=X_{1}$, we must take $X_{1}$ from whole subsets of $Y_1$ which are invariant under $4Y_{1}=Y_{1}$. These subsets partition $Y_{1}$ as $Y^\ast_1$. Therefore under the condition $4X_{1}=X_{1}$, the set $X_{1}$ must be a union of these parts. We shall therefore use $Y^\ast_1$ instead of $Y_1$.

Since we have no restriction on the multipliers here, cycle index for computing the size of $A_1$ is the product of the cycle indices $I_{(\mathbb{Z}^\ast_{27},Y^{\ast\ast}_0)}$, $I_{(\mathbb{Z}^\ast_{27},Y^{\ast}_1)}$ and $I_{(\mathbb{Z}^\ast_{27},Y_2)}$.

Let us now consider $A_2$. Here we have the condition that when $m_0=m_1=m_2$ then the non-invariance condition $R_{00}$ must hold. Since in this case we need to consider the non-invariance condition $R_{00}$, we will use the counting procedure described at the beginning of this section. Therefore this isomorphism condition will be split into the following two problems:
\begin{equation*}
\begin{split}
A_{21}:& \mbox{ The set of the action } (\mathbb{Z}^\ast_{27},\mathbb{Z}'_{27})\\
A_{22}:& \mbox{ The set resulting from an action with } \neg R_{00} \mbox{ that is with } 10X_{0}=X_{0}.
\end{split}
\end{equation*}
Again, since we have the condition $10X_{0}=X_{0}$ in $A_{22}$, we must take $X_{0}$ from whole subsets of $Y^\ast_0$. Since $m_0=m_1=m_2$, the action corresponding to $A_{22}$ is therefore $(\mathbb{Z}^\ast_{27},Y^\ast_0 \cup Y_1 \cup Y_2)$. The required result for $|A_2|$ is then given by $|A_{21}|-|A_{22}|$.

Let us now consider $A_3$. Here we have the conditions $R_{01}$, $\neg R_{00}$, $\neg R_{10}$ for $m_1=m_0$. Therefore in this case our isomorphism problem will again be divided into two problems, namely $A_{31}$ and $A_{32}$, where
$A_{31}$ is the set of non-isomorphic circulants resulting from the action with $\neg R_{00}$ and $\neg R_{10}$, that is, with layers arising from

\begin{equation*}
\begin{split}
10X_{0}&=X_{0} \mbox{ and }\\
4X_{1}&=X_{1}
\end{split}
\end{equation*}
and therefore we shall need to use $Y^{\ast}_0$ instead of $Y_0$ and $Y^\ast_1$ instead of $Y_1$, and
$A_{32}$ is the set of non-isomorphic circulants resulting from the action with $\neg R_{00}$, $\neg R_{10}$ and $\neg R_{01}$, that is, with blocks arising from
\begin{equation*}
\begin{split}
10X_{0}&=X_{0} \mbox{ and }\\
4X_{1}&=X_{1} \mbox{ and }\\
4X_{0}&=X_{0}.
\end{split}
\end{equation*}
In this case however, we know that $\neg(R_{01})\Rightarrow \neg(R_{00})$, therefore the first equation is redundant. Therefore for $A_{32}$ we shall use $Y^{\ast\ast}_0$ instead of $Y_0$ and $Y^\ast_1$ instead of $Y_1$.
Since $m_1=m_0$ and $m_2$ is independent, the cycle index of our action here, is the product of the cycle indices $I_{(\mathbb{Z}^\ast_{27},Y_0\cup Y_1)}$ and $I_{(\mathbb{Z}^\ast_{27},Y_2)}$ where, in both cases, the group $\mathbb{Z}^\ast_{27}$ is considered to act on the layers of the respective sets. This means that the cycle indices of $A_{31}$ and $A_{32}$ are:\\
\begin{equation*}
\begin{split}
\widetilde{A_{31}}(\mathbf{x}) = &I_{(\mathbb{Z}^\ast_{27},Y^\ast_0\cup Y^\ast_1)}\times I_{(\mathbb{Z}^\ast_{27},Y_2)}\\
\widetilde{A_{32}}(\mathbf{x}) = &I_{(\mathbb{Z}^\ast_{27},Y^{\ast\ast}_0\cup Y^\ast_1)}\times I_{(\mathbb{Z}^\ast_{27},Y_2)}\\
\end{split}
\end{equation*}
and $|A_3|=|A_{31}|-|A_{32}|$.

We shall now consider $A_4$. Here we have the conditions $R_{10}$ and $\neg R_{01}$ when $m_2=m_1$. Therefore we will now consider $A_{41}$ and $A_{42}$ as follows:

$A_{41}$ is the set of non-isomorphic circulants resulting from the action with $\neg R_{01}$, that is with layers arising from $4X_{0}=X_{0}$. Therefore in this action $X_{0}$ must be a union of parts in $Y^{\ast\ast}_0$, therefore we shall use $Y^{\ast\ast}_0$ instead of $Y_0$.

$A_{42}$ is the set of non-isomorphic circulants resulting from  the action with $\neg R_{01}$ and$\neg R_{10}$. This means the set $X_{0}$ must be a union of the
layers in $Y^{\ast\ast}_0$ and $X_{1}$ a union of layers in $Y^\ast_1$.
The required result for $|A_4|$ will then be $|A_{41}|-|A_{42}|$. This example also gives us an opportunity to illustrate our direct use of the cycle indices of the relevant group of multipliers.
Since $m_2=m_1$ while $m_0$ is independent, we require
\[I_{(\mathbb{Z}^\ast_{27},Y_1\cup Y_2)}\times I_{(\mathbb{Z}^\ast_{27},Y_0)},\]
blocked as required. Therefore we have these cycle indices:
\begin{equation*}
\begin{split}
\widetilde{A_{41}}(\mathbf{x}) = &I_{(\mathbb{Z}^\ast_{27},Y_1\cup Y_2)}\times I_{(\mathbb{Z}^\ast_{27},Y^{\ast\ast}_0)}\\
\widetilde{A_{42}}(\mathbf{x}) = &I_{(\mathbb{Z}^\ast_{27},Y^{\ast}_1\cup Y_2)}\times I_{(\mathbb{Z}^\ast_{27},Y^{\ast\ast}_0)}.\\
\end{split}
\end{equation*}

Finally, we consider $A_5$. Once again, $A_5$ will be divided into the problems $A_{51}$ and $A_{52}$. Although $A_{51}$ is determined in a manner similar to the previous cases, one should be cautious when determining $A_{52}$, since this time we have two non-invariance conditions. This means that we have to consider the following:

$A_{51}$ is the set of non--isomorphic circulants resulting from the action with $\neg R_{00}$, that is, we shall use $Y^\ast_0$ instead of $Y_0$,

$A_{52}$ is the set of non-isomorphic circulants resulting from the action with $\neg(R_{01} \mbox {and } R_{10} )$ and $\neg R_{00}$.

Now for $A_{52}$, by de Morgan's laws, we have that:
\begin{equation*}
\begin{split}
\neg(R_{01} \mbox { and } R_{10})\mbox{ and } \neg R_{00}=&(\neg R_{01} \mbox { or } \neg R_{10}) \mbox {and } \neg R_{00}\\
=&(\neg R_{01}\mbox { and } \neg R_{00}) \mbox { or } (\neg R_{10}\mbox { and } \neg R_{00})\\
|\neg(R_{01} \mbox { and } R_{10})\mbox{ and } \neg R_{00}|=&|(\neg R_{01}\mbox { and } \neg R_{00})|+|(\neg R_{10}\mbox { and } \neg R_{00})|-\\
&|(\neg R_{01}\mbox { and } \neg R_{10} \mbox{ and } \neg R_{00})|\\
\end{split}
\end{equation*}
Now $\neg(R_{01})\Rightarrow \neg(R_{00})$, therefore for $A_{52}$ we have:
\begin{equation*}
\begin{split}
|\neg(R_{01} \mbox { and } R_{10})\mbox{ and } \neg R_{00}|=&|\neg (R_{01})|+|(\neg (R_{10})\mbox { and } \neg (R_{00}))|-\\
&|(\neg (R_{01})\mbox { and } \neg (R_{10}))|.\\
\end{split}
\end{equation*}
Therefore we shall split $A_{52}$ into 3 enumeration subproblems, with the first problem enumerating under the condition $\neg (R_{01})$, the second under $\neg (R_{10})\mbox { and } \neg (R_{00})$ and the last under the invariance conditions $\neg (R_{01})\mbox { and } \neg (R_{10})$. These are the subproblems $A_{521}, A_{522}, A_{523}$, respectively, giving that 
\[ |A_5| = |A_{51}|-|A_{521}|-|A_{522}|+|A_{523}|.\]

Note that we now have two restrictions on the multipliers, an equality with $m_2=m_1$ and a congruence, $m_1\equiv m_0\bmod\  3$. Therefore in this case, we need to define a group $G$ which will act on $\lbrace Y_1\cup Y_2\cup Y_3 \rbrace$ (blocked as required according to the given invariance condition), such that:
\begin{enumerate}
\item[(1)] The same multiplier acts on all $\lbrace Y_1\cup Y_2\cup Y_3 \rbrace$\\
\item[(2)] Two different multipliers act:
\begin{itemize}
\item $a$ on $Y_1 \cup Y_2$
\item $a'$ on $Y_0$
\end{itemize}
with $a'\equiv a\bmod\ 3$.
\end{enumerate}
Now, since $a'\equiv a\bmod\ 3$ implies the possibility that $a' \equiv a$, the second possibility includes the first. Therefore we will construct $G$ as follows:
$G$ will contain all ordered pairs $(a,a')$ such that $a, a' \in \mathbb{Z}^\ast_{27}$ and $a'\equiv a\bmod\ 3$. Then $(a,a')$ will act as follows:
\[(a,a')(y)=\left\{
\begin{array}{ll}
&ay \mbox{ if } y \in Y_1 \cup Y_2\\
&a'y \mbox{ if } y \in Y_0
\end{array}\right.\]
once again, with $Y_0$, $Y_1$ and $Y_2$ blocked as required. This would give all actions as in (1) and (2) above. One may verify that $G$ is in fact a group since if $(a,a'),(b,b')\in G$ then $(ab,a'b')\in G$.

Finally we have the general 11-term formula for $|A|$ in terms of its subproblems:
{\small
\[
|A| = |A_1| + |A_{21}|-|A_{22}| + |A_{31}|- |A_{32}| + |A_{41}|- |A_{42}| +|A_{51}|-|A_{521}|-|A_{522}|+|A_{523}|. \qquad (**) \] }
We observe that (**) consists of \adm's term $A_{21}$, the enumerator up to Cayley isomorphism, together with ten less obvious small correcting terms. 

The reason why we have eleven subproblems arising from the five terms in $(*)$ following Table \ref{table:Table3} is because $11=1\cdot1+3\cdot2+1\cdot4$, which gives the 3rd (little) Schr\"{o}der number (see, for example, sequence A001003 in \cite{OEIS}). Note that, in this context, 5 is the third Catalan number.

\medskip\noindent
The above analysis can be carried out in an analogous way for $p=5$, and for directed or undirected graphs. To generalise a bit our notation, let us define $A[s;p^3]$ for $s=u$, undirected, $s=d$, directed, and $p=3,5$ to be the number of undirected/directed circulant graphs on $p^3$ vertices. We similarly define $A_i[s;p^3], A_{ij}[s;p^3]$ and $A_{ijk}[s;p^3]$ to be the number of undirected/directed circulant graphs on $p^3$ vertices making up the corresponding intermediate terms $|A_i|, |A_{ij}|$ or $|A_{ijk}|$, respectively.

With a slight abuse of notation, we also let $A_{\mathbf{w}}(t):=A_{\mathbf{w}}[s;p^3](t)$ (where ${\mathbf{w}}$ represents the subscript $i$, $ij$ or $ijk$) denote generically the generating function by valency of the type of circulant graph under
enumeration in the set $A_{\mathbf{w}}$ (intermediate or otherwise) being considered, where the coefficient of $t^r$ equals the number of circulant graphs under discussion having valency (out-valency, in the directed case) equal to $r$. This coefficient is sometimes denoted by $A_{\mathbf{w}}[s;p^3,r]$. Therefore the terms defined above, counting the circulant graphs regardless of valency, are equal to $A_{\mathbf{w}}[s;p^3](1)$.


\subsection{Numerical results for undirected and directed circulants} 

We first give our main results in Table \ref{tbl:allvalues} which shows the number of circulant (di)graphs on 27 and 125 vertices.

\bigskip
\tabcolsep=0.65ex
\begin{table}[htb!]
\caption{\label{tbl:allvalues} The number of undirected and directed $p^3$-circulant graphs, $p=3,5$}
\begin{tabular}{|l|r|r|}
\hline
        &  Undirected&    Directed\\
\hline
\hline
n=27    &  928                    & 3,728,891     \\
n=125    &  92,233,720,411,499,283 & 212,676,479,325,586,539,710,725,989,876,778,596            \\
\hline
\end{tabular}
\end{table}

Next, Table \ref{tbl:intermediatevalues} gives the values of some of the intermediate terms which, as defined and described in the previous section, jointly with (*) and (**) together yield the values in Table \ref{tbl:allvalues}. We note that it is \adm's term $A_{21}$ which is the greatest contributor to these values.

\bigskip
\tabcolsep=0.65ex
\begin{table}[htb!]
\caption{\label{tbl:intermediatevalues} The number of undirected and directed $p^3$-circulant graphs, $p=3,5$:
intermediate contributors and totals} 
\hspace*{-3ex}
\begin{tabular}{|l|r|r|r|r|}
\hline
Term&                  Undir&            Undir&     Dir&         Dir\\
        &                   n=27&            n=125&    n=27&       n=125\\
\hline
\hline
$A_1$&                         8&               27&      27&         216\\
\hline
$A_{21}$&                    944&92233720411833168& 3730584&         ($\flat$)\\
$A_{22}$&                     48&           419664&    2776&879609512976\\
$A_2=A_{21}-A_{22}$&         896&92233720411413504& 3727808&        ($\natural$)\\
\hline
$A_{31}$&                    16&              1272&     156&     5034768\\
$A_{32}$&                     8&                30&      30&         420\\
$A_3=A_{31}-A_{32}$&          8&              1242&     126&     5034348\\
\hline
$A_{41}$&                    16&              1272&     156&     5034768\\
$A_{42}$&                     8&                30&      30&         420\\
$A_4=A_{41}-A_{42}$&          8&              1242&     126&     5034348\\
\hline
$A_{51}$&                    32&             86592&    1168&175943379264\\
$A_{521}$&                   16&              1680&     200&    13423440\\
$A_{522}$&                   16&              1680&     200&    13423440\\
$A_{523}$&                    8&                36&      36&        1044\\
{\small ($A_{52}=A_{521}+A_{522}-A_{523}$)}&    24&    3324&364&26845836\\
{\small $A_5\!=\!A_{51}\!-\!A_{521}\!-\!A_{522}\!+\!A_{523}$}&8&83268&804&175916533428\\
\hline
\hline
$A\!=\!A_1\!+\!A_2\!+\!A_3\!+\!A_4\!+\!A_5$&928&92233720411499283&3728891&($\sharp$)\\
\hline
\end{tabular}

\medskip
($\flat$) 212676479325586539710726693559689232\\ 
($\natural$) 212676479325586539710725813950176256\\
($\sharp$) 212676479325586539710725989876778596
\end{table}

Finally, 
we give all the final generating functions $A[u;27](t), A[d;27](t), A[u;125](t)$ and $A[d;125](t)$, that is, the generating functions for all the circulants of \mbox{orders} $27$ and $125$, undirected and directed.

\bigskip
{\small
\hspace*{-10ex}
\begin{tabular}{l>{\raggedright\arraybackslash$}p{12cm}<{$}}
$A[u;27](t)=$ & 
t^{26}+3t^{24}+10t^{22}+34t^{20}+81t^{18}+143t^{16}+192t^{14}+192t^{12}+143t^{10}+81t^8+34t^6+10t^4+3t^2+1\\
\end{tabular}

\medskip
\hspace*{-10ex}
\begin{tabular}{l>{\raggedright\arraybackslash$}p{12cm}<{$}}
$A[d;27](t)= $& t^{26}+3t^{25}+23t^{24}+152t^{23}+844t^{22}+3662t^{21}+12814t^{20}+36548t^{19}+86837t^{18}+173593t^{17}+295172t^{16}+
429240t^{15}+536646t^{14}+577821t^{13}+536646t^{12}+429240t^{11}+295172t^{10}+173593t^9+86837t^8+36548t^7+12814t^6+3662t^5+
844t^4+152t^3+23t^2+3t+1\\
\end{tabular}


\bigskip
\hspace*{-10ex}
\begin{tabular}{l>{\raggedright\arraybackslash$}p{12cm}<{$}} \small
$A[u;125](t) =$ & 
t^{124}+3t^{122}+45t^{120}+774t^{118}+11207t^{116}+129485t^{114}+1229657t^{112}+9835988t^{110}+67622641t^{108}+405731843t^{106}+2150382085t^{104}+10165426468t^{102}+43203077195t^{100}+166165624857t^{98}+581579739591t^{96}+1861054998416t^{94}+5466849215583t^{92}+14792650391699t^{90}+36981626382405t^{88}+85641660162366t^{86}+184129570236171t^{84}+368259138698205t^{82}+686301123812811t^{80}+1193567168903172t^{78}+1939546652290065t^{76}+2948110907190899t^{74}+4195388602819760t^{72}+5593851464926268t^{70}+6992314336461413t^{68}+8197885767564289t^{66}+9017674350331611t^{64}+9308567065105337t^{62}+9017674350331611t^{60}+8197885767564289t^{58}+6992314336461413t^{56}+5593851464926268t^{54}+4195388602819760t^{52}+2948110907190899t^{50}+1939546652290065t^{48}+1193567168903172t^{46}+686301123812811t^{44}+368259138698205t^{42}+184129570236171t^{40}+85641660162366t^{38}+36981626382405t^{36}+14792650391699t^{34}+5466849215583t^{32}+1861054998416t^{30}+581579739591t^{28}+166165624857t^{26}+43203077195t^{24}+10165426468t^{22}+2150382085t^{20}+405731843t^{18}+67622641t^{16}+9835988t^{14}+1229657t^{12}+129485t^{10}+11207t^8+774t^6+45t^4+3t^2+1\\
\end{tabular}
%
}
\hspace*{-15ex}
{\scriptsize
\begin{tabular}{l>{\raggedright\arraybackslash$}p{14cm}<{$}}
$A[d;125](t)=$  & 
t^{124}+3t^{123}+90t^{122}+3183t^{121}+94261t^{120}+2253202t^{119}+44660526t^{118}+752765426t^{117}+11009026889t^{116}+141893725177t^{115}+1631777381270t^{114}+16911146021617t^{113}+159246624819695t^{112}+1371970915393992t^{111}+10877769404828584t^{110}+79770308932154652t^{109}+543435229633787791t^{108}+3452412046870263651t^{107}+20522671612153800248t^{106}+114494904782463078927t^{105}+601098250109006348605t^{104}+2976867524344040005968t^{103}+13937152500343088195264t^{102}+61808241523238111266288t^{101}+260109683076981986458211t^{100}+1040438732307841539649001t^{99}+3961670557633787406854497t^{98}+14379396838818630509486185t^{97}+49814339048764832197753486t^{96}+164902639609703309488079986t^{95}+522191692097394743906419238t^{94}+1583419969585645756697355513t^{93}+4601814286608285713826377107t^{92}+12829300435392788915059242212t^{91}+34337245282963060080748456842t^{90}+88295773584762135474231228583t^{89}+218286773584550853413218320189t^{88}+519168542579472256031665746576t^{87}+1188622715905633892156485804100t^{86}+2621065476099602847252984385458t^{85}+5569764136711656142477377352019t^{84}+11411224084970222152064822132446t^{83}+22550752358393534437019574995528t^{82}+43003760311355111831092568853122t^{81}+79166013300449183486752791791953t^{80}+140739579200798547810848407179486t^{79}+241704929497023593576284109748137t^{78}+401127329803571069200022148395060t^{77}+643475091559895257811475523869150t^{76}+998042999154123255509578916762251t^{75}+1497064498731184884738423828738826t^{74}+2172211233453091791403201171682487t^{73}+3049450385424532709259880563769891t^{72}+4142649580199365187088253032770042t^{71}+5446817040632498674709258423750752t^{70}+6932312597168634673342968879821758t^{69}+8541599450082782011722315254853863t^{68}+10189978291326827658936469293236257t^{67}+11771181819291335403215209136236758t^{66}+13167762713105561632909721570633733t^{65}+14265076272531025106827702213046251t^{64}+14966637400688288631941104298816390t^{63}+15208034778118744904852502416921288t^{62}+ 
%
%
+14966637400688288631941104298816390t^{61}+14265076272531025106827702213046251t^{60}+13167762713105561632909721570633733t^{59}+11771181819291335403215209136236758t^{58}+10189978291326827658936469293236257t^{57}+8541599450082782011722315254853863t^{56}+6932312597168634673342968879821758t^{55}+5446817040632498674709258423750752t^{54}+4142649580199365187088253032770042t^{53}+3049450385424532709259880563769891t^{52}+2172211233453091791403201171682487t^{51}+1497064498731184884738423828738826t^{50}+998042999154123255509578916762251t^{49}+643475091559895257811475523869150t^{48}+401127329803571069200022148395060t^{47}+241704929497023593576284109748137t^{46}+140739579200798547810848407179486t^{45}+79166013300449183486752791791953t^{44}+43003760311355111831092568853122t^{43}+22550752358393534437019574995528t^{42}+11411224084970222152064822132446t^{41}+5569764136711656142477377352019t^{40}+2621065476099602847252984385458t^{39}+1188622715905633892156485804100t^{38}+519168542579472256031665746576t^{37}+218286773584550853413218320189t^{36}+88295773584762135474231228583t^{35}+34337245282963060080748456842t^{34}+12829300435392788915059242212t^{33}+4601814286608285713826377107t^{32}+1583419969585645756697355513t^{31}+522191692097394743906419238t^{30}+164902639609703309488079986t^{29}+49814339048764832197753486t^{28}+14379396838818630509486185t^{27}+3961670557633787406854497t^{26}+1040438732307841539649001t^{25}+260109683076981986458211t^{24}+61808241523238111266288t^{23}+13937152500343088195264t^{22}+2976867524344040005968t^{21}+601098250109006348605t^{20}+114494904782463078927t^{19}+20522671612153800248t^{18}+3452412046870263651t^{17}+543435229633787791t^{16}+79770308932154652t^{15}+10877769404828584t^{14}+1371970915393992t^{13}+159246624819695t^{12}+16911146021617t^{11}+1631777381270t^{10}+141893725177t^9+11009026889t^8+752765426t^7+44660526t^6+2253202t^5+94261t^4+3183t^3+90t^2+3t+1 
\end{tabular}
%
}

\bigskip
Tables \ref{tbl:genfunctions27} to \ref{tbl:alldir27} in Appendix B give the generating functions for some of the intermediate terms given in Table \ref{tbl:intermediatevalues}. It will not be very illuminating for the reader should we give the generating functions for all the terms, so we present in  Appendix B only those for which interesting and sometimes surprising relationships occur and which we shall discuss in the next section. However, in order to give a complete result at least for one case, we give in Table \ref{tbl:alldir27} in Appendix B, the generating functions of all the intermediate terms appearing in the directed case of order 27.

We also observe that the number of nonisomorphic undirected circulants on 27 vertices is now verified by Matan Zif-Av's brute-force methods as well as that of McKay's from 1995, and by the structural and multiplier approaches described here. The generating functions for these circulant are now also verified by the structural (for the undirected case of order $27$) and multiplier approaches and also by Zif-Av's methods. The value $A[d;27]=3728891$ has recently appeared in \cite{OEIS} (Sequence A04929);
we also note that $A_{21}[d;27]=3730584$ is represented there as well, in A056391.

\section{Discussion of results: unexpected patterns} 



\subsection{Some observations and identities, and a conjecture \label{i}} 

Our tables
have been derived in Chapter 4 of \cite{arXivVg}. 
A direct phenomenological analysis of the main and intermediate analytical formulae shown in these tables
reveals some hidden patterns that need to be explained
in general, either combinatorially, algebraically or analytically.
First of all, in the four cases we are studying (undirected / directed, $p=3$ / $p=5$) the following three `coincidences' are observed.
\[
A_{31}[s;p^3]=A_{41}[s;p^3]\eqno(\ref{i}.1)
\]
\[
A_{32}[s;p^3]=A_{42}[s;p^3]\eqno(\ref{i}.2)
\]
\[
A_{521}[s;p^3]=A_{522}[s;p^3] \eqno(\ref{i}.3)
\]
and (as a corollary of the first two)
\[
A_3[s;p^3]=A_4[s;p^3]. \eqno(\ref{i}.4)
\]
for $p=3,5$ and $s=u,d$
(where $A_4[s;p^3]:=A_{41}[s;p^3]-A_{42}[s;p^3]$).

For example $A_{31}[u;125]=A_{41}[u;125]=1272$. Notice that their enumeration formulae are distinct.
Moreover, refined by valencies, the corresponding pair of generating functions are also distinct.
However, unexpectedly at first sight, the multisets of coefficients in these pairs of polynomials coincide.
A more thorough analysis
enabled us to reveal a simple pattern.
Namely, in all four cases, we observe the following identities:
\[
A_{31}[s;p^3](t)\equiv A_{41}[s;p^3](t^p)\,\,({\rm mod}\,\,t^{p^3-1}), \eqno(\ref{i}.1^t)
\]
\[
A_{32}[s;p^3](t)\equiv A_{42}[s;p^3](t^p)\,\,({\rm mod}\,\,t^{p^3-1}), \eqno(\ref{i}.2^t)
\]
\[
A_{522}[s;p^3](t)\equiv A_{521}[s;p^3](t^p)\,\,({\rm mod}\,\,t^{p^3-1}), \eqno(\ref{i}.3^t)
\]
and most spectacularly, as a corollary of the first two,
\[
A_3[s;p^3](t)\equiv A_4[s;p^3](t^p)\,\,({\rm mod}\,\,t^{p^3-1}). \eqno(\ref{i}.4^t)
\]
In particular, for the latter identity and undirected graphs we observe that their expressions are

\begin{tabular}{lc>{\raggedright\arraybackslash$}p{8cm}<{$}}
$A_4[u;27](t)$ & = & t^{24}+2t^{22}+t^{20}+t^6+2t^4+t^2,
\end{tabular}

\begin{tabular}{lcccr}
$A_3[u;27](t)$ &=& $t^{20}+t^{18}+2t^{14}+2t^{12}+t^8+t^6$ & $\equiv$ & $A_4[u;27](t^3)\,\,({\rm mod}\,\,t^{26})$
\end{tabular}

\begin{tabular}{lc>{\raggedright\arraybackslash$}p{8cm}<{$}}
$A_4[u;125](t)$ &=&
t^{122}+7t^{120}+22t^{118}+51t^{116}+79t^{114}+94t^{112}+79t^{110}+51t^{108}+22t^{106}+7t^{104}+t^{102}+t^{72}+7t^{70}+22t^{68}+51t^{66}+79t^{64}+94t^{62}+79t^{60}+51t^{58}+22t^{56}+7t^{54}+t^{52}+t^{22}+7t^{20}+22t^{18}+51t^{16}+79t^{14}+94t^{12}+79t^{10}+51t^8+22t^6+7t^4+t^2,
\end{tabular}

\begin{tabular}{lc>{\raggedright\arraybackslash$}p{8cm}<{$}}
$A_3[u;125](t)$ &=&
t^{114}+t^{112}+t^{110}+7t^{104}+7t^{102}+7t^{100}+22t^{94}+22t^{92}+22t^{90}+51t^{84}+51t^{82}+51t^{80}+79t^{74}+79t^{72}+79t^{70}+94t^{64}+94t^{62}+94t^{60}+79t^{54}+79t^{52}+79t^{50}+51t^{44}+51t^{42}+51t^{40}+22t^{34}+22t^{32}+22t^{30}+7t^{24}+7t^{22}+7t^{20}+t^{14}+t^{12}+t^{10},
\end{tabular}\newline
and therefore $A_3[u;125](t)$ is congruent to 
$A_4[u;125](t^5)\,\,({\rm mod}\,\,t^{124})$.

Thus, for example, $A_4[u;125]$ contributes $22$ circulant graphs of valency $6$ into the overall sum, and the same
number of circulant graphs of valency $30=5\times 6$ is contributed by $A_3[u;125]$.


Notice that the transformation
\[
\eta_{p,3}:\, t\to t^p\quad {\rm modulo}\quad t^{p^3-1} 
\]
in the ring of polynomials of $t$ over the rationals is periodic of order $3$. Thus,
${A_4[s;p^3](t)\equiv A_3[s;p^3](t^{p^2})\,\,({\rm mod}\,\,t^{p^3-1})}$, etc.
Besides, this operation fixes the terms $d\cdot t^{e(p^2+p+1)},\ e=0,1,\dots,p-1$.

Finally, more hidden identities of the same nature are valid in all four cases:
$A_{1}[s;p^3](t)$ and $A_{523}[s;p^3](t)$ are invariant with respect to $\eta_{p,3}$, that is,
\[
A_1[s;p^3](t)\equiv A_1[s;p^3](t^p)\,\,({\rm mod}\,\,t^{p^3-1}), \eqno(\ref{i}.5^t)
\]
and
\[
A_{523}[s;p^3](t)\equiv A_{523}[s;p^3](t^p)\,\,({\rm mod}\,\,t^{p^3-1}), \eqno(\ref{i}.6^t)
\]
as can be seen from our tables. 
Of course~(\ref{i}.5$^t$) and~(\ref{i}.6$^t$) make no sense for valency-unspecified circulants ($t=1$).


We conjecture that the above identities hold in general.

{\bf Conjecture.}
\textit{Identities~{\rm (\ref{i}.1$^t$) -- (\ref{i}.3$^t$)}, {\rm (\ref{i}.5$^t$)}
and~{\rm (\ref{i}.6$^t$)}
(and, consequently, identities~{\rm (\ref{i}.4$^t$)} and~{\rm (\ref{i}.1) -- (\ref{i}.4)}) are valid
in general for all odd prime $p$ and $s=u,d$.}

\medskip
If valid in general, this conjecture should have a simple formal analy\-ti\-cal proof.
For comparison, this is the case for two identities similar to~(\ref{i}.5$^t)$ and {\rm (\ref{i}.6$^t$)}
that are valid for intermediate enumerative polynomials for circulant graphs of prime-squared orders (cf. \cite{L03,KLP96});
see Appendix A. They promise a simple analytical proof and even suggest the existence of a direct
bijective proof of combinatorial or algebraic nature (cf. \cite{VL}).
One idea to guess such valency-violating bijections in our particular cases is the following: to extract from the relative
pairs of polynomials the corresponding terms with the coefficients equal to 1 and to compare their
corresponding single graphs. For example, as we see above, $A_4[u;125](t)$ contributes a unique circulant graph
of valency $2$, and it corresponds to a certain unique circulant graph of valency $5\times 2=10$
counted in $A_3[u;125](t)$. Likewise, the unique contributors of valencies $22$ and $52$
correspond to those of valencies $110$ and $12\equiv 5\times 52\ ({\rm mod}\ 124)$, respectively,
in $A_3[u;125](t)$ (furthermore, by complementarity, $110$ may be replaced with $14$). Hopefully the
structural approach can help here (maybe even within the rather elementarily framework of the Isomorphism
Theorem and related results?); in such a case the identities would gain some value for the structural theory of circulant graphs.
Perhaps a link can be established between these formulae and the figures discussed in Section \ref{sec:misha} whose valencies are multiples of $p^i, i=1,2$, while the graphs have clear homomorphic images of smaller size, though still belonging to the same variety of prime-power circulants.


Various formal identities are rather characteristic for the enumerators
of circulant graphs of prime or prime-squared orders (\cite{L03,KLP96}); 
but the present ones, if valid in general, are of a new nature:
valency-violating although by a simple rule.

It is interesting to note that, rather unexpectedly, they have served as a hint for the discovery of similar valency-violating identities for intermediate classes of circulants of order $p^2,\ p\geq3$. Details concerning these new identities and their analytical proof can be found in Appendix A. 

\subsection{Enumeration of self-com\-ple\-men\-ta\-ry circulants of orders $27$ and $125$}\label{S-C}

The generating functions $A[s;p^3](t)$ make it possible to calculate easily the numbers of the corresponding
self-complementary circulants (a graph is called
\textit{self-complementary} if it is isomorphic to its complementary graph).

\begin{Proposition}[\cite{VL17}]
	\label{Disc1}

{For any odd prime $p$ and integer $k\geq 1$,
$$
A[u_{sc};p^k] = A[u;p^k](t)|_{t^2:=-1},  \eqno(\ref{S-C}.1)
$$
$$
A[d_{sc};p^k] = A[d;p^k](t)|_{t:=-1},    \eqno(\ref{S-C}.2)
$$
where $u_{sc}$ and $d_{sc}$ stand for undirected and directed self-com\-ple\-mentary circulant graphs, resp.}
\end{Proposition}

\noindent
Proof.
In the framework of Redfield--P\'olya enumeration theory there is a well-known general approach
(going back to de Bruijn and even to Redfield himself) to counting self-dual configurations including
self-complementary graphs; see, e.g., sec.6.2 in~\cite{HP73} and~\cite{PR84,Fa99}. In short, manipulating with the
cycle index one needs to preserve variables corresponding to cycles of even order and to exclude (vanish) ones that
correspond to odd cycles, namely, to substitute ${x_{2r}:=2}$ and ${x_{2r-1}:=0}$, $r=1,2,\dots$ In terms of the standard
substitution $x_r:=1+t^r$ this is obtained by the subsequent substitution ${t:=-1}$. This approach cannot be applied
directly in our case since $A[u;p^k]$ and $A[d;p^k]$ are not solutions of single problems of this type. But we can
apply it to \textsl{each } subproblem out of the problems of the Redfield--P\'olya type to which our enumeration
is reduced according to~\cite{LP2000}. The last thing to be clarified is the usage of variables. In the current research
we specify circulant graphs of order $n=p^k$ by valency $r$ instead of the usual specification by the number of edges $N$.
But both parameters are related tightly: $N=nr$ for digraphs and $N=nr/2$ for undirected graphs ($r$ is even
in the latter case). Therefore $A[d;n](t^n)$ is the generating function for circulant digraphs by the number of edges
and $A[u;n](t^{n/2})$ is the same for undirected graphs. Thus, counting self-com\-ple\-men\-tary circulants we are to
substitute ${t:=-1}$ in these transformed polynomials. But since $n$ is odd this gives rise to the desired formulae
(\ref{S-C}.1) and (\ref{S-C}.2). \qed
\medskip

\textbf{Remark.} It was conjectured in \cite{L03}, Conjecture 6.1, that the same assertation is valid for arbitrary odd orders $n$. The right-hand-side expressions in (\ref{S-C}.1) and (\ref{S-C}.2) are the alternating sums
of the coefficients, and we may reformulate Proposition~\ref{Disc1} in terms of the general pattern (cf. A000171 in \cite{OEIS} and also \cite{NST09})
$$
sc(n) = e(n) - o(n),       \eqno(\ref{S-C}.3)
$$
where $e(n)$ and $o(n)$ stand for the number of non-isomorphic graphs of a certain class with even and odd
number of edges, resp., and $sc(n)$ stands for the number of self-com\-ple\-men\-tary graphs (sc-graphs for short)
 of the same class.


\begin{Corollary}\label{Disc2}
The values for $A[s_{sc};p^3]$ for $s=u,d$ and $p=3,5$ are given by
\begin{center}
\begin{tabular}{lcl}
 $A[u_{sc};27]$  &=&                 $0$,\\  
 $A[u_{sc};125]$ &=&          $42949681$,\\  
 $A[d_{sc};27]$  &=&               $457$,\\  
 $A[d_{sc};125]$ &=& $46116860227224068$.    
\end{tabular}
\end{center}

\end{Corollary}

The vanishing of $A[u_{sc};27]$ is obvious since any undirected sc-graph of order $n$ contains the median number
$n(n-1)/4$ of edges but $27\times26/4$ is fractional.

It makes sense to calculate the corresponding values for the intermediate generating functions (using
the same substitutions as in Proposition~\ref{Disc1}):
these are the numbers of self-complementary circulant graphs of the corresponding subclasses.
Table~\ref{SC-t1} contains the refined numeric data for self-com\-ple\-men\-tary circulant digraphs of order $27$.
All these digraphs are of valency $r=13$; for comparison in the last column we included the numbers of all
circulant digraphs of this valency. The intermediate contributors for all four classes of circulants are
represented in Table~\ref{SC-t2}.

\begin{table}[ht]
\caption{Intermediate contributors for counting self-com\-ple\-men\-tary circulant digraphs of order $n=27$}\label{SC-t1}
\centering
\begin{tabular}{|r|l|r||r|}
\hline
Contributor&Interconnection &\#(sc-circ)&\#(circ of val $13$)\\
\hline
\hline
$A_1[d; 27](-1)$&=&                                 1&       1\\
\hline
$A_{21}[d; 27](-1)$&=&                             472& 577996\\
$A_{22}[d; 27](-1)$&=&                              24&    276\\
\hline
$A_2[d;27](-1)$&$=A_{21}-A_{22}=$&                 448& 577720\\
\hline
$A_{31}[d; 27](-1)$&=&                               4&     14\\
$A_{32}[d; 27](-1)$&=&                               2&      2\\
\hline
$A_3[d;27](-1)$&$=A_{31}-A_{32}=$&                   2&     12\\
\hline
$A_{41}[d; 27](-1)$&=&                               4&     14\\
$A_{42}[d; 27](-1)$&=&                               2&      2\\
\hline
$A_4[d;27](-1)$&$=A_{41}-A_{42}=$&                   2&     12\\
\hline
$A_{51}[d; 27](-1)$&=&                              16&    124\\
$A_{521}[d; 27](-1)$&=&                              8&     26\\
$A_{522}[d; 27](-1)$&=&                              8&     26\\
$A_{523}[d; 27](-1)$&=&                              4&      4\\
\hline
$A_5[d;27](-1)$&=$A_{51}\!-\!A_{521}\!-\!A_{522}\!+\!A_{523}$=&4&76\\
\hline
$A[d; 27](-1)$&$=A_1\!+\!A_2\!+\!A_3\!+\!A_4\!+\!A_5=$&457& 577821\\
\hline
\end{tabular}
\end{table}

\tabcolsep=0.65ex
\begin{table}[htb!]
\caption{\label{tbl:s-c-intermediatevalues} The numbers of self-complementary undirected and directed
circulant graphs of orders $27$ and $125$; subclasses}\label{SC-t2}
\begin{tabular}{|l|r|r|r|r|}
\hline
Term (s-c)&            Undir&  Undir&Dir&     Dir          \\
        &             $n\!=\!27$&$n=125$&$n\!=\!27$& $n=125$       \\
\hline
\hline
$A_1$&                         0&       1&  1&                8\\
\hline
$A_{21}$&                      0&       42949840 &472&46116860227391504\\
$A_{22}$&                      0&       208 & 24&           209936\\
\hline
$A_2=A_{21}-A_{22}$&           0&42949632&448&46116860227181568\\
\hline
$A_{31}=A_{41}$&               0&       8&  4&              432\\
$A_{32}=A_{42}$&               0&       2&  2&               12\\
\hline
$A_3=A_{4}=A_{31}-A_{32}$&     0&       6&  2&              420\\
\hline
$A_{51}$&                      0&      64& 16&            43328\\
$A_{521}=A_{522}$&             0&      16&  8&              848\\
$A_{523}$&                     0&       4&  4&               20\\
\hline
$A_5=A_{51}\!-\!A_{521}\!-\!A_{522}\!+\!A_{523}$&0& 36&  4&41652\\
\hline
$A=A_1+A_2+A_3+A_4+A_5$&       0&42949681&457&46116860227224068\\
\hline
\end{tabular}
\medskip
\end{table}

\clearpage
\noindent
\textbf{Concluding comments.}

(1) All entries in the first row of Tab.~\ref{SC-t2} coincide with the median coefficients of the corresponding
generating functions $A_1[,](t)$; that is, all circulant graphs of the corresponding valency (which, moreover, are
unique in two cases, as we see) are self-com\-ple\-men\-tary. Is this a general pattern?

(2) For prime squared-orders there are several identities in which enu\-mera\-tors of sc-circulants are involved~\cite{L03}.
We may expect something similar for prime-cubed orders.

\section*{Acknowledgements}

The authors would like to thank Matan Zif-Av for several helpful suggestions in the writing of GAP programs and for verifying some of our numerical results using diverse methods.

This article is partially based upon work supported by the Project: Mobility - enhancing research, science and education at Matej Bel University, ITMS code 26110230082, under the Operational Program Education co-financed by the European Social Fund. The authors also gratefully acknowledge support from the  Scientific Grant
Agency of the Slovak Republic
under the number VEGA-1/0988/16.





\clearpage

\section*{ }
\addcontentsline{toc}{section}{References}

\bibliography{KlinLiskoArXivV0_0}

\begin{thebibliography}{10}

\bibitem{adam67}
{\'A}.~{\'A}d{\'a}m.
\newblock Research problem 2-10.
\newblock {\em J. Combin. Theory}, 2:393, 1967.

\bibitem{dixM96}
J.D. Dixon and B.~Mortimer.
\newblock {\em Permutation groups}, volume 163 of {\em Graduate Texts in
  Mathematics}.
\newblock Springer-Verlag, New York, 1996.

\bibitem{dob95}
E.~Dobson.
\newblock An extension of a classical result of {B}urnside.
\newblock Manuscript, 1995.

\bibitem{dob00}
E.~Dobson.
\newblock Classification of vertex-transitive graphs of order a prime cubed.
  {I}.
\newblock {\em Discrete Math.}, 224(1-3):99--106, 2000.

\bibitem{dob05}
E.~Dobson.
\newblock On groups of odd prime-power degree that contain a full cycle.
\newblock {\em Discrete Math.}, 299(1-3):65--78, 2005.

\bibitem{dob10}
E.~Dobson.
\newblock Asymptotic automorphism groups of {C}ayley digraphs and graphs of
  abelian groups of prime-power order.
\newblock {\em Ars Math. Contemp.}, 3(2):200--213, 2010.

\bibitem{dobM09}
E.~Dobson and J.~Morris.
\newblock Automorphism groups of wreath product digraphs.
\newblock {\em Electron. J. Combin.}, 16(1):Research Paper 17, 30, 2009.

\bibitem{dobW02}
E.~Dobson and D.~Witte.
\newblock Transitive permutation groups of prime-squared degree.
\newblock {\em J. Algebraic Combin.}, 16(1):43--69, 2002.

\bibitem{elspas&tur70}
B.~Elspas and J.~Turner.
\newblock Graphs with circulant adjacency matrices.
\newblock {\em J. Combin. Theory}, 9:297--307, 1970.

\bibitem{faradzev&l90}
I.A. Farad\v{z}ev, A.A. Ivanov, and M.~Klin.
\newblock Galois correspondence between permutation groups and cellular rings
  (association schemes).
\newblock {\em Graphs and Combin.}, (number):303--332, 1990.

\bibitem{faradzev&klin91}
I.A. Farad{\v{z}ev} and M.~Klin.
\newblock Computer package for computation with coherent configurations.
\newblock pages 219--223. Proc. ISSAC-91, ACM Press, 1991.

\bibitem{Fa99}
A.~Farrugia.
\newblock Self-complementary graphs and generalizations: a comprehensive
  reference manual.
\newblock Master's thesis, University of Malta, 1999.

\bibitem{arXivVg}
V.~Gatt.
\newblock On the enumeration of circulant graphs of prime-power order: the case
  of $p^3$. {W}ith {A}ppendix {C}: {C}omplete list of generating functions.
\newblock Master's thesis, University of Malta, 2013.
\newblock \url{https://arxiv.org/abs/1703.06038v2}.

\bibitem{golNP85}
J.J. Gol'fand, N.L. Najmark, and R.~P\"{o}schel.
\newblock The structure of {$S$}-rings over $\mathbb{Z}_{2^m}$.
\newblock Preprint {I}nst. f. {M}ath, AdWDDR, Berlin
  \url{http://math.tu-dresden.de/~poeschel/poePUBLICATIONSpdf/1985GolfandNajPoe.pdf},
  1985.

\bibitem{HP73}
F.~Harary and E.M. Palmer.
\newblock {\em Graphical Enumeration}.
\newblock Academic Press, 1973.
\newblock [pp. 139, 140, 243].

\bibitem{kli94}
M.~Klin.
\newblock Automorphism groups of $p^m$ circulant graphs.
\newblock In {\em The First China-Japan International Symposium on Algebraic
  Combinatorics}, pages 56--58. October 11--15, 1994.
\newblock Beijing.

\bibitem{kliK12}
M.~Klin and I.~Kov\'acs.
\newblock Automorphism groups of rational circulant graphs.
\newblock {\em Electron. J. Combin.}, 19(1):Paper 35, 52, 2012.

\bibitem{KLP96}
M.~Klin, V.~Liskovets, and R.~P\"oschel.
\newblock Analytic enumeration of circulant graphs with prime-squared number of
  vertices.
\newblock {\em S\'em. Lotharing. Combin.}, B36d:36, 1996.

\bibitem{KLP03}
M.~Klin, V.~Liskovets, and R.~P\"oschel.
\newblock On the analytical enumeration of circulant graphs.
\newblock Technical Report MATH-AL-07-2003, TU-Dresden, 2003.

\bibitem{kliNajP81}
M.~Klin, N.L. Najmark, and R.~P\"{o}schel.
\newblock Schur rings over $\mathbb{Z}_{2^m}$.
\newblock Preprint {P}-{MATH}-14/81 {A}kad. der {W}iss. der {DDR}, {ZIMM},
  {B}erlin.
  \url{http://math.tu-dresden.de/~poeschel/poePUBLICATIONSpdf/1981KlinNajmPoe.pdf},
  1981.

\bibitem{klinP75}
M.~Klin and R.~P\"{o}schel.
\newblock The isomorphism problem for cyclic graphs with $p^2$ or $pq$
  vertices.
\newblock Abstract presented at the {A}.{A}. {Z}ykov {S}eminar on {G}raph
  {T}heory, Odessa ({R}ussian), 1974.

\bibitem{kliP78}
M.~Klin and R.~P\"{o}schel.
\newblock The {K}\"{o}nig problem, the isomorphism problem for cyclic graphs
  and the characterization of {S}chur rings.
\newblock Preprint {A}d{W}d {DDR}, {ZIMM}, {B}erlin.
  \url{http://math.tu-dresden.de/~poeschel/poePUBLICATIONSpdf/1978\_Klin\_PoePP.pdf},
  1978.

\bibitem{KP80}
M.~Klin and R.~P\"{o}schel.
\newblock The isomorphism problem for circulant digraphs with $p^n$ vertices.
\newblock {\em Preprint P34/80. Berlin , AdWDDR, ZIMM}, 1980.

\bibitem{kliP80}
M.~Klin and R.~P\"{o}schel.
\newblock The isomorphism problem for circulant digraphs with $p^n$ vertices.
\newblock Preprint {P}-{MATH}-34/80 {A}kad. der {W}iss. der {DDR}, {ZIMM},
  {B}erlin, 1980.

\bibitem{klinP81}
M.~Klin and R.~P\"oschel.
\newblock The {K}\"onig problem, the isomorphism problem for cyclic graphs and
  the method of {S}chur rings.
\newblock In {\em Algebraic methods in graph theory, {V}ol. {I}, {II}
  ({S}zeged, 1978)}, volume~25 of {\em Colloq. Math. Soc. J\'anos Bolyai},
  pages 405--434. North-Holland, Amsterdam-New York, 1981.

\bibitem{klin&al91}
M.~Klin, C.~R\"{u}cker, G.~R\"{u}cker, and G.~Tinhofer.
\newblock Algebraic combinatorics in mathematical chemistry. methods and
  algorithms. {I}. permutation groups and coherent (cellular) algebras.
\newblock {\em Match}, 40:7--138, 1999.

\bibitem{kov02}
I.~Kov\'acs.
\newblock On automorphisms of circulant digraphs on {$p^m$} vertices, {$p$} an
  odd prime.
\newblock {\em Linear Algebra Appl.}, 356:231--252, 2002.
\newblock Special issue on algebraic graph theory (Edinburgh, 2001).

\bibitem{lauri&sca16}
J.~Lauri and R.~Scapellato.
\newblock {\em Topics in Graph Automorphisms and Reconstruction}.
\newblock Number 432 in London Mathematical Society Lecture Notes Series, 432.
  Cambridge University Press, 2016.

\bibitem{L03}
V.~Liskovets.
\newblock Some identities for enumerators of circulant graphs.
\newblock {\em J. Algebraic Combin.}, 18(3):189--209, 2003.

\bibitem{VL}
V.~Liskovets.
\newblock New identities in the enumeration of circulant graphs of
  prime-squared orders.
\newblock {\em Proc. XII Belarusian Mathem. Conf. (September 2016, Minsk).
  Inst. of Math. Nat. Acad. Sci. of Belarus}, P.4:73--74, 2016.

\bibitem{VL17}
V.~Liskovets.
\newblock Enumeration of self-complementary circulant graphs of prime-power
  order: old and new results.
\newblock {\em The {S}econd {M}alta {C}onference in {G}raph {T}heory and
  {C}ombinatorics (2{M}{C}{G}{T}{C}. June 2017, Qawra). {U}niversity of
  {M}alta, Abstracts}, pages 85--86, 2017.

\bibitem{LP2000}
V.~Liskovets and R.~P\"{o}schel.
\newblock Counting circulant graphs of prime-power order by decomposing into
  orbit enumeration problems.
\newblock {\em Discrete Math.}, 214:173--191, 2000.

\bibitem{mis15}
A.~Misseldine.
\newblock Counting {S}chur rings over cyclic groups.
\newblock Preprint. \url{https://arxiv.org/abs/1508.03757}, 2015.

\bibitem{muzychuk95}
M.~Muzychuk.
\newblock {\'A}d{\'a}m's conjecture is true in the square-free case.
\newblock {\em J. Combin. Theory (Ser. A)}, 72:118--134, 1995.

\bibitem{muz04}
M.~Muzychuk.
\newblock A solution of the isomorphism problem for circulant graphs.
\newblock {\em Proc. London Math. Soc. (3)}, 88(1):1--41, 2004.

\bibitem{muzKP01}
M.~Muzychuk, M.~Klin, and R.~P\"oschel.
\newblock The isomorphism problem for circulant graphs via {S}chur ring theory.
\newblock In {\em Codes and association schemes ({P}iscataway, {NJ}, 1999)},
  volume~56 of {\em DIMACS Ser. Discrete Math. Theoret. Comput. Sci.}, pages
  241--264. Amer. Math. Soc., Providence, RI, 2001.

\bibitem{NST09}
A.~Nakamoto, T.~Shirakura, and Sh. Tazawa.
\newblock An alternative enumeration of self-complementary graphs.
\newblock {\em Utilitas Mathematica}, 80:25--32, 2009.

\bibitem{PR84}
E.M. Palmer and R.W. Robinson.
\newblock Enumeration of self-dual configurations.
\newblock {\em Pacific J. of Math.}, 110(1):203--221, 1984.

\bibitem{pos74}
R.~P\"{o}schel.
\newblock Untersuchungen von {S}-{R}ingen, insbensondere im {G}ruppenring von
  $p$-{G}ruppen.
\newblock {\em Math. Nach.}, 60:1--27, 1974.

\bibitem{posK79}
R.~P\"oschel and L.~A. Kalu\v{z}nin.
\newblock {\em Funktionen- und {R}elationenalgebren}, volume~15 of {\em
  Mathematische Monographien [Mathematical Monographs]}.
\newblock VEB Deutscher Verlag der Wissenschaften, Berlin, 1979.
\newblock Ein Kapitel der diskreten Mathematik. [A chapter in discrete
  mathematics].

\bibitem{sch33}
I.~Schur.
\newblock Zur {T}heorie der einfach transitiven {P}ermutationgruppen.
\newblock {\em S.-B. Preuss. Akad. Wiss. phys.-math. Kl.}, 18/20:598--623,
  1933.

\bibitem{sch03}
I.~Schur.
\newblock On the theory of transitive permutation groups.
\newblock {\em Zap. Nauchn. Sem. S.-Peterburg. Otdel. Mat. Inst. Steklov.
  (POMI)}, 301(Teor. Predst. Din. Sist. Komb. i Algoritm. Metody. 9):5--34,
  243, 2003.
\newblock Russian translation of \cite{sch33}.

\bibitem{sco87}
W.~R. Scott.
\newblock {\em Group theory}.
\newblock Dover Publications, Inc., New York, second edition, 1987.

\bibitem{OEIS}
N.J.A. Sloane~({E}ditor).
\newblock {T}he {O}n-{L}ine {E}ncyclopedia of {I}nteger {S}equences.
\newblock Published electronically at \url{https://oeis.org}.

\bibitem{gap99}
{The GAP Group}.
\newblock {\em {GAP}---Groups, Algorithms and Programming, Version 4.1}.
\newblock Aachen, St. Andrews, \url{http://www-gap.dcs.st-and.ac.uk/~gap},
  1999.

\bibitem{vysKCh78}
V.~A. Vy\v{s}enski\u\i, M.~H. Klin, and N.~I. \v{C}eredni\v{c}enko.
\newblock Realization of an algorithm for construction of {$S$}-rings of cyclic
  groups of order {$p^{m}$}\ and its application for solution of the problem of
  cataloguing $p^m$-vertex cyclic graphs.
\newblock In {\em Computations in algebra and combinatorial analysis
  ({R}ussian)}, pages 73--86. Akad. Nauk Ukrain. SSR, Inst. Kibernet., Kiev,
  1978.

\bibitem{wie64}
H.~Wielandt.
\newblock {\em Finite permutation groups}.
\newblock Translated from the German by R. Bercov. Academic Press, New
  York-London, 1964.

\bibitem{wie69}
H.~Wielandt.
\newblock {\em Permutation groups through invariant relations and invariant
  functions}.
\newblock Lect. Notes, Det. Math. Ohio State University, Columbia (Ohio), 1969.

\end{thebibliography}


\newpage


\section*{Appendix A: New identities for order $p^2$ circulant graphs}
\label{sec:AppA}
\addcontentsline{toc}{section}{Appendix A:  New identities for order $p^2$ circulant graphs}


For a better understanding of the identities discussed in Section \ref{i} it will be useful and instructive to
consider the prime-squared circulants and new identities for them. They follow easily from the enumerative
formulae obtained in \cite{KLP96}. Moreover, a direct 1-to-$p$ correspondence between appropriate intermediate
subsets of $p^2$-circulants used in the proofs in \cite{KLP96} (see the first paragraph of Section 7.3, p.27) suggests
the existence of a transparent bijective proof of these identities (serving as a sample, in the future, for $p^3$-circulants).
For brevity, we concentrate on directed circulant graphs.

By \cite{KLP96} (in distinct designations), the enumerative generating functions for valency specified circulant
digraphs of order $p^2\ (p>2$, prime) satisfies the general equation
\[
 A[d;p^2](t)=A_1[d;p^2](t)+A_{21}[d;p^2](t)-A_{22}[d;p^2](t) \eqno(a1)
\]
where the RHS terms are the generating functions for appropriate intermediate types of
circulant graphs (similar to ones introduced in the $p^3$ case) and are calculated via the cycle index of the regular cyclic group
${\mathfrak I}_n({\bf x}):=\frac{1}{n}\sum\limits_{r|n}\phi(r)x_r^{n/r},$
where ${\bf x}:=\{x_1,x_2,\dots\}$. Namely, consider polynomials
$$
 {\mathfrak D}(p^2;{\bf x,y}):={\mathfrak I}_{p-1}({\bf x}){\mathfrak I}_{p-1}({\bf y}), \eqno(a2)
$$
$$
 \hspace*{-6ex}{\mathfrak B}(p^2;{\bf x,y}):={\mathfrak I}_{p-1}({\bf xy}), \eqno(a3)
$$
where ${\bf xy}:=\{x_1y_1,x_2y_2,\dots\}$. Then

$$ 
 A_1[d;p^2](t)={\mathfrak D}(p^2;{\bf x,y})|_{\{x_r:=1+t^r,\ y_r:=1+t^{pr}\}_{r=1,2,\ldots}} \eqno(a4)
$$
$$
 \hspace*{-1ex}A_{22}[d;p^2](t)={\mathfrak B}(p^2;{\bf x,y})|_{\{x_r:=1+t^r,\ y_r:=1+t^{pr}\}_{r=1,2,\ldots}} \eqno(a5)
$$ 
A certain similar formula holds for $A_{21}[d;p^2](t)$ as well but it does not seem to possess any interesting pattern.

Now, let ${\mathfrak F}({\bf x,y})$ be an arbitrary multivariate polynomial and suppose that it is symmetric
with respect to the interchange of variables ${\bf x}\leftrightarrow{\bf y}$ (that is,
$x_1\leftrightarrow y_1,x_2\leftrightarrow y_2,\dots$), i.e. satisfies the identity
$$
 {\mathfrak F}({\bf x,y})={\mathfrak F}({\bf y,x}).\eqno(a6)
$$

Consider the following substitution of variables
$$
 G(t):={\mathfrak F}({\bf x,y})|_{\{x_r:=g_r(t),\ y_r:=g_r(t^p)\}_{r=1,2,\ldots}} \eqno(a7)
$$
where $g_r(t), r=1,2,\dots$, are arbitrary polynomials and $p$ is an arbitrary positive integer.
Then the following polynomial congruence is valid:
$$
 G(t)\equiv G(t^p)\quad ({\rm mod}\,\,t^{p^2-1}).
$$
Indeed, combining the substitution (a7) with $t\to t^p$ we obtain
$$
 G(t^p)={\mathfrak F}({\bf x,y})|_{\{x_r:=g_r(t^p),\ y_r:=g_r(t^{p^2})\}_{r=1,2,\ldots}}.
$$
Since $t^{p^2}\equiv t$ modulo $t^{p^2-1}$, we have
$$
 G(t^p)\equiv \tilde{G}(t)\quad ({\rm mod}\,\,t^{p^2-1})
$$
where
$$
 \tilde{G}(t)={\mathfrak F}({\bf x,y})|_{\{x_r:=g_r(t^p),\ y_r:=g_r(t)\}_{r=1,2,\ldots}}.
$$
But $\tilde{G}(t)={G}(t)$ due to the symmetry property (a6). \qed

The polynomials ${\mathfrak D}$ and ${\mathfrak B}$ defined in (a2) and (a3) are symmetric.
Therefore we have the following identities.

\medskip
\begin{Proposition}
$$
 A_1[d;p^2](t)\equiv A_1[d;p^2](t^p)\quad ({\rm mod}\,\,t^{p^2-1}), \eqno(a8)
$$
$$
 A_{22}[d;p^2](t)\equiv A_{22}[d;p^2](t^p)\quad ({\rm mod}\,\,t^{p^2-1}), \eqno(a9)
$$
\end{Proposition}

For example, for $p=7$ we have the following polynomials (as calculations show) and
easily verifiable congruences:
$$
\begin{array}{ll}
A_1[d;49](t)\ =
& t^{48}+t^{47}+3t^{46}+4t^{45}+3t^{44}+t^{43}+t^{42}+t^{41}+t^{40}+\\
& 3t^{39}+4t^{38}+3t^{37}+t^{36}+t^{35}+3t^{34}+3t^{33}+9t^{32}+12t^{31}+\\
& 9t^{30}+3t^{29}+3t^{28}+4t^{27}+4t^{26}+12t^{25}+16t^{24}+12t^{23}+\\
& 4t^{22}+4t^{21}+3t^{20}+3t^{19}+9t^{18}+12t^{17}+9t^{16}+3t^{15}+3t^{14}+\\
& t^{13}+t^{12}+3t^{11}+4t^{10}+3t^{9}+t^{8}+t^{7}+t^{6}+t^{5}+3t^{4}+4t^{3}+\\
& 3t^{2}+t+1\\
\hspace{12.5ex} \equiv & A_1[d;49](t^7)\,\,({\rm mod}\,\,t^{48})
\end{array}
$$
$$
\begin{array}{ll}
A_{22}[d;49](t)=
& t^{48}+t^{47}+3t^{46}+4t^{45}+3t^{44}+t^{43}+t^{42}+t^{41}+6t^{40}+\\
& 15t^{39}+20t^{38}+15t^{37}+6t^{36}+t^{35}+3t^{34}+15t^{33}+39t^{32}+\\
& 50t^{31}+39t^{30}+15t^{29}+3t^{28}+4t^{27}+20t^{26}+50t^{25}+68t^{24}+\\
& 50t^{23}+20t^{22}+4t^{21}+3t^{20}+15t^{19}+39t^{18}+50t^{17}+39t^{16}+\\
& 15t^{15}+3t^{14}+t^{13}+6t^{12}+15t^{11}+20t^{10}+15t^{9}+6t^{8}+t^{7}+\\
& t^{6}+t^{5}+3t^{4}+4t^{3}+3t^{2}+t+1\\
\hspace{12.5ex} \equiv & A_{22}[d;49](t^7)\,\,({\rm mod}\,\,t^{48})
\end{array}
$$

For the enumerators $A[u;p^2](t)$ and $A[o;p^2](t)$ of undirected and oriented
circulant graphs, respectively, of order $p^2$ the expressions similar to (a1) hold.
Accordingly for them the identities similar to (a8) and (a9) hold.
They follow likewise from the next equations, proven, respectively, in \cite{KLP96}:
$$ 
A_1[u;p^2](t)\ ={\mathfrak D}^*(p^2;{\bf x,y})|_{\{x_r:=1+t^{2r},\ y_r:=1+t^{2pr}\}_{r=1,2,\ldots}}, 
$$
$$
A_{22}[u;p^2](t)={\mathfrak B}^*(p^2;{\bf x,y})|_{\{x_r:=1+t^{2r},\ y_r:=1+t^{2pr}\}_{r=1,2,\ldots}}, 
$$
where
$
{\mathfrak D}^*(p^2;{\bf x,y}):={\mathfrak I}_\frac{p-1}2({\bf x}){\mathfrak I}_\frac{p-1}2({\bf y}), 
$
and
${\mathfrak B}^*(p^2;{\bf x,y}):={\mathfrak I}_\frac{p-1}2({\bf xy})$
and in \cite{KLP03}:
$$
A_1[o;p^2](t)\ ={\mathfrak D}(p^2;{\bf x,y})|_{\{x_r:=1,\ y_r:=1\}_{r\ {\rm even}},\ \{{x^2_r:=1+2t^r,\ y^2_r:=1+2t^{pr}\}_{r\ {\rm odd}} }}, 
$$
$$
A_{22}[o;p^2](t)={\mathfrak B}(p^2;{\bf x,y})|_{\{x_r:=1,\ y_r:=1\}_{r\ {\rm even}},\ \{{x^2_r:=1+2t^r,\ y^2_r:=1+2t^{pr}\}_{r\ {\rm odd}} }}. 
$$ 

\newpage


\section*{Appendix B: Tables of generating functions}
\label{sec:AppB}
\addcontentsline{toc}{section}{Appendix B: Tables of generating functions}

\begin{table}[htb!] 
\caption{\label{tbl:genfunctions27}  Intermediate generating functions for $n=27$}
\begin{tabular}{|l|>{\raggedright\arraybackslash$}p{12cm}<{$}|}
\hline
Term & \mbox{Generating function}\\
\hline
\hline
$A_1[u;27](t)$ & 
t^{26}+t^{24}+t^{20}+t^{18}+t^{8}+t^{6}+t^{2}+1\\
\hline
$A_1[d;27](t)$& t^{26}+t^{25}+t^{24}+t^{23}+t^{22}+t^{21}+t^{20}+t^{19}+t^{18}+t^{17}+t^{16}+
t^{15}+t^{14}+t^{13}+t^{12}+t^{11}+t^{10}+t^9+t^8+t^7+t^6+t^5+t^4+t^3+t^2+t+1\\
\hline
$A_{31}[u;27](t)$ & 
t^{26}+t^{24}+2t^{20}+2t^{18}+2t^{14}+2t^{12}+2t^8+2t^6+t^2+1\\
\hline
$A_{31}[d;27](t)$ & 
t^{26}+t^{25}+t^{24}+2t^{23}+2t^{22}+2t^{21}+6t^{20}+6t^{19}+6t^{18}+10t^{17}+10t^{16}+10t^{15}+14t^{14}+14t^{13}+14t^{12}+10t^{11}+10t^{10}+10t^9+6t^8+6t^7+6t^6+2t^5+2t^4+2t^3+t^2+t+1\\
\hline
$A_{32}[u;27](t)$ & 
t^{26}+t^{24}+t^{20}+t^{18}+t^8+t^6+t^2+1\\
\hline
$A_{32}[d;27](t)$ & 
t^{26}+t^{25}+t^{24}+t^{23}+t^{22}+t^{21}+t^{20}+t^{19}+t^{18}+t^{17}+t^{16}+t^{15}+2t^{14}+2t^{13}+2t^{12}+t^{11}+t^{10}+t^9+t^8+t^7+t^6+t^5+t^4+t^3+t^2+t+1\\
\hline
$A_{41}[u;27](t)$ & 
t^{26}+2t^{24}+2t^{22}+2t^{20}+t^{18}+t^8+2t^6+2t^4+2t^2+1\\
\hline
$A_{41}[d;27](t)$ & 
t^{26}+2t^{25}+6t^{24}+10t^{23}+14t^{22}+10t^{21}+6t^{20}+2t^{19}+t^{18}+t^{17}+2t^{16}+6t^{15}+10t^{14}+14t^{13}+10t^{12}+6t^{11}+2t^{10}+t^9+t^8+2t^7+6t^6+10t^5+14t^4+10t^3+6t^2+2t+1\\
\hline
$A_{42}[u;27](t)$ & 
t^{26}+t^{24}+t^{20}+t^{18}+t^8+t^6+t^2+1\\
\hline
$A_{42}[d;27](t)$ & 
t^{26}+t^{25}+t^{24}+t^{23}+2t^{22}+t^{21}+t^{20}+t^{19}+t^{18}+t^{17}+t^{16}+t^{15}+t^{14}+2t^{13}+t^{12}+t^{11}+t^{10}+t^9+t^8+t^7+t^6+t^5+2t^4+t^3+t^2+t+1\\
\hline
$A_{521}[u;27](t)$ & 
t^{26}+2t^{24}+2t^{22}+2t^{20}+t^{18}+t^8+2t^6+2t^4+2t^2+1\\
\hline
$A_{521}[d;27](t)$ 
& t^{26}+2t^{25}+6t^{24}+10t^{23}+14t^{22}+10t^{21}+6t^{20}+
2t^{19}+t^{18}+t^{17}+4t^{16}+10t^{15}+20t^{14}+26t^{13}+20t^{12}+
10t^{11}+4t^{10}+t^9+t^8+2t^7+6t^6+10t^5+14t^4+10t^3+6t^2+2t+1\\
\hline
$A_{522}[u;27](t)$ & 
t^{26}+t^{24}+2t^{20}+2t^{18}+2t^{14}+2t^{12}+2t^8+2t^6+t^2+1
\\
\hline
$A_{522}[d;27](t)$ 
& 
t^{26}+t^{25}+t^{24}+2t^{23}+4t^{22}+2t^{21}+6t^{20}+10t^{19}+
6t^{18}+10t^{17}+20t^{16}+10t^{15}+14t^{14}+26t^{13}+14t^{12}+
10t^{11}+20t^{10}+10t^9+6t^8+10t^7+6t^6+2t^5+4t^4+2t^3+t^2+t+1
\\
\hline
$A_{523}[u;27](t)$ & 
t^{26}+t^{24}+t^{20}+t^{18}+t^8+t^6+t^2+1
\\
\hline
$A_{523}[d;27](t)$ & 
t^{26}+t^{25}+t^{24}+t^{23}+2t^{22}+t^{21}+t^{20}+t^{19}+t^{18}+t^{17}+2t^{16}+t^{15}+2t^{14}+4t^{13}+2t^{12}+t^{11}+2t^{10}+t^9+t^8+t^7+t^6+t^5+2t^4+t^3+t^2+t+1
\\
\hline
\end{tabular}
\end{table}

\newpage

\begin{table}[htb!] 
\caption{\label{tbl:genfunctions125} Intermediate generating functions for $n=125$ (Part 1)}
\begin{tabular}{|l|>{\raggedright\arraybackslash$}p{12cm}<{$}|}
\hline
Term & \mbox{Generating function}\\
\hline
\hline
$A_1[u;125](t)$ &
t^{124}+t^{122}+t^{120}+t^{114}+t^{112}+t^{110}+t^{104}+t^{102}+t^{100}+t^{74}+t^{72}+t^{70}+t^{64}+t^{62}+t^{60}+t^{54}+t^{52}+t^{50}+t^{24}+t^{22}+t^{20}+t^{14}+t^{12}+t^{10}+t^{4}+t^{2}+1\\
\hline
$A_1[d;125](t)$ &
t^{124}+t^{123}+2t^{122}+t^{121}+t^{120}+t^{119}+t^{118}+2t^{117}+t^{116}+t^{115}+2t^{114}+2t^{113}+4t^{112}+2t^{111}+2t^{110}+t^{109}+t^{108}+2t^{107}+t^{106}+t^{105}+t^{104}+t^{103}+2t^{102}+t^{101}+t^{100}+t^{99}+t^{98}+2t^{97}+t^{96}+t^{95}+t^{94}+t^{93}+2t^{92}+t^{91}+t^{90}+2t^{89}+2t^{88}+4t^{87}+2t^{86}+2t^{85}+t^{84}+t^{83}+2t^{82}+t^{81}+t^{80}+t^{79}+t^{78}+2t^{77}+t^{76}+t^{75}+2t^{74}+2t^{73}+4t^{72}+2t^{71}+2t^{70}+2t^{69}+2t^{68}+4t^{67}+2t^{66}+2t^{65}+4t^{64}+4t^{63}+8t^{62}+4t^{61}+4t^{60}+2t^{59}+2t^{58}+4t^{57}+2t^{56}+2t^{55}+2t^{54}+2t^{53}+4t^{52}+2t^{51}+2t^{50}+t^{49}+t^{48}+2t^{47}+t^{46}+t^{45}+t^{44}+t^{43}+2t^{42}+t^{41}+t^{40}+2t^{39}+2t^{38}+4t^{37}+2t^{36}+2t^{35}+t^{34}+t^{33}+2t^{32}+t^{31}+t^{30}+t^{29}+t^{28}+2t^{27}+t^{26}+t^{25}+t^{24}+t^{23}+2t^{22}+t^{21}+t^{20}+t^{19}+t^{18}+2t^{17}+t^{16}+t^{15}+2t^{14}+2t^{13}+4t^{12}+2t^{11}+2t^{10}+t^{9}+t^{8}+2t^{7}+t^{6}+t^{5}+t^{4}+t^{3}+2t^{2}+t+1\\
\hline
$A_{31}[u;125](t)$ & 
t^{124}+t^{122}+t^{120}+2t^{114}+2t^{112}+2t^{110}+8t^{104}+8t^{102}+8t^{100}+22t^{94}+22t^{92}+22t^{90}+51t^{84}+51t^{82}+51t^{80}+80t^{74}+80t^{72}+80t^{70}+96t^{64}+96t^{62}+96t^{60}+80t^{54}+80t^{52}+80t^{50}+51t^{44}+51t^{42}+51t^{40}+22t^{34}+22t^{32}+22t^{30}+8t^{24}+8t^{22}+8t^{20}+2t^{14}+2t^{12}+2t^{10}+t^4+t^2+1\\
\hline
$A_{31}[d;125](t)$ & 
t^{124}+t^{123}+2t^{122}+t^{121}+t^{120}+2t^{119}+2t^{118}+4t^{117}+2t^{116}+2t^{115}+16t^{114}+16t^{113}+32t^{112}+16t^{111}+16t^{110}+102t^{109}+102t^{108}+204t^{107}+102t^{106}+102t^{105}+536t^{104}+536t^{103}+1072t^{102}+536t^{101}+536t^{100}+2126t^{99}+2126t^{98}+4252t^{97}+2126t^{96}+2126t^{95}+6744t^{94}+6744t^{93}+13488t^{92}+6744t^{91}+6744t^{90}+17310t^{89}+17310t^{88}+34620t^{87}+17310t^{86}+17310t^{85}+36803t^{84}+36803t^{83}+73606t^{82}+36803t^{81}+36803t^{80}+65376t^{79}+65376t^{78}+130752t^{77}+65376t^{76}+65376t^{75}+98104t^{74}+98104t^{73}+196208t^{72}+98104t^{71}+98104t^{70}+124812t^{69}+124812t^{68}+249624t^{67}+124812t^{66}+124812t^{65}+135264t^{64}+135264t^{63}+270528t^{62}+135264t^{61}+135264t^{60}+124812t^{59}+124812t^{58}+249624t^{57}+124812t^{56}+124812t^{55}+98104t^{54}+98104t^{53}+196208t^{52}+98104t^{51}+98104t^{50}+65376t^{49}+65376t^{48}+130752t^{47}+65376t^{46}+65376t^{45}+36803t^{44}+36803t^{43}+73606t^{42}+36803t^{41}+36803t^{40}+17310t^{39}+17310t^{38}+34620t^{37}+17310t^{36}+17310t^{35}+6744t^{34}+6744t^{33}+13488t^{32}+6744t^{31}+6744t^{30}+2126t^{29}+2126t^{28}+4252t^{27}+2126t^{26}+2126t^{25}+536t^{24}+536t^{23}+1072t^{22}+536t^{21}+536t^{20}+102t^{19}+102t^{18}+204t^{17}+102t^{16}+102t^{15}+16t^{14}+16t^{13}+32t^{12}+16t^{11}+16t^{10}+2t^9+2t^8+4t^7+2t^6+2t^5+t^4+t^3+2t^2+t+1
\\
\hline
$A_{32}[u;125](t)$ & 
t^{124}+t^{122}+t^{120}+t^{114}+t^{112}+t^{110}+t^{104}+t^{102}+t^{100}+t^{74}+t^{72}+t^{70}+2t^{64}+2t^{62}+2t^{60}+t^{54}+t^{52}+t^{50}+t^{24}+t^{22}+t^{20}+t^{14}+t^{12}+t^{10}+t^4+t^2+1\\
\hline
$A_{32}[d;125](t)$ & 
t^{124}+t^{123}+2t^{122}+t^{121}+t^{120}+t^{119}+t^{118}+2t^{117}+ t^{116}+t^{115}+2t^{114}+2t^{113}+4t^{112}+2t^{111}+2t^{110}+ t^{109}+t^{108}+2t^{107}+t^{106}+t^{105}+t^{104}+t^{103}+2t^{102}+ t^{101}+t^{100}+t^{99}+t^{98}+2t^{97}+t^{96}+t^{95}+4t^{94}+4t^{93}+8t^{92}+4t^{91}+4t^{90}+6t^{89}+6t^{88}+12t^{87}+6t^{86}+6t^{85}+ 4t^{84}+4t^{83}+8t^{82}+4t^{81}+4t^{80}+t^{79}+t^{78}+2t^{77}+ t^{76}+t^{75}+2t^{74}+2t^{73}+4t^{72}+2t^{71}+2t^{70}+6t^{69}+ 6t^{68}+12t^{67}+6t^{66}+6t^{65}+10t^{64}+10t^{63}+20t^{62}+ 10t^{61}+10t^{60}+6t^{59}+6t^{58}+12t^{57}+6t^{56}+6t^{55}+2t^{54}+2t^{53}+4t^{52}+2t^{51}+2t^{50}+t^{49}+t^{48}+2t^{47}+t^{46}+t^{45}+4t^{44}+4t^{43}+8t^{42}+4t^{41}+4t^{40}+6t^{39}+6t^{38}+12t^{37}+ 6t^{36}+6t^{35}+4t^{34}+4t^{33}+8t^{32}+4t^{31}+4t^{30}+t^{29}+ t^{28}+2t^{27}+t^{26}+t^{25}+t^{24}+t^{23}+2t^{22}+t^{21}+t^{20}+ t^{19}+t^{18}+2t^{17}+t^{16}+t^{15}+2t^{14}+2t^{13}+4t^{12}+2t^{11}+2t^{10}+t^{9}+t^{8}+2t^{7}+t^{6}+t^{5}+t^{4}+t^{3}+2t^{2}+t+1\\
\hline
\hline
\end{tabular}
\end{table}

\newpage

\begin{table}[htb!] 
\caption{\label{tbl:genfunctions125cont}  Intermediate generating functions for $n=125$ (Part 2)}
\begin{tabular}{|l|>{\raggedright\arraybackslash$}p{12cm}<{$}|}
\hline
Term & \mbox{Generating function}\\
\hline
\hline
$A_{41}[u;125](t)$ & 
t^{124}+2t^{122}+8t^{120}+22t^{118}+51t^{116}+80t^{114}+96t^{112}+80t^{110}+51t^{108}+22t^{106}+8t^{104}+2t^{102}+t^{100}+t^{74}+2t^{72}+8t^{70}+22t^{68}+51t^{66}+80t^{64}+96t^{62}+80t^{60}+51t^{58}+22t^{56}+8t^{54}+2t^{52}+t^{50}+t^{24}+2t^{22}+8t^{20}+22t^{18}+51t^{16}+80t^{14}+96t^{12}+80t^{10}+51t^8+22t^6+8t^4+2t^2+1\\
\hline
$A_{41}[d;125](t)$ & 
t^{124}+2t^{123}+16t^{122}+102t^{121}+536t^{120}+2126t^{119}+6744t^{118}+17310t^{117}+36803t^{116}+65376t^{115}+98104t^{114}+124812t^{113}+135264t^{112}+124812t^{111}+98104t^{110}+65376t^{109}+36803t^{108}+17310t^{107}+6744t^{106}+2126t^{105}+536t^{104}+102t^{103}+16t^{102}+2t^{101}+t^{100}+t^{99}+2t^{98}+16t^{97}+102t^{96}+536t^{95}+2126t^{94}+6744t^{93}+17310t^{92}+36803t^{91}+65376t^{90}+98104t^{89}+124812t^{88}+135264t^{87}+124812t^{86}+98104t^{85}+65376t^{84}+36803t^{83}+17310t^{82}+6744t^{81}+2126t^{80}+536t^{79}+102t^{78}+16t^{77}+2t^{76}+t^{75}+2t^{74}+4t^{73}+32t^{72}+204t^{71}+1072t^{70}+4252t^{69}+13488t^{68}+34620t^{67}+73606t^{66}+130752t^{65}+196208t^{64}+249624t^{63}+270528t^{62}+249624t^{61}+196208t^{60}+130752t^{59}+73606t^{58}+34620t^{57}+13488t^{56}+4252t^{55}+1072t^{54}+204t^{53}+32t^{52}+4t^{51}+2t^{50}+t^{49}+2t^{48}+16t^{47}+102t^{46}+536t^{45}+2126t^{44}+6744t^{43}+17310t^{42}+36803t^{41}+65376t^{40}+98104t^{39}+124812t^{38}+135264t^{37}+124812t^{36}+98104t^{35}+65376t^{34}+36803t^{33}+17310t^{32}+6744t^{31}+2126t^{30}+536t^{29}+102t^{28}+16t^{27}+2t^{26}+t^{25}+t^{24}+2t^{23}+16t^{22}+102t^{21}+536t^{20}+2126t^{19}+6744t^{18}+17310t^{17}+36803t^{16}+65376t^{15}+98104t^{14}+124812t^{13}+135264t^{12}+124812t^{11}+98104t^{10}+65376t^9+36803t^8+17310t^7+6744t^6+2126t^5+536t^4+102t^3+16t^2+2t+1
\\
\hline
$A_{42}[u;125](t)$ & 
t^{124}+t^{122}+t^{120}+t^{114}+2t^{112}+t^{110}+t^{104}+t^{102}+t^{100}+t^{74}+t^{72}+t^{70}+t^{64}+2t^{62}+t^{60}+t^{54}+t^{52}+t^{50}+t^{24}+t^{22}+t^{20}+t^{14}+2t^{12}+t^{10}+t^4+t^2+1\\
\hline
$A_{42}[d;125](t)$ & 
t^{124}+t^{123}+2t^{122}+t^{121}+t^{120}+t^{119}+4t^{118}+6t^{117}+4t^{116}+t^{115}+2t^{114}+6t^{113}+10t^{112}+6t^{111}+2t^{110}+ t^{109}+4t^{108}+6t^{107}+4t^{106}+t^{105}+t^{104}+t^{103}+2t^{102}+t^{101}+t^{100}+t^{99}+t^{98}+2t^{97}+t^{96}+t^{95}+t^{94}+4t^{93}+6t^{92}+4t^{91}+t^{90}+2t^{89}+6t^{88}+10t^{87}+6t^{86}+2t^{85}+
t^{84}+4t^{83}+6t^{82}+4t^{81}+t^{80}+t^{79}+t^{78}+2t^{77}+t^{76}+
t^{75}+2t^{74}+2t^{73}+4t^{72}+2t^{71}+2t^{70}+2t^{69}+8t^{68}+ 12t^{67}+8t^{66}+2t^{65}+4t^{64}+12t^{63}+20t^{62}+12t^{61}+4t^{60}+2t^{59}+8t^{58}+12t^{57}+8t^{56}+2t^{55}+2t^{54}+2t^{53}+4t^{52}+ 2t^{51}+2t^{50}+t^{49}+t^{48}+2t^{47}+t^{46}+t^{45}+t^{44}+4t^{43}+6t^{42}+4t^{41}+t^{40}+2t^{39}+6t^{38}+10t^{37}+6t^{36}+2t^{35}+ t^{34}+4t^{33}+6t^{32}+4t^{31}+t^{30}+t^{29}+t^{28}+2t^{27}+t^{26}+t^{25}+t^{24}+t^{23}+2t^{22}+t^{21}+t^{20}+t^{19}+4t^{18}+6t^{17}+
4t^{16}+t^{15}+2t^{14}+6t^{13}+10t^{12}+6t^{11}+2t^{10}+t^{9}+
4t^{8}+6t^{7}+4t^{6}+t^{5}+t^{4}+t^{3}+2t^{2}+t+1\\
\hline
$A_{521}[u;125](t)$ & 
t^{124}+2t^{122}+8t^{120}+22t^{118}+51t^{116}+80t^{114}+96t^{112}+80t^{110}+51t^{108}+22t^{106}+8t^{104}+2t^{102}+t^{100}+t^{74}+4t^{72}+14t^{70}+44t^{68}+99t^{66}+160t^{64}+188t^{62}+160t^{60}+99t^{58}+44t^{56}+14t^{54}+4t^{52}+t^{50}+t^{24}+2t^{22}+8t^{20}+22t^{18}+51t^{16}+80t^{14}+96t^{12}+80t^{10}+51t^8+22t^6+8t^4+2t^2+1\\
\hline
\hline
\end{tabular}
\end{table}

\newpage

\begin{table}[htb!] 
\caption{\label{tbl:genfunctions125cont2} Intermediate generating functions for $n=125$ (Part 3)}
\begin{tabular}{|l|>{\raggedright\arraybackslash$}p{12cm}<{$}|}
\hline
Term & \mbox{Generating function}\\
\hline
\hline
$A_{521}[d;125](t)$ & 
t^{124}+2t^{123}+16t^{122}+102t^{121}+536t^{120}+2126t^{119}+6744t^{118}+17310t^{117}+36803t^{116}+65376t^{115}+98104t^{114}+124812t^{113}+135264t^{112}+124812t^{111}+98104t^{110}+65376t^{109}+36803t^{108}+17310t^{107}+6744t^{106}+2126t^{105}+536t^{104}+102t^{103}+16t^{102}+2t^{101}+t^{100}+t^{99}+8t^{98}+60t^{97}+408t^{96}+2126t^{95}+8504t^{94}+26932t^{93}+69240t^{92}+147107t^{91}+261504t^{90}+392256t^{89}+499248t^{88}+540860t^{87}+499248t^{86}+392256t^{85}+261504t^{84}+147107t^{83}+69240t^{82}+26932t^{81}+8504t^{80}+2126t^{79}+408t^{78}+60t^{77}+8t^{76}+t^{75}+2t^{74}+12t^{73}+92t^{72}+612t^{71}+3196t^{70}+12756t^{69}+40420t^{68}+103860t^{67}+220710t^{66}+392256t^{65}+588464t^{64}+748872t^{63}+811384t^{62}+748872t^{61}+588464t^{60}+392256t^{59}+220710t^{58}+103860t^{57}+40420t^{56}+12756t^{55}+3196t^{54}+612t^{53}+92t^{52}+12t^{51}+2t^{50}+t^{49}+8t^{48}+60t^{47}+408t^{46}+2126t^{45}+8504t^{44}+26932t^{43}+69240t^{42}+147107t^{41}+261504t^{40}+392256t^{39}+499248t^{38}+540860t^{37}+499248t^{36}+392256t^{35}+261504t^{34}+147107t^{33}+69240t^{32}+26932t^{31}+8504t^{30}+2126t^{29}+408t^{28}+60t^{27}+8t^{26}+t^{25}+t^{24}+2t^{23}+16t^{22}+102t^{21}+536t^{20}+2126t^{19}+6744t^{18}+17310t^{17}+36803t^{16}+65376t^{15}+98104t^{14}+124812t^{13}+135264t^{12}+124812t^{11}+98104t^{10}+65376t^9+36803t^8+17310t^7+6744t^6+2126t^5+536t^4+102t^3+16t^2+2t+1\\
\hline
$A_{522}[u;125](t)$ & 
t^{124}+t^{122}+t^{120}+2t^{114}+4t^{112}+2t^{110}+8t^{104}+14t^{102}+8t^{100}+22t^{94}+44t^{92}+22t^{90}+51t^{84}+99t^{82}+51t^{80}+80t^{74}+160t^{72}+80t^{70}+96t^{64}+188t^{62}+96t^{60}+80t^{54}+160t^{52}+80t^{50}+51t^{44}+99t^{42}+51t^{40}+22t^{34}+44t^{32}+22t^{30}+8t^{24}+14t^{22}+8t^{20}+2t^{14}+4t^{12}+2t^{10}+t^4+t^2+1\\
\hline
$A_{522}[d;125](t)$ & 
t^{124}+t^{123}+2t^{122}+t^{121}+t^{120}+2t^{119}+8t^{118}+12t^{117}+8t^{116}+2t^{115}+16t^{114}+60t^{113}+92t^{112}+60t^{111}+16t^{110}+102t^{109}+408t^{108}+612t^{107}+408t^{106}+102t^{105}+536t^{104}+2126t^{103}+3196t^{102}+2126t^{101}+536t^{100}+2126t^{99}+8504t^{98}+12756t^{97}+8504t^{96}+2126t^{95}+6744t^{94}+26932t^{93}+40420t^{92}+26932t^{91}+6744t^{90}+17310t^{89}+69240t^{88}+103860t^{87}+69240t^{86}+17310t^{85}+36803t^{84}+147107t^{83}+220710t^{82}+147107t^{81}+36803t^{80}+65376t^{79}+261504t^{78}+392256t^{77}+261504t^{76}+65376t^{75}+98104t^{74}+392256t^{73}+588464t^{72}+392256t^{71}+98104t^{70}+124812t^{69}+499248t^{68}+748872t^{67}+499248t^{66}+124812t^{65}+135264t^{64}+540860t^{63}+811384t^{62}+540860t^{61}+135264t^{60}+124812t^{59}+499248t^{58}+748872t^{57}+499248t^{56}+124812t^{55}+98104t^{54}+392256t^{53}+588464t^{52}+392256t^{51}+98104t^{50}+65376t^{49}+261504t^{48}+392256t^{47}+261504t^{46}+65376t^{45}+36803t^{44}
+147107t^{43}+220710t^{42}+147107t^{41}+36803t^{40}+17310t^{39}+69240t^{38}+103860t^{37}+69240t^{36}+17310t^{35}+6744t^{34}+26932t^{33}+40420t^{32}+26932t^{31}+6744t^{30}+2126t^{29}+8504t^{28}+12756t^{27}+8504t^{26}+2126t^{25}+536t^{24}+2126t^{23}+3196t^{22}+2126t^{21}+536t^{20}+102t^{19}+408t^{18}+612t^{17}+408t^{16}+102t^{15}+16t^{14}+60t^{13}+92t^{12}+60t^{11}+16t^{10}+2t^9+8t^8+12t^7+8t^6+2t^5+t^4+t^3+2t^2+t+1
\\
\hline
\hline
\end{tabular}
\end{table}

\newpage

\begin{table}[htb!] 
\caption{\label{tbl:genfunctions125cont3} Intermediate generating functions for $n=125$ (Part 4)}
\begin{tabular}{|l|>{\raggedright\arraybackslash$}p{12cm}<{$}|}
\hline
Term & \mbox{Generating function}\\
\hline
\hline
$A_{523}[u;125](t)$ & 
t^{124}+t^{122}+t^{120}+t^{114}+2t^{112}+t^{110}+t^{104}+t^{102}+t^{100}+t^{74}+2t^{72}+t^{70}+2t^{64}+4t^{62}+2t^{60}+t^{54}+2t^{52}+t^{50}+t^{24}+t^{22}+t^{20}+t^{14}+2t^{12}+t^{10}+t^4+t^2+1
\\
\hline
$A_{523}[d;125](t)$ & 
t^{124}+t^{123}+2t^{122}+t^{121}+t^{120}+t^{119}+4t^{118}+6t^{117}+4t^{116}+t^{115}+2t^{114}+6t^{113}+10t^{112}+6t^{111}+2t^{110}+
t^{109}+4t^{108}+6t^{107}+4t^{106}+t^{105}+t^{104}+t^{103}+2t^{102}+t^{101}+t^{100}+t^{99}+4t^{98}+6t^{97}+4t^{96}+t^{95}+4t^{94}+
16t^{93}+24t^{92}+16t^{91}+4t^{90}+6t^{89}+24t^{88}+36t^{87}+
24t^{86}+6t^{85}+4t^{84}+16t^{83}+24t^{82}+16t^{81}+4t^{80}+t^{79}+4t^{78}+6t^{77}+4t^{76}+t^{75}+2t^{74}+6t^{73}+10t^{72}+6t^{71}+ 2t^{70}+6t^{69}+24t^{68}+36t^{67}+24t^{66}+6t^{65}+10t^{64}+ 36t^{63}+56t^{62}+36t^{61}+10t^{60}+6t^{59}+24t^{58}+36t^{57}+ 24t^{56}+6t^{55}+2t^{54}+6t^{53}+10t^{52}+6t^{51}+2t^{50}+t^{49}+ 4t^{48}+6t^{47}+4t^{46}+t^{45}+4t^{44}+16t^{43}+24t^{42}+16t^{41}+ 4t^{40}+6t^{39}+24t^{38}+36t^{37}+24t^{36}+6t^{35}+4t^{34}+ 16t^{33}+24t^{32}+16t^{31}+4t^{30}+t^{29}+4t^{28}+6t^{27}+4t^{26}+ t^{25}+t^{24}+t^{23}+2t^{22}+t^{21}+t^{20}+t^{19}+4t^{18}+6t^{17}+ 4t^{16}+t^{15}+2t^{14}+6t^{13}+10t^{12}+6t^{11}+2t^{10}+t^{9}+
4t^{8}+6t^{7}+4t^{6}+t^{5}+t^{4}+t^{3}+2t^{2}+t+1\\
\hline
\hline
\end{tabular}
\end{table}

\newpage

\begin{table}[htb!] \small
\caption{\label{tbl:alldir27} The generating functions for directed circulant graphs of order $27$} 
\hspace*{-3ex}
\begin{tabular}{|l||>{\raggedright\arraybackslash$}p{11cm}<{$}|}
\hline
$A_1$& t^{26}+t^{25}+t^{24}+t^{23}+t^{22}+t^{21}+t^{20}+t^{19}+t^{18}+t^{17}+t^{16}+
t^{15}+t^{14}+t^{13}+t^{12}+t^{11}+t^{10}+t^9+t^8+t^7+t^6+t^5+t^4+t^3+t^2+t+1\\
\hline
\hline
$A_{21}$& t^{26}+3t^{25}+23t^{24}+152t^{23}+850t^{22}+3680t^{21}+12850t^{20}+36606t^{19}+86919t^{18}+173701t^{17}+295311t^{16}+429388t^{15}+
536810t^{14}+577996t^{13}+536810t^{12}+429388t^{11}+295311t^{10}+173701t^9+86919t^8+36606t^7+12850t^6+3680t^5+850t^4
+152t^3+23t^2+3t+1\\
\hline
$A_{22}$& t^{26}+2t^{25}+6t^{24}+11t^{23}+22t^{22}+38t^{21}+65t^{20}+92t^{19}+129t^{18}+172t^{17}+214t^{16}+235t^{15}+
263t^{14}+276t^{13}+263t^{12}+235t^{11}+214t^{10}+172t^9+129t^8+92t^7+65t^6+38t^5+22t^4+11t^3+6t^2+2t+1\\
\hline
$A_2=A_{21}-A_{22}$& t^{25}+17t^{24}+141t^{23}+828t^{22}+3642t^{21}+12785t^{20}+36514t^{19}+86790t^{18}+173529t^{17}+295097t^{16}+429153t^{15}+
536547t^{14}+577720t^{13}+536547t^{12}+429153t^{11}+295097t^{10}+173529t^9+86790t^8+36514t^7+12785t^6+3642t^5+828t^4+141t^3+
17t^2+t\\
\hline
\hline
$A_{31}$& t^{26}+t^{25}+t^{24}+2t^{23}+2t^{22}+2t^{21}+6t^{20}+6t^{19}+6t^{18}+10t^{17}+10t^{16}+10t^{15}+14t^{14}+14t^{13}+14t^{12}+10t^{11}+
10t^{10}+10t^9+6t^8+6t^7+6t^6+2t^5+2t^4+2t^3+t^2+t+1\\
\hline
$A_{32}$& t^{26}+t^{25}+t^{24}+t^{23}+t^{22}+t^{21}+t^{20}+t^{19}+t^{18}+t^{17}+t^{16}+t^{15}+2t^{14}+2t^{13}+2t^{12}+t^{11}+t^{10}+t^9+t^8+t^7+t^6+t^5+t^4+
t^3+t^2+t+1\\
\hline
$A_3=A_{31}-A_{32}$& t^{23}+t^{22}+t^{21}+5t^{20}+5t^{19}+5t^{18}+9t^{17}+9t^{16}+9t^{15}+12t^{14}+12t^{13}+12t^{12}+9t^{11}+9t^{10}+9t^9+5t^8+5t^7+
5t^6+t^5+t^4+t^3\\
\hline
\hline
$A_{41}$& t^{26}+2t^{25}+6t^{24}+10t^{23}+14t^{22}+10t^{21}+6t^{20}+2t^{19}+t^{18}+t^{17}+2t^{16}+6t^{15}+10t^{14}+14t^{13}+10t^{12}+
6t^{11}+2t^{10}+t^9+t^8+2t^7+6t^6+10t^5+14t^4+10t^3+6t^2+2t+1\\
\hline
$A_{42}$& t^{26}+t^{25}+t^{24}+t^{23}+2t^{22}+t^{21}+t^{20}+t^{19}+t^{18}+t^{17}+t^{16}+t^{15}+t^{14}+2t^{13}+t^{12}+t^{11}+t^{10}+t^9+t^8+t^7+t^6+t^5+
2t^4+t^3+t^2+t+1\\
\hline
$A_4=A_{41}-A_{42}$& t^{25}+5t^{24}+9t^{23}+12t^{22}+9t^{21}+5t^{20}+t^{19}+t^{16}+5t^{15}+9t^{14}+12t^{13}+9t^{12}+5t^{11}+t^{10}+t^7+5t^6+9t^5+
12t^4+9t^3+5t^2+t\\
\hline
\hline
$A_{51}$& t^{26}+2t^{25}+6t^{24}+11t^{23}+18t^{22}+20t^{21}+29t^{20}+38t^{19}+47t^{18}+64t^{17}+86t^{16}+91t^{15}+109t^{14}+124t^{13}+
109t^{12}+91t^{11}+86t^{10}+64t^9+47t^8+38t^7+29t^6+20t^5+18t^4+11t^3+6t^2+2t+1\\
\hline
$A_{521}$& t^{26}+2t^{25}+6t^{24}+10t^{23}+14t^{22}+10t^{21}+6t^{20}+2t^{19}+t^{18}+t^{17}+4t^{16}+10t^{15}+20t^{14}+26t^{13}+20t^{12}+
10t^{11}+4t^{10}+t^9+t^8+2t^7+6t^6+10t^5+14t^4+10t^3+6t^2+2t+1\\
\hline
$A_{522}$& t^{26}+t^{25}+t^{24}+2t^{23}+4t^{22}+2t^{21}+6t^{20}+10t^{19}+6t^{18}+10t^{17}+20t^{16}+10t^{15}+14t^{14}+26t^{13}+14t^{12}+
10t^{11}+20t^{10}+10t^9+6t^8+10t^7+6t^6+2t^5+4t^4+2t^3+t^2+t+1\\
\hline
$A_{523}$& t^{26}+t^{25}+t^{24}+t^{23}+2t^{22}+t^{21}+t^{20}+t^{19}+t^{18}+t^{17}+2t^{16}+t^{15}+2t^{14}+4t^{13}+2t^{12}+t^{11}+2t^{10}+t^9+t^8+
t^7+t^6+t^5+2t^4+t^3+t^2+t+1\\
\hline
{\small $A_5\!=\!A_{51}\!-\!A_{521}\!$} & \\
{\small $-\!A_{522}\!+\!A_{523}$}& 2t^{22}+9t^{21}+18t^{20}+27t^{19}+41t^{18}+54t^{17}+64t^{16}+72t^{15}+77t^{14}+76t^{13}+77t^{12}+72t^{11}+64t^{10}+54t^9+41t^8+
27t^7+18t^6+9t^5+2t^4\\
\hline
\hline
{\small $A\!=\!A_1\!+\!A_2\!\!$} & \\
{\small $+\!A_3\!+\!A_4\!+\!A_5$}& t^{26}+3t^{25}+23t^{24}+152t^{23}+844t^{22}+3662t^{21}+12814t^{20}+36548t^{19}+86837t^{18}+173593t^{17}+295172t^{16}+
429240t^{15}+536646t^{14}+577821t^{13}+536646t^{12}+429240t^{11}+295172t^{10}+173593t^9+86837t^8+36548t^7+12814t^6+3662t^5+
844t^4+152t^3+23t^2+3t+1\\
\hline
\end{tabular}
\end{table}

\end{document}